\documentclass{amsart}

\usepackage[left=1.5in, right=1.5in]{geometry}
\usepackage[latin1]{inputenc}
\usepackage{amsfonts}
\usepackage[toc,page,title,titletoc,header]{appendix}
\usepackage{graphicx,psfrag,epsfig,multirow,caption,subcaption}
\usepackage{amssymb,amsmath,amscd,amsthm,amssymb,verbatim,setspace}
\numberwithin{equation}{section}
\usepackage{mathrsfs,wrapfig}
\usepackage{indentfirst}
\usepackage{extarrows}
\usepackage[colorlinks=true]{hyperref}

\usepackage{marginnote}
\usepackage{color}
\numberwithin{equation}{section}
\newtheorem{Thm}{Theorem}[section]
\newtheorem{Lm}[Thm]{Lemma}
\newtheorem{Prop}[Thm]{Proposition}

\newtheorem{Def}[Thm]{Definition}
\newtheorem{Not}[Thm]{Notation}

\newcommand{\Id}{\mathrm{Id}}
\newcommand{\A}{\mathbb{A}}

\newcommand{\cS}{\mathcal{S}}
\newcommand{\cR}{\mathcal{R}}

\newcommand{\cM}{\mathcal{M}}

\newcommand{\cD}{\mathcal{D}}
\newcommand{\cL}{\mathcal{L}}

\newcommand{\R}{\mathbb{R}}

\newcommand{\cC}{\mathcal{C}}
\newcommand{\N}{\mathbb{N}}

\newcommand{\bu}{\mathbf{u}}
\newcommand{\bl}{\mathbf{l}}

\newcommand{\dt}{\delta}

\newcommand{\eps}{\epsilon}

\def\sx{\mathsf{x}}

\def\bx{\mathbf{x}}
\def\bp{\mathbf{p}}
\def\bs{\mathbf{s}}
\def\bw{\mathbf{w}}
\def\be{\mathbf{e}}

\title[NCS II]{A New Mechanism for Noncollision Singularities}
\begin{document}
\author{Joseph Gerver, Guan Huang, Jinxin Xue}
\email{gerver@camden.rutgers.edu}
\address{University of Rutgers-Camden.}

\email{huangguan@tsinghua.edu.cn}
\address{Yau Mathematical Sciences Center, Jinchunyuan West Building 304, Tsinghua University, Beijing, China, 100084}
\email{jxue@tsinghua.edu.cn}
\address{Yau Mathematical Sciences Center \& Department of Mathematics, Jingzhai 310, Tsinghua University, Beijing, China, 100084}

\maketitle
\begin{abstract}
In this paper, we prove the existence of noncollision singularities in the planar four-body problem with a model different from \cite{X}. In this model, the acceleration can be arbitrarily fast and the masses can be comparable. This work provides a general principle to construct noncollision singularities as well as other related orbits with complicited dynamics. It  not only answers a question in \cite{X} but also solves an analogous version of a conjecture of Anosov.
\end{abstract}

%\tableofcontents
%\renewcommand\contentsname{Index}

\section{Introduction}

In this paper, we prove the existence of noncollision singularities in a model of the four-body problem  which is drastically different from the model of \cite{X}. Noncollision singularities are singularities of the $N$-body problem for which no collision occurs. Their dynamical behaviors are very wild. Indeed, when approaching the singular time, the orbit determined by a noncollision singualrity has to oscillate infinitely often between smaller and smaller neighborhoods of collision and infinity. It had been a longstanding conjecture that noncollision singularities exist in the Newtonian $N$-body problem for $N>3$. The conjecture was finally answered positively in \cite{X} after a century \cite{MM, G0, Xi} etc. In this paper, we use a new acceleration mechanism, which is much faster than the one in \cite{X} and more importantly  gives a very general principle to construct noncollision singularities and other orbits with complicated dynamical behaviors. The new model in this paper does not rely on small perturbations of Kepler motions so it allows comparable masses, including in particular four equal masses.

\subsection{Statements of results}
The configuration is as follows: there are two bodies $Q_1$ and $Q_2$ moving opposite to each other nearly along the $x$-axis and a binary $Q_3$-$Q_4$ with masses $m_3=m_4$ performing nearly Kepler elliptic motion when it is away from  $Q_1$ and $Q_2$. The mass center of the binary moves nearly along the $x$-axis. When the binary gets close to either $Q_1$ or $Q_2$, a triple is formed and the triple will pass close to a triple collision after which the binary is ejected and it will perform nearly Kepler elliptic motion again when they are away from the two exterior bodies. Passing close to a triple collision allows us to greatly accelerate the binary.
\begin{figure}[ht]
\begin{center}
\includegraphics[width=0.6\textwidth]{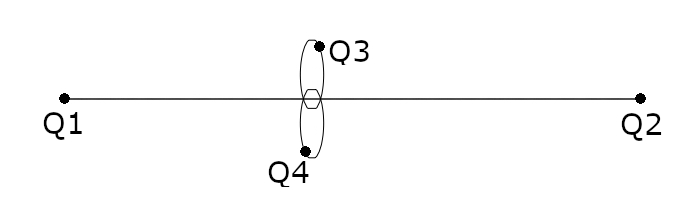}
\end{center}
\caption{The configuration of the four-body problem}
\end{figure}

The main mechanism of acceleration is from the near triple collision dynamics. Near collision, the potential energy goes to $-\infty$ so the kinetic energy becomes large. There is a case of triple collision analyzed by Devaney called the isosceles three-body problem \cite{D}. This is a configuration with two equal masses having the same distance to the third mass. Assume the relative position of the two equal masses is parallel to the $y$-axis and the third mass lies on the $x$-axis.  Devaney classified two types of final motion after approaching near triple collision: either the two equal masses move apart along the $y$-axis while the third body oscillates on the $x$-axis, or the third mass moves along the $x$-axis towards infinity with arbitrarily large velocity and the binary moves opposite to the third body, experiencing repeated binary collisions.  In our model, the near-triple-collision orbits  shadows the latter type.

Using near triple collisions as the acceleration mechanism was pioneered by \cite{MM} and later used in \cite{Xi}. Let us briefly introduce the model of \cite{MM} for comparison with our work. In \cite{MM}, the authors considered a collinear four-body problem labelled from left to right by $Q_1,Q_2,Q_3,Q_4$. The binary $Q_3$-$Q_4$ can be considered a degenerate Kepler elliptic problem while $Q_2$ travels back and forth between $Q_1$ and the binary. Double collisions cannot be avoided. However, it is known that double collisions can be regularized (just imagine replacing a Kepler double collision by an ellastic collision). Allowing countably many double collisions, initial conditions exist such that all four bodies escape to infinity within finite time. This result was the first one that made people believe noncollision singularities exist in the $N$-body problem. However, the main defect is the double collisions. To avoid double collisions, we have to break the symmetry imposed by the collinearity.

Our first result can be considered as an analogue of  \cite{MM}. We consider the isosceles four-body problem by requiring $|Q_i-Q_j|=|Q_{7-i}-Q_j|$, $i=3,4,j=1,2$, where $Q_1$, $Q_2$ and the mass center of $Q_3$-$Q_4$ always have horizontal velocities. In this case, double collision between $Q_3$ and $Q_4$ is not avoidable, but it is regularized and allowed as in \cite{MM}.
\begin{Thm}\label{ThmMM}
 There is an open subset  $\mathscr M$ of $(b_1,\infty)^2\subset \R_+^2$ with $b_1\simeq 0.379$ such that for the isosceles four-body problem with masses $(m_1,m_2)$ in the open set $\mathscr M$ and $m_3=m_4=1$, there exists a nonempty set $\Sigma$ of initial conditions such that for all $x\in \Sigma$, there exists $t_x<\infty$ such that the orbit starting from $x$ satisfies
$$\limsup_{t\to t_x} |Q_i(t)|\to \infty,\quad i=1,2,3,4, $$
and there is no triple collision on the time interval $[0,t_x).$ Moreover,  there exists $x\in \Sigma$ such that $\frac{|\dot{Q}_3(t_{n+1})+\dot{Q}_4(t_{n+1})|}{|\dot{Q}_3(t_n)+\dot{Q}_4(t_n)|}\to\infty$,
where $t_n\to t_x$ is a sequence of times when the binary visits a fixed neighborhood of the origin for the $n$-th time. %Moreover, for any sequence of positive numbers $\{a_n\to\infty\}$, there exists $x\in\Sigma$ and a sequence $0<t_1<t_1'<t_2<t_2'<\cdots<t_x$ such that for the  orbit starting from $x$, we have
%$$\frac{|\dot{Q}_3(t_n')+\dot{Q}_4(t_n')|}{|\dot{Q}_3(t_n)+\dot{Q}_4(t_n)|}>a_n,$$
%where $t_n<t_n'$ are two consecutive moments when the mass center of the pair $Q_3$-$Q_4$ passes through the origin.
%Moreover, for any sequence of numbers $\{a_n\}\to \infty$, we have $\frac{1}{2}|\dot Q_3(t_n)+\dot Q_4(t_n)|\geq a_n$ where $t_n(<t_x)$ is the $n$-th time when the mass center of $Q_3$ and $Q_4$ passes through the origin.
\end{Thm}

 %The  proof of the above theorem combines the triple collision blow up techniques from \cite{D} as well as the derivative control and regularization techniques in \cite{X}, which is very simple in the current setting.  The main difficulty is to control the timing so that when the binary gets close to either $Q_1$ and $Q_2$, it is in a correct position so that when ejected from the near triple collision the velocity of the binary gets significantly increased. %This is handled by the fact that the relative phase of the binary depends sensitively on the initial conditions.

 After the work of \cite{MM}, it was conjectured by Anosov \cite{A}  that a genuine noncollision singularity can be found in a neighborhood of the solution in \cite{MM}.  The conjecture of Anosov as stated is very vague, for instance, the term ``neighborhood" needs a formulation. However, essentially the conjecture of Anosov requires a perturbative analysis of the solution of \cite{MM}, in parituclar, the study of the collinearity symmetry broken. In the next theorem, we show that perturbing the isosceles configuration in Theorem \ref{ThmMM} slightly, a genuine noncollision singularity can be constructed.  In dynamical systems, when a special solution is discovered, it is natural and important to study the stability of the solution. The {\it survival of the stable} principle implies that only stable solutions are the main concern. %Perturbing slightly away from the above solution, we shall allow the binary to gain certain nonzero but very tiny angular of momentum. In this way, we avoid double collision. We obtain the following result on the existence of noncollision singularities in the current model.

%In this paper we solve two conjectures.

%The first conjecture is the following conjecture of Saari and Xia.

%Moreover, since the masses are comparable, which is different from the model of \cite{G,X} where there are arbitrarily small masses.
\begin{Thm}\label{ThmMain}
There is a residual subset $\mathscr R$ of $\mathscr M$ in Theorem \ref{ThmMM} such that for each $(m_1,m_2)\in \mathscr R$ and $m_3=m_4=1$,  there exists a nonempty set $\Sigma$ of initial conditions such that for all $x\in \Sigma$, there exists $t_x<\infty$ such that the orbit starting from $x$ satisfies
$$\limsup_{t\to t_x} |Q_i(t)|\to \infty,\quad i=1,2,3,4, $$
and there is no collision on the time interval $[0,t_x). $  Moreover,  there exists $x\in \Sigma$ such that $\frac{|\dot{Q}_3(t_{n+1})+\dot{Q}_4(t_{n+1})|}{|\dot{Q}_3(t_n)+\dot{Q}_4(t_n)|}\to\infty$,
where $t_n\to t_x$ is a sequence of times when the binary visits a fixed neighborhood of the origin for the $n$-th time.%Moreover, for any sequence of positive  numbers $\{a_n\to\infty\}$, there exists $x\in\Sigma$ and a sequence $0<t_1<t_1'<t_2<t_2'<\cdots<t_x$ such that for the orbit starting from $x$, we have
%$$\frac{|\dot{Q}_3(t_n')+\dot{Q}_4(t_n')|}{|\dot{Q}_3(t_n)+\dot{Q}_4(t_n)|}>a_n,$$
%where $t_n<t_n'$ are two consecutive moments when the mass center of the pair $Q_3$-$Q_4$ passes through the origin.
%Moreover, for any sequence of numbers $\{a_n\}\to \infty$, we have $\frac{1}{2}|\dot Q_3(t_n)+\dot Q_4(t_n)|\geq a_n$ where $t_n(<t_x)$ is the $n$-th time when the mass center of $Q_3$ and $Q_4$ passes through the origin.
%

\end{Thm}
The essential  difficulty in the proof is to control the dynamics after breaking the symmetry of Theorem \ref{ThmMM}, which is a problem that has never appeared in the literature before. %Perturbing away from the isosceles case breaks the symmetry hence significantly increases the number of degrees of freedom. The crucial ingredient in the proof is  certain transversality condition (Propositions \ref{PropDG2}, \ref{LmTrans} and \ref{PropTrans1}) which controls various angular momenta and arguments of periapsis so that all the four bodies move closer and closer to a straight line and the binary gets more and more perpendicular to the line. This fact is revealed by heavy derivative calculations along an orbit from Theorem~\ref{ThmMM}. 
If we consider Theorem~\ref{ThmMM} as an analogue of the result of \cite{MM}, considering the commonality of using the near triple collision as the mechanism of acceleration, then Theorem \ref{ThmMain} is  an affirmative answer to the analogue of the conjecture of Anosov.

\subsection{Backgrounds and perspectives}
The existence of noncollision singularities for the Newtonian $N$-body problem with $N\geq 4$ has been known as the Painlev\'e conjecture since 1897. After the work of von Zeipel, we know that for a noncollision singularity to occur, when the time approaches the singular time, we have 
$$\lim_t\max_{i\neq j}|Q_i(t)-Q_j(t)|=\infty,\quad \lim_t\min_{i\neq j}|Q_i(t)-Q_j(t)|=0.$$
The result of von Zeipel reveals the remarkable nonlocal and oscillating nature of noncollision singularities. Since the orbit has to approach collision repeatedly,  to construct a noncollision singularity, we need a local model for the dynamics of near collision and a way to glue local pieces together. The models \cite{MM,Xi,G0} impose strong symmetries on the model, which almost trivializes the latter problem of gluing different local pieces. However, the symmetries also make the models sensitive to perturbations, which leaves open the conjecture of Anosov. %The local model concerning the dynamics near collision always gives strong hyperbolicity. To glue different local pieces, it is crucial to establish the transversality of the stable and unstable manifolds of the local models. %The hyperbolicity-transversality framework in this paper is a development of the works  \cite{Do,DX,X} etc. We believe this is the key ingredient in establishing the analogue of Anosov conjecture. 

It is of central  importance to find a local model on the dynamics near collision in order to construct a noncollision singularity, since near collision means intensive energy exchange and an acceleration mechanism. The local model of \cite{X} is an ideal model constructed by the first named author \cite{G}, which involves only Kepler motions and achieves self-similarity after a two-step procedure. This is a rather mysterious fact of the two-body problem which deserves more attentions. It seems to us that there is a lack of a general principle to reproduce more such models. However, the local model that we use in this paper has a general feature. Indeed, when a collision occurs, it is standard to blow up the spacetime to study the dynamics approaching a collision. The blowup procedure in general yields homothetic solutions of the $N$-body problem modeled on central configurations. In other words, we can model the local dynamics near collision on central configurations. % In this paper, we use the Lagrange central configuration in three-body problem that is a configuration of an equilateral triangle formed by the three bodies. It turns out that the Lagrange central configuration is a hyperbolic fixed point of the blowup dynamics, which has stable and unstabla manifolds. Then the proof of Theorem \ref{ThmMain} amounts to establish transversality of these invariant manifolds from one local model to the next.
  From this perspective, the noncollision singularities constructed in this paper are rather different from that of \cite{X} since the local model of \cite{X} is nowhere close to a central configuration. 

The result also gives a general principle to construct complicated dynamics in the $N$-body problem by studying central configurations and the transversal intersection of their stable and unstable manifolds. In the $N$-body problem, there are $n\choose k$ choices of $k$-collisions, $k\geq 3$. Blowing up these $k$-collisions, we will arrive at fixed points corresponding to central configurations with various indices, and there are heteroclinic orbits connecting saddle type fixed points arising in this way, and there are also connections between different saddles on the same collision manifold. In general, we get rich dynamics by shadowing these heteroclinic orbits (c.f. \cite{Mo1}-\cite{Mo4},\cite{MoMo}) and the methods in this paper are expected to be applicable.

Compared to \cite{X}, there are some further new features in this model. First, the model in \cite{X} has to have two very small masses, since the proof in \cite{X} is essentially perturbative. In~\cite{X}, the author asked if a four-body problem with comparable masses can admit noncollision singularities, which is answered by  Theorem~\ref{ThmMain}. Second, the acceleration rate of the noncollision singularities during each return in \cite{X} is roughly constant, however, in Theorem~\ref{ThmMain} (the ``moreover" part), the acceleration rate can be arbitrarily large.

%Finally, we have the following speculations. The noncollision singularities are constructed by perturbing heteroclinic orbits between the hyperbolic fixed points in the blowups of two different triple collisions. Most of the work is devoted to showing certain transversality, more precisely,  the stable manifold of one fixed point intersects transversally the unstable manifold of the other.
  
%In this section, we give  several technical lemmas and complete the proof of the main theorem. %The main coordinates that we use are Delaunay for the shuttling between $Q_1$ and $Q_2$ and the blowup coordinates near triple collisions. The two systems are related to each other through polar coordinates and Cartesian coordinates.

The paper is organized as follows. In Section \ref{SCoord}, we introduce various coordinates that we use later to describe the dynamics. In Section \ref{SMM}, we study the isosceles four-body problem and prove Theorem \ref{ThmMM}. In Section \ref{SMainProof}, we give the proof of Theorem \ref{ThmMain} assuming certain propositions to be proved in later sections. In Section \ref{SJacobi}, we write down the Hamiltonian in various coordinates and study the isosceles three-body problem (I3BP). In Section \ref{SBlowupMM}, we provide some preliminary estimates of the  isosceles four-body problem (I4BP). In Section \ref{SDerLl}, we study the $C^1$ dynamics of the I4BP. In Section \ref{local-gs}, we give the $C^1$ estimate for the local map in the F4BP. In Section \ref{SG0g0}, we study the $C^1$ estimate of the global map in the F4BP.  Finally, we have seven appendices. In Appendix \ref{app-blowup}, we write down the equations of motion of the 3BP in the blowup coordinates. In Appendix \ref{App-delaunay}, we give the Delaunay coordinates and their relation to the Cartesian coordinates. In Appendix \ref{SSDG}, we compute the derivative of the global map in the I4BP. In Appendix \ref{App-2-derviatives}, we compute the second order derivatives of the potentials with respect to the Delaunay variables  $G$s and $g$s. In Appendix \ref{appendix-D}, we compute the derivatives of  the  Delaunay variables in the left-right transition.  In Appendix \ref{STrans}, we present a numerical verification of  the transversality. 

\section{The coordinates}\label{SCoord}
In this section we introduce several systems of coordinates for the four-body problem that we will use later.
Without loss of generality, we assume $m_3=m_4=1$ throughout the paper and we always assume $Q_1$ is on the left and $Q_2$ is on the right.
\subsection{The Jacobi-Cartesian coordinates}\label{subsection-jacobi}
The first step is to remove the translation invariance. We introduce $q_i=Q_i-Q_4,\ i=3,1,2$. So we get that the following symplectic form is preserved by assuming $\sum_{i=1}^4 P_i=0$
$$\sum_{i=1}^4 dP_i\wedge dQ_i=\sum_{i=1}^3 dP_i\wedge dq_i.$$
We next assume that the pair is closer to $Q_1$ than $Q_2$ and introduce the following, which we call the {\it left Jacobi coordinates}
\begin{equation}\label{EqJacobi}
\begin{cases}
\bx_0&=q_3\\
\bx_1&=q_1-\frac{q_3}{2}\\
\bx_2&=q_2-\frac{m_1 q_1+ q_3 }{m_1+2}
\end{cases},\quad
\begin{cases}
\bp_0&=P_3+\frac{ P_1}{2}+\frac{ P_2}{2}\\
\bp_1&=P_1+\frac{ m_1 P_2}{m_1+2}\\
\bp_2&=P_2
\end{cases}.
\end{equation}
It can be easily checked that we have the following reduced symplectic form preserved by the coordinate change
$$\sum_{i=3,1,2} dP_i\wedge dq_i=\sum_{i=0,1,2}d\bp_i\wedge d\bx_i.$$
When the pair is closer to $Q_2$ than $Q_1$, we introduce the {\it right Jacobi coordinates} by permuting the subscripts $1$ and $2$.

\subsection{The polar coordinates}\label{def-jacobi-polar}
We next  introduce polar coordinates for each pair $(\bx_i,\bp_i),\ i=0,1,2$.  The polar coordinates are given by the following relation
\begin{equation}
\begin{aligned}
\begin{cases}
\bx_i&=(r_i\cos\theta_i,r_i\sin\theta_i),\\
\bp_i&=(R_i\cos\theta_i-\frac{\Theta_i}{r_i}\sin\theta_i,R_i\sin\theta_i+\frac{\Theta_i}{r_i}\cos\theta_i).\\
\end{cases}
\end{aligned}
\end{equation}
It is clear that $r_i$ has the geometric meaning of length and $\theta_i$ is the polar angle, $R_i$ is the radial momentum conjugate to $r_i$ and $\Theta_i$ is the angular momentum.  It can be verified that the following symplectic form is preserved: $$\sum_{i=0,1,2}d\bp_i\wedge d\bx_i=\sum_{i=0,1,2}(dR_i\wedge dr_i+d\Theta_i\wedge d\theta_i).$$
We use the following list of polar variables to describe the dynamics
$$(R_0,r_0,\Theta_0,\theta_0; R_1,r_1,\Theta_1, \theta_1; R_2,r_2,\Theta_2, \theta_2).$$
\subsection{The Delaunay coordinates}
We  introduce Delaunay coordinates $(L_i,\ell_i,G_i,g_i)$ for each pair $(\bx_i,\bp_i)$ (see Appendix \ref{App-delaunay}). Delaunay coordinates are symplectic, since we have $d\bp_i\wedge d\bx_i=dL_i\wedge d\ell_i+dG_i\wedge dg_i.$
The Delaunay variables have explicit geometric meanings. For Kepler elliptic motion $L_i^2$ is the semimajor axis, $L_iG_i$ is the semiminor axis, $G_i=\bx_i\times \bp_i$ is the angular momentum hence is the same as $\Theta_i$, and $g_i$ is the argument of periapsis. We use the following list of twelve Delaunay variables to describe the dynamics
$$(L_0,\ell_0,G_0,g_0;L_1,\ell_1,G_1,g_1; L_2,\ell_2,G_2,  g_2).$$
%We shall always assume that the total angular momentum vanishes $G_0+G_1+G_2=0$. Note that the Hamiltonian system depends only on the relative angles $g_i-g_j,\ i,j=0,1,2,$ but does not depend on each individual angle $g_i,\ i=0,1,2.$ Then we can further simplify our set of variables as $(L_0,\ell_0,L_1,\ell_1,,L_2,\ell_2, G_0,g_0-g_1, G_2,  g_2-g_1)$ for the left case and $(L_0,\ell_0,L_1,\ell_1,,L_2,\ell_2, G_0,g_0-g_2, G_1,  g_1-g_2)$ for the right case. 

\subsection{The triple collision blowup coordinates}\label{SSSBlowupCoord}
When $Q_1$-$Q_3$-$Q_4$ or $Q_2$-$Q_3$-$Q_4$ is close to a triple collision, we shall use a blowup coordinate system to study the dynamics. We give the triple collision blowup in Appendix \ref{app-blowup} and will work on the isosceles case in the next section.
%When necessary, we  use the upper indexes $L$ and $R$ to indicate  the blow-up coordinates for the triples $Q_1$-$Q_3$-$Q_4$ and $Q_2$-$Q_3$-$Q_3$, respectively.
\section{The isosceles four-body probem and proof of Theorem \ref{ThmMM}}\label{SMM}
In this section we consider the isosceles four-body problem (I4BP) and prove Theorem \ref{ThmMM}. We assume that $Q_1$ and $Q_2$ move on the $x$-axis off to infinity and $Q_3$-$Q_4$ stays on collision-ejection orbit perpendicular to the $x$-axis. Thus we have $G_0=g_0-\pi/2=G_1=g_1=G_2=g_2=0$ and the system has three degrees of freedom. %This system can be completely described by the six variables
%$$(r,v,\psi,w,L_2,\ell_2)\quad \text{or}\quad(L_0,\ell_0,L_1,\ell_1,L_2,\ell_2).$$

\subsection{The isosceles three-body problem (I3BP)}
\begin{figure}[ht]
\begin{center}
\includegraphics[width=0.8\textwidth]{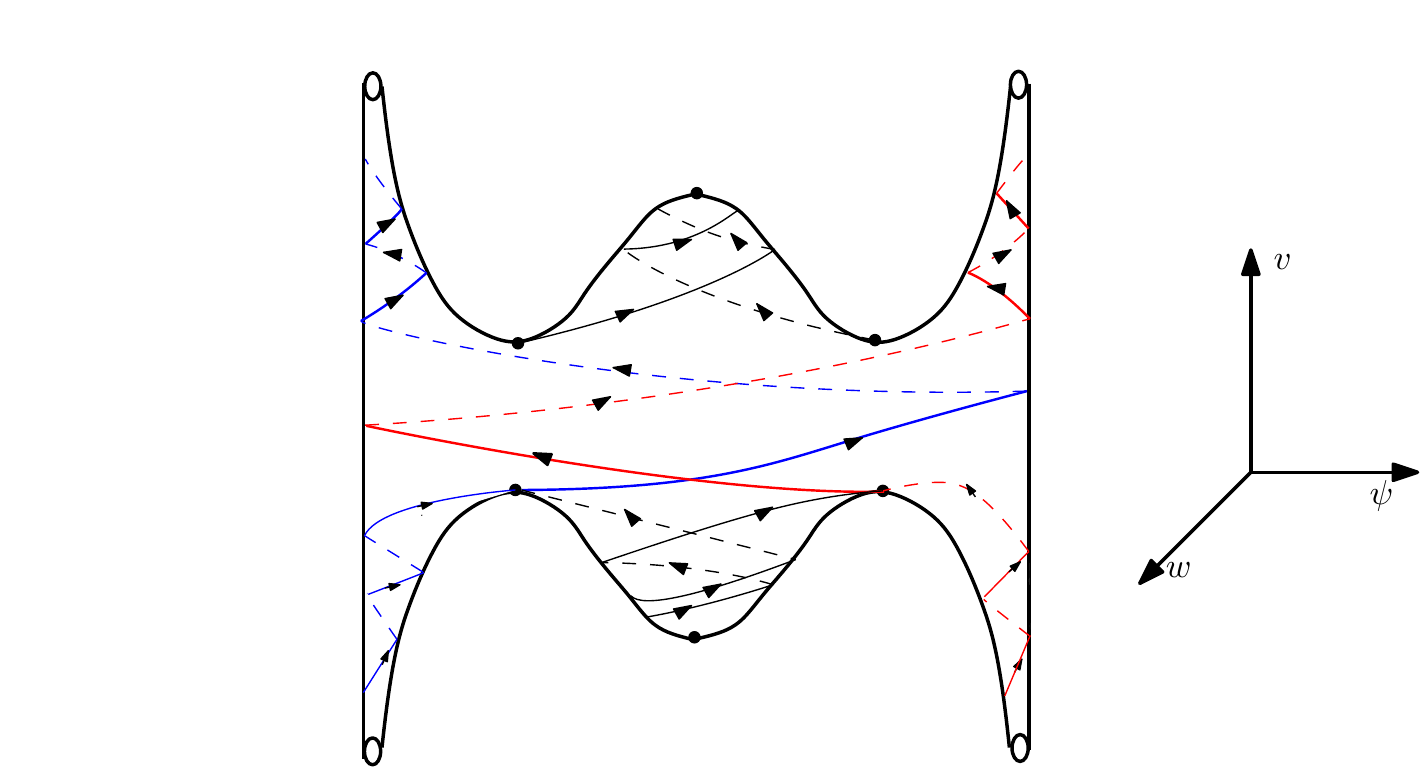}
\label{fig-collision}
\end{center}
\caption{The collision manifold $\cM_0$}
\label{collision-manifold}
\end{figure}

 We first analyze the case when the triple $Q_1$-$Q_3$-$Q_4$ is close to triple collision and $Q_2$ is far apart. As a first approximation, we ignore $Q_2$ and focus only on the $Q_1$-$Q_3$-$Q_4$ three-body problem near triple collision. The problem is called the isosceles three-body problem first studied by Devaney \cite{D}. In Jacobi coordinates  we have
\[\bx_0=(0,x_0)^t,\quad \bx_1=(x_1,0)^t,\quad \bp_0=(0,p_0)^t,\quad \bp_1=(p_1,0)^t,\]
We denote $\bx=(M_0^{1/2}\bx_0,M_1^{1/2}\bx_1)\in \R^4$ and $\bp=(M_0^{-1/2}\bp_0,M_1^{-1/2}\bp_1)\in \R^4$ and by $|\cdot|$ the Euclidean norm, and by $M_0=\frac12,\  M_1=\frac{2m_1}{m_1+2}$ the reduced masses.
We next introduce the blowup coordinates
\begin{equation}\label{EqBlowup}
\begin{aligned}
\begin{cases}
r&=|\bx|= \sqrt {M_0|\bx_0|^2 + M_1|\bx_1|^2},\\
v&=r^{-1/2}(\bx\cdot\bp) ,\\
\psi&=\arctan\frac{\sqrt{M_1}x_1}{\sqrt{M_0}x_0},\\
w& =r^{-1/2}(\sqrt{M_0/M_1} x_0 p_1-\sqrt{M_1/M_0} x_1p_0).
\end{cases}
\end{aligned}
\end{equation} 
The physical meanings of the variables are as follows. The variable $r$ measures the size of the I3BP, $v$ is the projection of the rescaled momentum $r^{1/2}\bp$ to the radial componnet $\bx$, $w$ is the scalar part of $\mathbf w$ (see equation \eqref{EqBlowup1}) where $\mathbf w=r^{1/2}\bp-v r^{-1}\bx$ is the projection of $r^{1/2}\bp$ to the tangential of the sphere $|r^{-1}\bx|=1$ and $\psi$ measures the relative size of the positions $x_0$ and $x_1$.
All variables except $r$ are rescaling invariant and the triple collision corresponds to $r=0$. The coordinate change is accompanied by a time reparametrization $dt= r^{3/2}d\tau$. Equations of motion can be derived in terms of these coordinates (see Section \ref{SBlowupMM}). In particular the $r$-equation is $\frac{dr}{d\tau}=rv$, so we see that $\{r=0\}$ is an invariant submanifold. All the energy level sets are the same in the limit $r\to 0$, denoted by $\mathcal M_0$, which has two dimensions. We call $\mathcal{M}_0$ the collision manifold and the dynamics on $\mathcal{M}_0$ is   illustrated in Figure \ref{collision-manifold}.

\begin{Thm}[\cite{D}] \label{collision-m1}
\begin{enumerate}
\item The collision manifold $\mathcal M_0$ is a topological 2-sphere with four punctures $($called arms$)$, each puncture corresponding to $\psi=\pm\pi/2$. The manifold $\cM_0$ is symmetric with respect to $\psi\mapsto -\psi$ or $w\mapsto -w$.
\item The variable  $v$ is a Lyapunov function in the sense that $\frac{d}{d\tau}v\geq 0$ along the flow on $\mathcal M_0$ and $\frac{d}{d\tau}v=0$ only at the fixed points.
\item There are six fixed points on $\mathcal M_0$: the two with $(v,\psi,w)=(\pm v_\dagger,0,0)$ correspond to Euler (collinear) central configurations and the four with $(v,\psi,w)=(\pm v_*,\pm\psi_*,0)$ correspond to Lagrange (equilateral) central configurations, where $v_\dagger,v_*$ are numbers depending only on the masses and $\psi_*=\arctan\left(\sqrt\frac{3m_1}{2+m_1}\right)$. $($In the following, since we consider only the lower Lagrange point, we fix the convention $v_*<0$$)$. 
\item The Euler fixed point with $v>0$ is a sink and with $v<0$ is a source.  For some masses, the eigenvalues are real and for others complex.  The Lagrange fixed points are saddles.
%\item For a Lagrange fixed point with $v<0$, one of the stable manifolds is homo-clinic to the Euler fixed point with $v<0$; one stable manifold  comes along one lower arm and one unstable manifold  climbs to one upper arm.
\end{enumerate}
\end{Thm}
We are only interested in the Lagrange fixed points in the lower half space $v<0$. The main observation is that we can arrange the orbit to stay close to the saddle for as long a time as we wish, by selecting  the initial condition sufficiently close to its stable manifold. Since in the lower half space we have $v_*<0$,  from the equation $\frac{d}{d\tau}r=rv$ we get that $r$ will decrease to as small as we wish before leaving a neighborhood of triple collision.
%We denote by $W^s_0$ the stable manifold of left lower Lagrange fixed point on $\mathcal M_0$ and by $W^s_{iso}$ its stable manifold in the isosceles three-body problem. We see that $W^s_0$ has one-dimension and is a submanifold of $W^s_{iso}$.
We shall consider the situation that the binary comes from the right to have a near triple collision with $Q_1$. We hope that after the near triple collision the binary moves to the right and $Q_1$ moves to the left. In this case the relative position $x_1$ from the mass center of the pair to $Q_1$ has a negative sign before and after the near triple collision, therefore the variable $\psi$ should also carry a negative sign before and after the near triple collision.  Thus we need the left lower Lagrange fixed point to have a stable manifold coming from the left lower arm and an unstable manifold escaping to the left upper arm.  The case of the near triple collision of $Q_2$-$Q_3$-$Q_4$ is similar. This time the pair comes from the left towards $Q_2$ so that we will need the lower right Lagrange fixed point to have a stable manifold coming from the right lower arm and an unstable manifold escaping to the right upper arm. The latter case follows from the former case by symmetry. So we focus only on the former case. The existence of the stable manifold coming from the left lower arm follows directly from the fact that $v$ is a Lyapunov function. However, the existence of an unstable manifold escaping to the left upper arm is nontrivial, and depends on the mass ratio $m_1/m_3=m_1$. This was studied by \cite{S, SM}. Here we only cite the relevant statements.
\begin{Thm}[\cite{S, SM}]\label{collision-m2} Assume $m_1<55/4$, in which case the Euler fixed points are respectively sink and source with complex eigenvalues. Then there exist $b_1\simeq 0.379$ and $b_2\simeq 2.662$ such that the following holds for the left lower Lagrange fixed point:
\begin{enumerate}
\item when $b_1<m_1<b_2$, the two  unstable manifolds escape from the two upper arms respectively;
\item when $m_1>b_2$, one of the unstable manifold escapes from the left upper arm and the other dies at the upper Euler fixed point.
\end{enumerate}
\end{Thm}
We consider only the masses $(m_1,m_2)\in (b_1,\infty)^2\subset\mathbb{R}^2$. %Note that in these cases the Lagrange fixed points in the lower half have unstable manifolds running  up both arms.

%We are only interested in the Lagrange fixed points in the lower half space $v<0$. The main observation is that if the initial condition is sufficiently close to the stable manifold, then the orbit will stay close to the saddle for as long time as we wish.
%Since in the lower half space we have $v_*<0$, so from the equation $r'=rv$ we get that $r$ will decrease to as small as we wish.

\begin{Not}\label{NotWs}
\begin{enumerate}
\item Denote by $O_1$ $($or $O_L)$ the left lower Lagrange fixed point on the collision manifold of $Q_1$-$Q_3$-$Q_4$, and by $O_2$ $($or $O_R)$ the right lower Lagrange fixed point on the collision manifold of $Q_2$-$Q_3$-$Q_4$.
\item Denote  by $W^s_{\therefore}(O_i)$ the stable manifold of $O_i$, in the isosceles three-body problem $Q_1$-$Q_3$-$Q_4$, $i=1,2$. We see that $W^s_{\therefore}(O_i)$, $i=1,2$ has dimension two.
\item Denote by $W^u_0(O_L)$ the branch of the unstable manifold of $O_L$ on the collision manifold $\mathcal{M}_0$ that escapes from the left upper arm. Similarly define $W^u_0(O_R)$ by changing left to right.
%\item A small neighborhood of $W^s_0(O_i)$ on $\mathcal M_0$ has two sides, $i=1,2$. Initial condition chosen on one side will shadow one of the two unstable manifolds. Choose an orbit which first shadows $W_0^s(O_i)$ and then shadows $W_0^u(O_i)$.  We say that a point of $\mathcal M_0$ lies on the ``correct side" of $W_0^s(O_i)$ if it lies in the region bounded by $W_0^s(O_i), W_0^u(O_i)$, and the chosen orbit.
%
\end{enumerate}
\end{Not}
Note the fact that an orbit starting sufficiently close to the stable manifold $W^s_\therefore(O_i)$ of a Lagrange fixed point will follow the unstable manifold of that point very far. There are two branches of the unstable manifold on $\mathcal M_0$ and there are two sides of the two dimensional manifold $W^s_\therefore(O_i)$ in the three dimensional energy level set. In order for the exiting orbit to follow the correct branch of the unstable manifold $W^u_0(O_L)$ defined above, we have to choose our initial condition on the correct side of $W^s_\therefore(O_i)$.  In the following we use the phrase ``{\it correct side}" of the stable manifold $W^s_{\therefore}(O_i)$, $i=1,2$ with this meaning. %When restricted on $\mathcal{M}_0$, choosing an orbit which first shadows $W^s_{iso}(O_i)\cap\mathcal{M}_0$ then shadows $W_0^u(O_i)$, the points  inside  the region bounded by $W^s_{iso}(O_i)\cap\mathcal{M}_0$, $W_0^u(O_i)$ and the orbit are in the ``correct side".
%and by $W^s_{iso}$ its stable manifold in the isosceles three-body problem. We see that $W^s_0$ has one-dimension and is a submanifold of $W^s_{iso}$.% two dimensions with a new $r$ component.

% In this paper, we always require an orbit to escape to the left upper arm, if it comes in from the left lower arm. So we always choose our initial condition on the correct side of $W^s_0$ when projected to $\mathcal M_0$.

We next introduce the Poincar\'e sections $\mathcal S_\pm$ to separate the near triple collision regime and the perturbed Kepler motions regime.
\begin{Def}
Let us fix $\epsilon>0$ to be a sufficiently small number whose meaning will be reserved throughout the paper. We introduce a Poincar\'e section
$\mathcal S_-:=\{r=\eps^{-1}\}$ before the near triple collision and $\cS_+=\{v=\eps^{-1/2}\}$  after the near triple collision. We introduce superscript $1,2$ to indicate that the section is chosen close to the near triple collision $Q_1$-$Q_3$-$Q_4$ or $Q_2$-$Q_3$-$Q_4$ respectively.  We also use the notation $\cS^{L}_\pm $ and $\cS^{R}_\pm $ interchangeably with $\cS^{1}_\pm$ and $\cS^{2}_\pm$ respectively.
\end{Def}

\subsection{The local and global maps and the key derivative estimate}
\begin{Def}[Local and global maps]\label{DefLocalGlobal}
We define the \emph{local map} $\mathbb L$  to be the Poincar\'e map going from $\cS^{i}_- $ to $\cS^{i}_+$, and \emph{global map} $\mathbb G$  to be the Poincar\'e map going from $\cS^{i}_+ $ to $\cS^{3-i}_-$, $i=1,2$. By default, we use the blowup coordinates for the local map and Delaunay coordinates for the global map, if not otherwise stated.
\end{Def}

\begin{Not} We use the super-script $i$ and $f$ to stand for the corresponding variables on the \emph{initial} and \emph{final} section respectively.
\end{Not}

The following lemma says that the local map is well approximated by the dynamics of the isosceles three-body problem $Q_1$-$Q_3$-$Q_4$ (the $Q_2$-$Q_3$-$Q_4$ case is similar).

\begin{Prop}\label{PropLocalMM} For the isosceles four-body problem,  we have that
there exist $C_1>1,D>0, \chi_0\gg 1$ and $\delta_0>0$ such that for any $\delta\in(\chi^{-D},\delta_0)$ and $\chi>\chi_0$ the following holds: Let  $\zeta:\ [0,1]\to \cS_-^L$ be a smooth curve on the section $\cS_-^L$ satisfying
\begin{enumerate}
\item
 $\zeta$  is on the correct side of $W^s_{\therefore}(O_L)$;
 \item  $\zeta(0)\in W^s_{\therefore}(O_L)$ and $\mathrm{dist}(\zeta(1),W^s_{\therefore}(O_L))=\epsilon^2$;
\item on the curve $\zeta$,  $|\bx_2|\geq \chi$ and  $|\bp_2|=O(1)$.
%\item dist$(\zeta_-,\gamma_I)<\dt,\quad \dt/2<|\psi_--\psi_*|<\delta$ where $\psi_*$ is the $\psi$-coordinate of the point of intersection  between $\gamma_I$ and $\cS_-^L$;
%\item $|w_{0,-}|,|g_0-\frac{\pi}{2}|,|w_1|,|g_1|\leqslant \delta^{C_2}$;
%\item The curve $\xi_-$  is transverse to $W^s_{\therefore}$ on a given energy level and its projection to the $(w_0,g_0)$-plane is transverse to that of $W^s_\therefore$;
%\item $|G_{2,-}|\leqslant \delta^{C_3}$, $|\dot\bx|\sim1$ and  $|\bx_{2,-}|\simeq\chi\geqslant \delta^{-C_3}.$
\end{enumerate}
Then \begin{enumerate}
\item  there exists a subsegment $\bar \zeta$ of $\zeta$, such that  for any point $\mathbf{x}$ on the image of  the curve $\bar\zeta$ under the local map, its $r$-coordinate satisfies
$r^f\in(\delta^{3/2},\delta )$.
\item the oscillation of $L_2$ is estimated as $|L_2^f-L_2^i|\leq \frac{C_ 1}{ \chi^3}$, and
\item The travel time for points on $\bar\zeta$ between the two sections satisfies $C_1^{-1}\log\dt^{-1}<|\tau^f-\tau^i|<C_1\log\dt^{-1}$, and for the original time coordinate $t$, we have $|t^f-t^i|\leqslant C_1$.
\end{enumerate}
\end{Prop}
\begin{Not}
Let $\bu_1,\ldots,\bu_n$ be a tuple of linearly independent vectors. We denote by $\mathcal C_\eta(\bu_1,\ldots,\bu_n)$ the $\eta$-cone around $\bu_1,\ldots,\bu_n$, that is, the set of vectors forming an angle at most $\eta$ with the plane span$\{\bu_1,\ldots,\bu_n\}$. 
\end{Not}
%We next consider the I3BP $Q_1$-$Q_3$-$Q_4$. In blowup coordinates and restricted to a fixed energy level, we induce a map $\mathbb L$ from $\cS_-^1$ to $\cS_+^1$. We use coordinates $(r,\psi)$ for the map by eliminating $v,w$ using the energy condition as well as the section conditions. 
%\begin{Prop}\label{PropdL} Let $\zeta$ be a curve on $\cS_-^1$ whose tangent lies in the cone $\mathcal C_\eta(\hat u)$ and intersects $W^s_\therefore(O_1)$. Then we have the following for the image $T(\cR\mathbb L\zeta)\subset \mathcal C_\eta((1,0)).$  
%\end{Prop}
Immediately after the local map, on the sections $\cS^{L,R}_+$, we introduce the renormalization map as follows.
\begin{Def}[Renormalization map]\label{renormalize-map}
We  denote by $\cR$ the renormalization map, which occurs on the sections $\cS_+^{L,R}$ immediately after the application of the local map. To define the renormalization $\cR$, we partition the $r$-interval $(0,e)$\footnote{The number $e$ is not essential, any real number $a>1$ will do the job.} by points $\{e^{-2n-1}\}_{n\in \{0\}\cup\N}$. For each $r$-interval $[e^{-2n-1}, e^{-2n+1})$ we define $\lambda=e^{2n}/\eps>1$ and define the renormalization map on the section $\cS^L_+$ to be
$$(r, v,\psi,w, L_2,\ell_2)\mapsto (\lambda r, v,\psi,w, \sqrt\lambda L_2,\ell_2).$$
The effect of the renormalization on the Delaunay coordinates and the Cartesian coordinates respectively are as follows 
\begin{equation}
\begin{aligned}
&(L_0,\ell_0,L_1,\ell_1, L_2,\ell_2)\mapsto (\sqrt\lambda L_0,\ell_0,\sqrt\lambda L_1,\ell_1,\sqrt\lambda L_2,\ell_2);\\
&(x_0,p_0;x_1,p_1;x_2,p_2)\mapsto (\lambda x_0,\lambda^{-1/2} p_0;\lambda x_1,\lambda^{-1/2} p_1;\lambda x_2,\lambda^{-1/2}  p_2).
\end{aligned}
\end{equation}
Moreover, we also make the time change $t\mapsto t\lambda^{3/2}$ and the change of Hamiltonian $H\mapsto H/ \lambda$ when applying the renormalization. 
\end{Def}

The main effect of the renormalization is to rescale the semimajor axis of the binary's elliptic motion to order 1 size.
Applying  the renormalization to the outcome of the local map in Proposition \ref{PropLocalMM}, we have the following, whose proof is given in Section \ref{SSRenorm}.
\begin{Prop}\label{PropRenorm}
There exists a constant $C_2>1$ such that after applying the renormalization on the section $\cS_+^{L}$, for the image $\mathcal{R}\mathbb{L}\bar\zeta$ where $\bar\zeta$ is the resulting segment in Proposition \ref{PropLocalMM},  we have
\begin{equation}\label{EqRenorm}\{  L_0, L_1, r_1\eps,L_2/\sqrt\lambda\}\subset (C_2^{-1},C_2),\ r_0<C_2. 
\end{equation}
%Moreover, the tangent vectors $T(\mathcal{R}\mathbb{L}(\zeta_0))\subset\mathcal{C}_\eta((1,0,0,0))\subset T\cS_+^L$, providing the $\sqrt{\lambda}\gg\eta^{-1}$. 
A similar result holds for the right section $\cS_+^R$.
\end{Prop}

On the section $\cS^L_+$, we perform the standard energy reduction to eliminate $(L_1,\ell_1)$ from Delaunay coordinates by fixing the total energy and treating the variable $\ell_1$ as the new time.  So we use  the four variables $(L_0,\ell_0,L_2,\ell_2)$ to parametrize the section $\cS_+^L$. Similarly, on the section $\cS^R_-$, we eliminate $(L_2,\ell_2)$  and use the four variables   $(L_0,\ell_0,L_1,\ell_1)$.

We have the following estimates of the derivatives of the local and global maps. In both cases, we use Delaunay coordinates. Denote by $\cM$ the phase space of the full four-body problem (F4BP), which has ten dimensions after allowing for conservation of momentum and angular momentum.
%
%\begin{Prop}\label{PropGlobC0}
%Let $(X,Y)\in \mathcal S_+^j$ satisfy the conclusion of the previous lemma. Then after the global map, we have on the section $\mathcal S_-^{3-j}$
%$$ L^f_i-L^i_i=O(\eps),\ i=0,1,2. $$
%\end{Prop}
%{\bf A proof is needed somewhere. }

\begin{Prop}\label{PropDL1}
Let $\gamma:\ [0,T]\to \cM$ be an orbit of the I4BP with $\gamma(0)\in \cS_-^L$ and $\gamma(T)\in \cS_+^L$ with $T=\tau^f-\tau^i<C_1\log\delta^{-1}$. 
Then for $\dt$ sufficiently small, we have the following estimate of the derivatives of the local map along $\gamma$
\begin{equation}
d\mathcal Rd\mathbb L=e^{-Tv_*} \mathbf u_1\otimes \mathbf l_1+O(1), %O\left[\begin{array}{cccc}\eps^2\chi&\eps^4\chi&\beta\eps^2\chi&\frac{\eps^2}{\beta^2}\\
%\chi&\eps^2\chi&\beta\chi&\frac{1}{\beta^2}\\
%1&\eps^2&\beta&\frac{1}{\beta^2\chi}\\
%\chi&\eps^2\chi&\beta\chi&\frac{1}{\beta^2}\end{array}\right].
\end{equation}
where using the variables $(L_0,\ell_0,L_2,\ell_2),$ we have $$\bu_1=(1,o(1),o(1),o(1)),\ \bl_1=(o(1),1,o(1),o(1))$$ as $\dt\to0$ and $\chi\to\infty$.
\end{Prop}

\begin{Prop}\label{PropDG1}
There exist $C_2>1$ and $\chi_0\gg1$ such that the following holds:  Let $\gamma:\ [0,T]\to \cM$ be an orbit of the I4BP with $\gamma(0)\in \cS_+^i$ and $\gamma(T)\in \cS_-^{3-i},\ i=1,2$ and  the initial condition $\gamma(0)$ satisfying \eqref{EqRenorm} with $C_2$, and in addition
$|x_2|= \chi\geq \chi_0.$
Then we have the following estimate of the derivatives of the global map along $\gamma$ (with $\beta=1/\sqrt\lambda$)
\begin{equation}
d\mathbb G=\chi \bar\bu_1\otimes \bar\bl_1+O((\eps^2+\beta)\chi), %O\left[\begin{array}{cccc}\eps^2\chi&\eps^4\chi&\beta\eps^2\chi&\frac{\eps^2}{\beta^2}\\
%\chi&\eps^2\chi&\beta\chi&\frac{1}{\beta^2}\\
%1&\eps^2&\beta&\frac{1}{\beta^2\chi}\\
%\chi&\eps^2\chi&\beta\chi&\frac{1}{\beta^2}\end{array}\right].
\end{equation}
where $\bar\bu_1=(0,1,0,O(1)),\ \bar\bl_1=(1,0,0,0)$.
\end{Prop}
We will prove  this proposition in Section \ref{SDerLl}. The key point in this estimate is the $(2,1)$-entry $\frac{\partial \ell_0}{\partial L_0}=O(\chi)$, which implies that a small change in the initial $L_0$ will lead to a huge change in the final phase variable $\ell_0$.

From the last two propositions, it is clear that we have the transversality conditions $\bl_1\cdot \bar\bu_1\neq 0$ and $\bar\bl_1\cdot \bu_1\neq 0$.  We define the Poincar\'e map $\mathbb P:=\mathbb G\cR\mathbb L:\ \cS_-^i\to \cS_-^{3-i},\ i=1,2$. Thus we have the following cone preservation property. We fix a small number $\eta>0$ and choose $\dt$ small and $\chi$ large accordingly. 
\begin{Prop}Let $\bx\in \mathcal S_-^i$ be the initial condition for an orbit of the I4BP such that $\mathbb P(\bx)\in \mathcal S_-^{3-i}, \ i=1,2$. Suppose the assumptions of  Proposition \ref{PropDL1} and \ref{PropDG1} are satisfied along the orbit, then we have 
$$(d_\bx\mathbb P)\mathcal C_\eta(\bar\bu_1)\subsetneq \mathcal C_\eta(\bar\bu_1).$$
Moreover, for each $v\in \mathcal C_\eta(\bar\bu_1)$, we have $|d_\bx\mathbb P(v)|\geq e^{-Tv_*}\chi |v|. $
\end{Prop}

\subsection{The I3BP near infinity}
%We shall also consider solutions of the isosceles three-body problem with the binary and the third body moving off to infinity in opposite directions. We use another coordinate change of McGehee.

%We may view the entire $y$-axis as a normally hyperbolic invariant manifold.

%In this coordinates, we get that the infinity is a normally hyperbolic invariant manifold topologically a cylinder parametrized by $(L_0,\ell_0)$.
We will need the following special orbits $\gamma_{I}$ and $\gamma_O$ as guides of our noncollision singular orbits. Let us motivate the definition briefly.   Suppose the triple $Q_1$-$Q_3$-$Q_4$ completes its near triple collision so the energy of subsystem $(\bx_0,\bp_0)$ (the binary) and $(\bx_1, \bp_1)$ (the relative motion of $Q_1$ and the mass center of the binary) both have very large $(\gg 1)$ modulus but add up to a total energy of order 1 by the local energy almost conservation. Now the energy of $(\bx_2,\bp_2)$ is also of order 1 since $Q_2$ was not involved in the near triple collision. We denote by $\lambda(\gg 1)$ the order of magnitude of the energy of $(\bx_1,\bp_1).$ We renormalize the whole system so that the energy of $(\bx_0,\bp_0)$ and $(\bx_1,\bp_1)$ are of order 1, so they almost cancel each other and the energy of $(\bx_2,\bp_2)$ becomes almost zero (so the total energy is also almost zero). Now the binary moves towards $Q_2$ and we need to transform to right Jacobi coordinates. The relative motion of $Q_2$ and the mass center of the binary has order 1 velocity so in the right Jacobi coordinates, both $(\bx_0,\bp_0)$ and $(\bx_2,\bp_2)$ have order 1 energy, and so does $(\bx_1,\bp_1)$. The energy partition is almost determined by the mass ratios. Suppose the triple $Q_2$-$Q_3$-$Q_4$ goes to triple collision. We expect that in the limit $\lambda\to \infty$, the triple collision orbit converges to a certain orbit $\gamma_I$ on the stable manifold $W^s_{\therefore}(O_2)$, with a fixed energy partition at infinity.

Next, if the triple $Q_2$-$Q_3$-$Q_4$ stays close to the right lower Lagrange fixed point for a long time and escapes to the right upper arm along the unstable manifold, we expect $r$ remains small for a long time. So we use $\gamma_O$, that is the unstable manifold of the right lower Lagrange fixed point on the collision manifold $\cM_0$, to guide the orbits leaving the triple collision. The existence of $\gamma_I$ is given by the following theorem, whose proof is given in Section \ref{SSgammaI}.

%\begin{Def} \label{NotTrans}
 %Denote by $\gamma^2_I:\ \R\to \R^{6}$ $(``I"$ means $``in"$, $2$ records $O_2$ and $\R^6$ is the phase space for the full three-body problem$)$ the orbit of the isosceles three-body problem $Q_2$-$Q_3$-$Q_4$  lying on the stable manifold of $O_2$ and with the boundary condition $$\left(\frac{L_0(-\infty)}{L_2(-\infty)}\right)^2=\frac{m_1}{16m_2^2(m_1+2)},\quad \gamma^2_I(0)\in \cS^2_-.$$ Note that $\gamma^2_I$ does not lie on the collision manifold $\cM_0$. Denote by $\gamma^1_O:\ (-\infty,0]\to \R^6$ ($``O"$ means $``out"$) the right branch of unstable manifold of $O_1$ on the collision manifold $\cM_0$ of the isosceles three body problem $Q_1$-$Q_3$-$Q_4$ such that $\gamma_O^1(0)\in \cS^1_+$. Similarly, we introduce $\gamma^1_I$ and $\gamma^2_O$.
%\end{Def}

\begin{Prop} \label{ThmGammaI} There exists an orbit of the I3BP $Q_2$-$Q_3$-$Q_4$ lying on the stable manifold $W_\therefore^s(O_2)$ of the Lagrange fixed point $O_2$ and satisfying the boundary condition $\left(\frac{L_0(-\infty)}{L_2(-\infty)}\right)^2=\frac{m_1}{16m_2^2(m_1+2)}$ as $t\to-\infty$.
\end{Prop}
\begin{Not}
\begin{enumerate}
\item We will call the orbit in Proposition \ref{ThmGammaI} $\gamma_I$ and the unstable manifold $W^u_0$ in Notation \ref{NotWs} $\gamma_O$.
\item We shall use the superscript $1$ or $2$ to indicate the orbits $\gamma_I$ and $\gamma_O$ defined for the I3BP $Q_1$-$Q_3$-$Q_4$ or $Q_2$-$Q_3$-$Q_4$ respectively.
\item For $\gamma_I$, we assume $\gamma_I(0)\in \cS_-$ and for $\gamma_O$ we assume $\gamma_O(0)\in \cS_+$, after time translations if necessary.
\end{enumerate}
\end{Not}

We may view infinity as a center manifold and orbits approaching or leaving it hyperbolically as stable and unstable manifolds (c.f. \cite{DMMY}). Then $\gamma_I$ is the intersection of the unstable manifold of infinity and the stable manifold of the Lagrange fixed point. However, we do not adopt this viewpoint in this paper. People have studied a similar problem of finding orbits on $W^s_\therefore(O_i)$ approaching infinity parabolically as $t\to-\infty$ (c.f. \cite{Mo4,SM}).
\begin{Not}
\begin{enumerate}
\item Denote by $\cM_\therefore$ the phase space of the full three-body problem (F3BP), which has six dimensions after reducing the conservation of momentum and angular momentum. We shall use a superscript $i=1,2$ to stand for the triple $Q_i$-$Q_3$-$Q_4$, if necessary;
\item Denote by $\pi^i_\therefore$ the projections to the $\cM^i_\therefore,\ i=1,2$ subsystem respectively.
\end{enumerate}
When the triple $Q_i$-$Q_3$-$Q_4$ is close to triple collision in the sense that it lies between the two sections $\cS^i_-$ and $\cS^i_+,\ i=1,2$, we always use the blowup coordinates.
\end{Not}
\begin{Prop}\label{PropDevaneyOrbit}
Let $\dt$ and $\chi$ be as in Proposition \ref{PropLocalMM}, satisfying $\dt>\chi^{-D}$.
Let $\gamma:\ [0,T]\to \cM$ be an orbit of the I4BP such that
\begin{enumerate}
\item  $\gamma(0)\in \cS_+^L$ with $r(0)<\dt;$
\item  $(\mathcal R\gamma)(T)\in \cS_-^R$ satisfies $\mathrm{dist}(\pi^2_\therefore(\mathcal R\gamma)(T),W^s_\therefore(O_2))<\dt$.
\end{enumerate}
Then as  $\dt\to 0$, we have
%\begin{enumerate}
%\item
$\mathrm{dist}(\pi^2_\therefore(\mathcal R\gamma)(T),\gamma_I(0))\to 0$.
%\item denote by $T''$ the time $\gamma(T'')\in \cS^R_-$ and choose the renormalization such that $L_0(T')=1$ on the section $\cS^L_+$, then we have that  $\pi^2_\therefore(\cR(\gamma(\cdot+T'')))\to \gamma^2_I$ uniformly on every compact set $[c,d]\subset [-T'',T-T'']$.
%\end{enumerate}
The same thing holds when we permute $1$ and $2$.
\end{Prop}

\subsection{Proof of Theorem \ref{ThmMM}}

We now give the proof of Theorem  \ref{ThmMM}, assuming the propositions in the previous subsection.
The proof consists of two parts. In the first part, we show that there exists a Cantor set of  initial conditions such that the map $\mathbb P$ can be iterated for infinitely many steps. In the second part, we show that such  initial conditions indeed lead to  unbounded orbits in finite time.

We shall prove the following quantitative version of Theorem \ref{ThmMM}.

\begin{Def}[Sojourn time]\label{DefSojourn}
Let $\{T_n\}$ be a sequence satisfying $T_n>C_1^{-1}\log \dt_0^{-1}$, where $\delta_0$ is as in Proposition~\ref{PropLocalMM}. We say that $\{T_n\}$ is a sequence of sojourn times of type-$(c,C),\ c>0,C>0$ if it satisfies in addition
$$T_{n+1}\leq c\sum_{i=1}^n T_i+C,\ \forall \ n\geq 1. $$
\end{Def}

\begin{Prop}\label{ThmMMQ}
There exist constants $c_0$ and $C_0$ such that for any sequence $\{T_n\}$ of sojourn times of type-$(c_0,C_0)$, there exists an initial condition $\sx$ whose forward orbit $\{\mathbb P^n\sx\}$ is well-defined and the time duration $\tau_n$ $($in blowup coordinate$)$ defining the local map in the $n$-th iteration satisfies $|\tau_n-T_n|<C_2$ where $C_2$ is a uniform constant depending on $\eps$, $C_1$,and $D$.
\end{Prop}

%For  the proof of Theorem \ref{ThmMM}, we needs only Proposition \ref{PropDLl}.Since we consider the isosceles problem, we have $G_0=g_0-\pi/2=G_1=g_1=G_2=g_2=0$. Then  the dynamics is completely described by the six  Delaunay coordinates $(L_0,\ell_0,L_1,\ell_1,L_2,\ell_2)$.

%We consider the  two sections $\cS_+^L$ and $\cS_-^R$.
\begin{proof}
From now on we only consider the phase points on $\cS_+^L$ such that condition \eqref{EqRenorm} is satisfied. Clearly, they form an open set on the section $\cS^L_+$.  Fix a small number $\eta>0$.
For any such a phase point $X=(L_0,\ell_0,L_2,\ell_2)\in \cS^L_+$, we introduce the cone $\mathcal C_\eta(\bu_1)\subset T_X\cS^L_+$. %to be the set of vectors forming an angle less than $\eta$ with the direction  $(1,0,0,0)$. We assume $0<\eta\ll\epsilon^2$.

%Let us consider a section $\cS^R_-=\{v_-=v_*- C_\epsilon\}$ near the triple collision of the triple $Q_2$-$Q_3$-$Q_4$. Since $Q_1$ is far, the system almost decouples into a product system of the triple $(x_0,p_0,x_2,p_2)$ and the two-body problem $(x_1,p_1)$. We denote by $\pi_{0,2} \cS_-$ the projection of the section $\cS_-$ to the subspace of $(x_0,p_0,x_2,p_2)$.

We have the following non-degeneracy property for the global map, proved in Section \ref{SSNondeg}.
\begin{Lm}\label{LmNondeg}
Let $\xi$ be an initial segment of length $\frac{1}{\epsilon\chi}$ on the section $\cS^L_+$, with all its tangent vectors lying inside the cone $\mathcal C_\eta(\bu_1)$  for each point $X$ on $\xi$.  Then its image under the global map  on the section $\cS_-^R$  winds  around the cylinder formed by $(L_0,\ell_0)$ along the direction $\ell_0$ many times.  In particular, it intersects transversely with the stable manifold $W_{\therefore}^s(O_2)$ when projected to $\cM_\therefore$.
 \end{Lm}

Let us start with an initial segment $\zeta_{0}$ on the section $\cS_-^L$ satisfying the assumption of Proposition~\ref{PropLocalMM}. %We  parametrize the curve $\gamma_{0,-}$ by $\psi\in \psi_*+[\psi_0,\psi_1]$ with $\frac{\delta}{2}\leqslant \psi_0<\psi_1\leqslant\delta$ and assume further that its tangent vectors when projected to the component $(L_2,\ell_2)$ are of order $O(1)$.
By Proposition \ref{PropLocalMM}, we know that for  the image $\mathbb L\zeta_0$ of the curve $\zeta_{0}$ under  the local map to  the section $\cS_+^L$, its  $r$-component ranges from $\delta^{3/2}$ to $\delta $.  %Therefore, for the tangent vector of the curve $\gamma_{0,+}$, the ratio between its $L_2$ component and $L_0$ component is of order $\delta^{\frac{4}{5}v_*}\ll\eta$.

Applying the renormalization map, we focus on one segment of the partition of $\mathbb L \zeta_{0}$. Each such segment is rescaled to
 length of order one and the $r$-component has values in the interval  $(e^{-1}\eps^{-1},e\eps^{-1})$. Each segment is a smooth curve
  and the phase points on it satisfy condition \eqref{EqRenorm}.  The renormalization map stretches the $L_0$ and $L_2$ variables
  with the same ratio $\sqrt\lambda$ and leaves $\ell_0$, $\ell_2$  untouched. Therefore  after renormalization, each segment $\mathcal R\mathbb L\zeta_{0}$ has tangent vectors lying in the cone $\mathcal C_\eta(\bu_1)$ by Proposition \ref{PropDL1}, provided the rescaling factor $\lambda\gg \eta^{-1}$.

Next, we apply the global map to obtain the strong expansion in Lemma \ref{LmNondeg}. The strong expansion shows that after applying the global map $\mathbb G$ and arriving at the section  $\cS^R_-$, the resulting curve
$\mathbb P\zeta_0$ winds around the $\ell_0$-circle many times.

Then on the section $\cS_-^R$, we are in a position to pick from the image  $\mathbb P\zeta_{0}$ segments with phase points satisfying the assumptions of Proposition  \ref{PropLocalMM}.   This involves deleting many open intervals from the  segment $\mathcal{R}\mathbb L\zeta_{0}$, then in the original segment $\zeta_{0}$. We then repeat the above procedure to  the newly picked segments on the section $\cS_-^R$.  Repeating this procedure infinitely many times,   in the limit, we get a Cantor set $\mathcal{I}_{\zeta_{0}}$ on the curve $\zeta_{0}$ after deleting open intervals for infinitely many steps.

We finally prove the statement on the sojourn time. From Proposition \ref{PropLocalMM}, we estimate each travel time $\tau_n$ of the $n$-th application of the local map by $\log \dt^{-1}$ up to a multiple $C_1$, where $\dt$ satisfies $\chi^{-D}<\dt<\dt_0$. This shows that $C_1^{-1}\log \dt_0^{-1}<\tau_n<C_1D\log \chi_n$. The distance $\chi_n$ is obtained by multiplying $\chi_0$ by a renormalization factor $\lambda_k$ during each step of renormalization where $\lambda_k r/\eps \in (1,e)$ and $r$ is the value of $r$ immediately after the $k$-th application of the local map. We thus have $\chi_n\geq\prod_{k=1}^{n-1} \lambda_k  \chi_0$. We then relate $\lambda_n$ to $\dt$ by the estimate $r\in (\dt^{3/2},\dt)$ then to $\tau_n$ as above. This gives the sojourn time estimate if we choose $C_0=C_1D\log \chi_0$ and $c_0=C_1^2D$.
\end{proof}

 \begin{proof}[Proof of Theorem \ref{ThmMM}] We next show that each initial condition in the Cantor set $\mathcal{I}_{\zeta_{0}}$ leads to an unbounded orbit in finite time. Note that here  we allow double collisions between the binary $Q_3$-$Q_4$ but do not allow any other type of double or triple collision.

Consider $\sx_0\in\mathcal{I}_{\zeta_{0}}$ and assume that the initial distance between the particles $Q_1$ and $Q_2$ is $\chi_0\gg1$. Let $\lambda_i$, $i=1,2,\dots$ be the renormalization factor of the  orbit  starting at $\sx_0$ when the renormalization map is applied the $i$-th time.
We  estimate the time span of the orbit by induction. %From Proposition \ref{PropLocalMM}, we require $\frac{1}{\chi}\leqslant \delta_1^{C'}$.
Suppose the binary comes to a near triple collision with $Q_1$.
Then after the local map and the first renormalization, we have that
the distance between $Q_1$ and $Q_2$ is $\chi_1\sim \lambda_1\chi_0$, while $Q_1$ and $Q_2$ are moving in opposite directions with  $|\dot Q_1|\sim1$,  $|\dot Q_2|\sim\lambda_1^{-1/2}$, and  the binary $Q_3$-$Q_4$ moves towards $Q_2$ with speed of order~$O(1)$  by Proposition \ref{PropRenorm}. Here the notation $\sim$ means up to a scalar multiple that is bounded away from both 0 and $\infty$ uniformly.
 Then when the orbit reaches $\cS_-^R$, the distance between $Q_1$ and $Q_2$ is $\sim (1+C)\lambda_1\chi_0$ with $C>0$, and  for the time $t_1$  that (after the renormalization) the orbit spends between the sections $\cS_{+}^L$ and $\cS_-^R$, we have $t_1\sim (1+C) (1+\lambda^{-1}_1)\lambda_1\chi_0$.
In the original spacetime scale without the renormalization, the distance between $Q_1$ and $Q_2$ is $\bar d_1\sim (1+C)\chi_0$ and  the traveling time is  $\bar{t}_1\sim C\lambda_1^{-1/2}\chi_0$.

We next perform the induction. From Proposition \ref{PropLocalMM}, we see that $\chi_i\ge\lambda_i\chi_{i-1}\to \infty$ hence we can take $\dt_i\to 0$ and by the definition of the renormalization, we can take the rescaling factor $\lambda_i\to \infty$. Since the rate of acceleration is proportional to $\lambda_i$, it can be chosen to be arbitrary large.

Then, arguing inductively, we have that in the original spacetime scale without the renormalization,  the distance  between $Q_1$ and $Q_2$ is $\bar d_{i}\sim (1+C)^i\chi_0$ and the time $\bar t_{i+1}$ that the orbit spends between  the $i$-th and $(i+1)$-th renormalization satisfies $\bar t_{i+1}\sim C^i\prod_{k=1}^{i}\lambda^{-1/2}_k\chi_0$.  Clearly, $\bar d_{i}\to \infty$ as $i\to\infty$ and   $\sum_{i=1}^\infty \bar t_i<+\infty$ if  $C<\lambda_i^{1/2}$ for all $i$. Therefore, the orbit starting from $\sx_0$ becomes unbounded in finite time.

The proof of  Theorem \ref{ThmMM} is now complete.
\end{proof}
\section{Proof of Theorem \ref{ThmMain}}\label{SMainProof}
In this section, we prove Theorem \ref{ThmMain} for the F4BP.
The definitions of the sections $\cS^{L/R}_\pm$, and the local, global and renormalization maps for the F4BP, are similar to the isosceles case. We will show in this section that the noncollision singularity of the F4BP is obtained by slightly perturbing the I4BP orbit from Theorem~\ref{ThmMM}. The main issue is to control the extra variables
$(G_0,g_0,G_1,g_1,G_2,g_2).$ %which represents the main difficulty of the Anosov conjecture. 
\subsection{The coordinates and the diagonal form}
We always assume the total angular momentum is vanishing, i.e.
$G_0+G_1+G_2=0.$ Since the Hamiltonian does not depend on each individual angle but depends only on the relative angles, this reduces the number of degrees of freedom by 1.  By further restricting to the zero energy level and to the Poincar\'e sections, we can reduce two more variables when defining the local and global maps. The variables to be removed are chosen to be $L_1$ amd $\ell_1$ on the section $\cS_\pm^L$ and $L_2$ amd $\ell_2$ on the section $\cS^R_\pm$. We denote by $$X^L=(L_0,\ell_0,L_2,\ell_2), \ Y^L=(G_0,\mathsf g_{01}:=g_0-g_1-\frac{\pi}{2},G_2,\mathsf g_{21}:=g_2-g_1-\pi),$$
$$X^R=(L_0,\ell_0,L_1,\ell_1),\  Y^R=(G_0,\mathsf g_{02}:=g_0-g_2-\frac{\pi}{2},G_1,\mathsf g_{12}:=g_1-g_2-\pi)$$
the variables on the sections $\cS^L_\pm$ and $\cS^R_\pm$ respectively. 
 We shall estimate the matrices $d\mathbb G=\frac{\partial (X^R,Y^R)|_{\cS^R_-}}{\partial (X^L,Y^L)|_{\cS^L_+}}$,$\ d\mathbb L=\frac{\partial (X^R,Y^R)|_{\cS^R_+}}{\partial (X^R,Y^R)|_{\cS^R_-}}$ and the same maps with $L$ and $R$ switched. As we state in the next lemma (proved in Section \ref{SG0g0}),  the derivative matrix has a product structure in the isosceles limit, i.e. we evaluate the derivative matrix in the limit $|Y|\to 0$ (along an I4BP orbit). %which significantly simplifies the calculations.  %in the sense that instead of estimating $100$ matrix elements, it is enough to estimate \mbox{$4\times 4+4\times 4+2\times 2=36$} elements.

\begin{Lm}\label{LmDiagonal} The derivative matrix $\frac{\partial (X^R,Y^R)|_{\cS^R_-}}{\partial (X^L,Y^L)|_{\cS^L_+}}$  evaluated along an I4BP orbit has the following structure $$\frac{\partial (X^R,Y^R)|_{\cS^R_-}}{\partial (X^L,Y^L)|_{\cS^L_+}}=\left[\begin{array}{ccc}
\frac{\partial (X^R)|_{\cS^R_-}}{\partial (X^L)|_{\cS^L_+}}&0\\
0&\frac{\partial Y^R|_{\cS^R_-}}{\partial Y^L|_{\cS^L_+}}\\
%0&0&\frac{\partial (Y_1,Y_2)|_{\cS^R_-}}{\partial (Y_1,Y_2)|_{\cS^L_+}}\\
\end{array}\right].$$
The same diagonal form of the derivative matrix also holds for the local map, hence for the Poincar\'e return map. 
\end{Lm}
The first diagonal block  is estimated in Proposition \ref{PropDL1} for the local map and \ref{PropDG1} for the global map. We next show that how to estimate the second diagonal block. 
\subsection{The local map}\label{SSLoc}
When the subsystem $Q_1$-$Q_3$-$Q_4$ is near triple collision, we shall treat the F4BP as a small perturbation of the product system of the F3BP and a Kepler problem $(\bx_2,\bp_2)$.  Near triple collision of the  $Q_1$-$Q_3$-$Q_4$ subsystem, besides the blowup coordinates $(r,v,\psi,w)$ that we used in the I3BP, we introduce the new scaling invariant variables $Y:=(w_0=r^{-1/2}G_0,\mathsf g_{01},w_2=r^{-1/2}G_2,\mathsf g_{21})$ (see Appendix \ref{app-blowup}). %Though the Delaunay variables $(L_0,\ell_0,L_1,\ell_1)$ are not well-defined near triple collision, the variables $(w_0,\mathsf g_{01},w_2,\mathsf g_{21})$ are well-defined there.

We first consider the F3BP $Q_1$-$Q_3$-$Q_4$. Assuming that the total angular momentum $G_0+G_1=0$, the system has three degrees of freedom, and we introduce the variables $w_0,\mathsf g_{01}$ in addition to the blowup coordinates $(r,v,\psi,w)$ of the I3BP. The diagonal structure of Lemma \ref{LmDiagonal} also holds for the F3BP, so we can talk about the second diagonal block that is $2\times 2$ by projecting to the $(w_0,\mathsf g_{01})$-plane. %The phase space of the I3BP is an invariant submanifold of codimension 2 of that of the F3BP with $\mathsf g_{01}=\pi/2$ and $w_0=0$. 
\begin{Lm}\label{LmDF3BPLoc}
When linearized at each Lagrange fixed point, the equations of motion of the F3BP projected to the $(w_0,\mathsf g_{01})$-plane have  one  positive and one negative eigenvalue, $\mu_u$ and $\mu_s$ respectively. 
\end{Lm}
 %Therefore, the Lagrange fixed point $O_1$ of the F3BP has three dimensional stable manifold $W_3^s(O_1)$ (containing $\gamma_I$) and two dimensional unstable manifold $W_3^u(O_1)$ (containing $\gamma_O$).

%For the local map, the second diagonal is the same as $d(\widehat{\cR \mathbb L})$.  

\begin{Not}
 Denote by $E_3^s(\gamma_I(\tau)),\ \tau\in [0,\infty)$ the tangent vector $($normalized to have length $1)$ of $W_3^s(\gamma_I(\tau))$ at $(0,\frac{\pi}{2})$ in the $(w_0,\mathsf g_{01})$-plane $($for the local map$)$. Similarly, we define $E_3^u(\gamma_O(\tau))$, $\tau\in (-\infty,0]$. We shall sometimes refer to the subspaces $E_3^s(\gamma_I(\tau))$ or $E_3^u(\gamma_O(\tau))$, by which we mean the spaces spanned by these vectors.
\end{Not}

We next consider the F4BP, taking into consideration $Q_2$.  Assuming the total angular momentum vanishes, we use coordinates $Y$ in addition to the blowup coordinates for the I3BP. Again by the diagonal structure of Lemma \ref{LmDiagonal}, we can talk about the second diagonal block  by projecting to the $Y$-components.  %Let $\gamma:[0,T]\to\cM$ be an orbit of the I4BP with $\gamma(0)\in \cS_-^1,\ \gamma(T)\in \cS_+^1$ and $\gamma_(0)$ is sufficiently close to $\gamma_I(0)$.
 
 \begin{Lm}\label{LmDF4BPLoc}\begin{enumerate}
\item  The linearized dynamics at the lower Lagrange fixed point, when projected to the $Y$-component, has four eigenvalues $\mu_u>-\frac{v_*}{2}>0>\mu_s$, with the same $\mu_u$ and $\mu_s$ as in Lemma \ref{LmDF3BPLoc}. Denoting by $\mathbf e_1,\mathbf e_2,\mathbf e_3,\mathbf e_4\in \R^4$ the respective eigenvectors, we have that the projection of $\mathbf e_1$ and $\mathbf e_4$ to the $(w_0,\mathsf g_{01})$-plane are eigenvectors for the linearized dynamics in Lemma \ref{LmDF3BPLoc}. Moreover, the vector $\mathbf e_2$ has nontrivial projection to the $G_2$ component, and so does its pushforward along $\gamma_O$ by the tangent dynamics. 
\item The tangent dynamics along $\gamma_I$ and $\gamma_O$ commutes with the projection to the $(w_0,\mathsf g_{01})$-plane $($the F3BP$)$. 
\end{enumerate}
 \end{Lm}
 The proof of the last two lemmas will  be given in Section \ref{local-gs}.
   We can construct the 2 dimensional plane $E^s_{4}(\gamma_I(0))\subset T_{\gamma_I(0)}\cS_-$ such that its pushforward by the tangent flow along $\gamma_I(\tau)$ as $\tau\to\infty$ is $\mathrm{span}\{\mathbf e_3,\mathbf e_4\}$. Similarly, we construct the 2 dimensional plane $E^u_{4}(\gamma_O(0))\subset T_{\gamma_O(0)}\cS_+$ whose pushforward by the tangent flow along $\gamma_O(\tau)$ as $\tau\to-\infty$ is $\mathrm{span}\{\mathbf e_1,\mathbf e_2\}$. Moreover, by Lemma \ref{LmDF4BPLoc}, we know that the projection to the $(w_0,\mathsf g_{01})$-component  of the pushforward of $\mathbf e_1$ to $T_{\gamma_I(0)}\cS_-^1$(respectively  the pushforward of $\mathbf e_4$ to $T_{\gamma_O(0)}\cS_+^1$) is $E^s_3(\gamma_I(0))$ (respectively $E^u_3(\gamma_O(0))$). %(c.f. Lemma \ref{LmFundLoc}). 

After the application of the local map, we perform a renormalization on the section $\cS_+^{L/R}$. The definition of the renormalization is to extend Definition \ref{renormalize-map} to the four-body problem.  The renormalization does not change the values of the variables $(w_0,\mathsf{g}_{01},w_2,\mathsf{g}_{21})$, but we need to use the variables $(G_0,\mathsf{g}_{01},G_2,\mathsf{g}_{21})$ in order to apply the global map. By the definition of the renormalization, the $r$-variable is rescaled to a scalar of order $1/\eps$, so the value of $G_i$ is about $1/\eps^{1/2}$ times of that of $w_i,\ i=0,2$. 

%Next, we consider  the F3BP that is a small perturbation of the I3BP with nonisosceles variables $(G_0,g_0)$. Linearizing at the lower Lagrangian fixed point, we get two eigenvalues $\mu_u>0>\mu_s$. We define $E^s_\therefore(\gamma_I(0))$ and $E_\therefore^u(\gamma_O(0))$ similarly as above in the current setting. 
%We denote by $d(\widetilde{\cR \mathbb L})$ the projection/restriction of $d\cR\mathbb L$ to the  $(G_0,g_0)$-components in the F3BP when evaluated along an orbit of the I3BP. 

In the following propositions, without danger of confusion, when talking about the vectors $\bu_1,\bl_1$ that were defined in Proposition \ref{PropDL1} for the I4BP, we naturally embed them into the phase space of the F4BP by adding zeros to the extra dimensions. Similarly for the subspaces $E^{u/s}_{4}$, $E^{u/s}_{3}$, etc. Moreover, though the basepoint of the subspaces $E^{s}_{4}$, $E^{s}_{3}$ is $\gamma_I(0)$, we naturally parallel transport these subspaces to nearby points of $\gamma_I(0)$. Similarly for $E^{u}_{4}$, $E^{u}_{3}$. 
\begin{Prop}\label{PropDL2} 
For all $\eta>0$, there exists a small number $d_0>0$ such that the following holds: Let $\gamma:\ [0,T]\to \cM$ be an orbit of the F4BP in the blowup coordinates with $\gamma(0)\in \cS_-^1,\ \gamma(T)\in \cS_+^1$ and $\gamma(0)$ is sufficiently close to $\gamma_I(0)$ that
$|Y|<d_0 e^{-\mu_u T},\  \chi>d_0^{-1}e^{\mu_u T}$ along $\gamma$. Then we have $$d(\cR \mathbb L)=(e^{-T v_*} \bu_1\otimes \bl_1)\oplus d(\widehat{\cR \mathbb L})+O(d_0),$$
\begin{equation}\label{EqdRL}
d(\widehat{\cR \mathbb L})P\subsetneq \cC_\eta(E^u_{4}(\gamma_O(0))),\quad \forall\ \textrm{2-plane}\ P\  \mathrm{with}\ P\cap \cC_\eta(E^s_{4}(\gamma_I(0)))=\{0\}. \end{equation}
where 
$\bu_1$ and $\bl_1$ are as in Proposition \ref{PropDL1}. Moreover, for all $v\in P$, we have$$  \|d(\widehat{\cR \mathbb L})v\|\geq d_0e^{-Tv_*/2}\|v\|.$$
%Moreover, we have \begin{equation}\label{EqPCRP}\Pi(d\widehat{\cR \mathbb L})\Pi=d\widetilde{\cR \mathbb L},\end{equation} where $\Pi$ means the projection to the $(G_0,\mathsf g_{02})$-components. 
%for each vector $v\in \cC_\eta(\bu_1, \bu_2,\bu_3)$, we have $\|d\cR d\mathbb L v\|\geq e^{-Tv_*/2}\|v\|$.
\end{Prop}

\subsection{The global map}We next provide the statement for the  global map. 

\begin{Prop} \label{PropDG2}  
For all $\eta>0$, there exist $\chi_0,\ \eps_0, \nu_0$ such that for all $\chi>\chi_0,\eps<\eps_0, \nu<\nu_0$,  we have the following: Let $\gamma:\ [0,T]\to\cM$ be an orbit of the F4BP with $\gamma(0)\in \cS^L_+$, $\gamma(T)\in 
\cS^R_-$ and $|Y|\leq \nu$ along the orbit.  Then there exist vector fields $\bar\bu_i,\ i=1,2,3,$ and $\bar\bl_i$,  $i=1,2$, such that we have 
$$d\mathbb G=\chi \bar \bu_1\otimes \bar\bl_1+\chi \bar \bu_2\otimes \bar\bl_2+O((\eps+\nu+\lambda^{-1})\chi), $$
%Moreover,  for each vector $v\in \cC_{\eta/2}(\bv_1, \bv_2,\bv_3)$, we have $\|d\mathbb G v\|\geq \frac12\|v\|$.
and 
\begin{equation}\label{EqDG}d\mathbb G \cC_{\eta}(\bu_1, E_{4}^u(\gamma_O(0)))\subsetneq\cC_{2\eta}(\bar\bu_1, \bar\bu_2,\bar \bu_3),\end{equation}
and for each vector $v\in \cC_{\eta}(\bu_1, E_{4}^u(\gamma_O(0)))$, we have
\begin{equation}\label{EqDG>}
\|d\mathbb G v\|\geq \frac12\|v\|,\end{equation}
where %we have
%\begin{enumerate}
$\bu_1$ is as in Proposition \ref{PropDL1}, and $\bar\bu_1$ and $\bar\bl_1$ are as in Proposition \ref{PropDG1}. These vectors have definite limits in the isosceles limit $Y\to 0$ and $\chi\to\infty.$ We have $dG_2 \cdot \bar \bu_2\neq 0$. 
%\item $\bu_2$ is a vector determined by $ \mathrm{span}\{\bu_2,\bu_3\}\cap (\bar\bl_3)^\perp$ and $\bar\bu_2$ is the projection of $\bu_2$ to the $(G_0,\mathsf g_{02})$-components;
%\item $\bu_3=$,\ $\bu_3\cdot\bar\bl_3\neq 0$;
%\item $ \bar\bu_2=(0_{1\times 4},0,-c_2,1,c_1-c_2),\ \bar\bl_2=(0_{1\times 4},-\beta,0,\frac{1}{L_1^L}-\beta,1)$ {\color{red} define $c_1,c_2$ here. }
%\end{enumerate}
\end{Prop}
The proof of this proposition will be given in Section \ref{non-isosceles-sec}. The proof is highly nontrivial since there is no expansion in the $\bar\bu_3$ direction (this is why we have \eqref{EqDG>}), so we have to show that the $o(\chi)$ error in $d\mathbb G$ does not spoil $\bar\bu_3$, which involves a careful analysis of the structure of the $o(\chi)$ error. 
All the above vectors can be given explicitly.  We shall be interested in controlling the plane $\mathrm{span}(\bar\bu_2,\bar\bu_3)$ as it is pushed forward by the flow, of which, one direction is used to control the $(G_0,\mathsf{g}_{01})$, i.e. the motion of be binary, and the other is used to control $(G_2,\mathsf{g}_{12})$, i.e. to make sure that the binary comes to a near triple collision configuration with $Q_2$ after application of the global map.
\subsection{The transversality condition}
When we compose $\mathcal R\mathbb L$ and $\mathbb G$ to get the Poincar\'e map $\mathbb P=\mathbb G\mathcal R\mathbb L$, we hope that $d\mathbb P$ preserves the cone $\cC_{2\eta}(\bar\bu_1, \bar\bu_2,\bar \bu_3)$ and expands each vector inside. For this purpose, we have to verify that $\mathrm{span}(\bar\bu_2,\bar \bu_3)$ satisfies \eqref{EqdRL} for a plane $P$, i.e. $\mathrm{span}\{\bar\bu_2,\bar\bu_3\}\cap \cC_\eta(E_{4}^s(\gamma_I(0)))=\{0\}. $ 

Let $\mathbf e_3',\mathbf e_4'$ be two vectors in $E_{4}^s(\gamma_I(0))$ spanning unit area, then the transversality condition is given explicitly as the nonvanishing of the determinant $\det(\bar\bu_2,\bar\bu_3, \mathbf e_3',\mathbf e_4')\neq 0$. We have the following lemma as a corollary of Lemma \ref{LmUlam} in Appendix \ref{SUlam}.  
\begin{Lm}\label{LmTrans} The determinant  $\det(\bar\bu_2,\bar\bu_3, \mathbf e_3',\mathbf e_4')$ is an analytic function
of the masses $m_1,m_2$, so it is nonvanishing for generic masses if it is so for one choice of $(m_1,m_2)$.  In particular, for such masses we have $\text{span}\{\bar{\bu}_2,\bar{\bu}_3\}\cap\mathcal{C}_\eta(E^s_{4}(\gamma_I(0)))=\{0\}$ for some $\eta>0$. 
\end{Lm}
We shall show in Appendix \ref{STrans} that the determinant is indeed nonvanishing for the masses $(m_1,m_2)=(1,1)$. This gives the residual set in Theorem \ref{ThmMain}. We point out that the verification of the nonvanishing for the masses $(m_1,m_2)=(1,1)$ in Appendix \ref{STrans} is based on a numerical calculation, since the I3BP is nonperturbative. The labor of numerics is used to solve initial and boundary value problems of ODEs.

\begin{Not}\label{NotEs}
\begin{enumerate}
%\item 

\item Denote by $\mathcal{B}^i_{O}\subset \mathbb{R}$ the set of all times $t$ where the $w_0$ or $G_0$ component of $E^u_3(\gamma^i_O(t))$ vanishes, $i=1,2$. Similarly, we define $\mathcal{B}^i_{I}$ to be that of $E^s_3(\gamma_I^i(t))$, $i=1,2$.
\item The sets $\Psi_O^i$ and $\Psi_I^i$, $i=1,2$, are, respectively,  the set of all times $t$ when the $\psi$-component of $\gamma_O^i$ and $\gamma_I^i$, $i=1,2$, have absolute value $\pi/2$.
% \item Let $\Gamma_O^i(t)$ be the solution of the variational equation of the $(G_0,g_0)$ or $(w_0,g_0)$ components along the orbit $\gamma_O^i(t)$ with $\Gamma^i_O(0)=E^s(\gamma_I^{3-i}(-\infty))$, $i=1,2$. Denote by $\mathbb{B}_O^i\subset\mathbb{R}$ the set of all the instant $t$ where the $w_0$ or $G_0$ component of $\Gamma_O^i(t)$ vanishes, $i=1,2$. Similarly, we define $\mathbb{B}_I^i$, $i=1,2$, by interchanging  ``I" and ``O" and changing $-\infty$ to $\infty$.
\end{enumerate}
\end{Not}
%The transversality conditions in the following two propositions are necessary when we compose local and global maps.
\begin{Prop}\label{PropTrans1} For the choice of masses $m_1=m_2=m_3=m_4=1$ the following transversality conditions are satisfied:
%\begin{enumerate}
%\item $\bar \bu_1\cdot \bar \bl_2\neq0$ and $\bar \bu_2\cdot \bar \bl_1\neq0$.
%Treating $E^u(\gamma^i_O(\infty))$ and $E^s(\gamma^{3-i}_I(-\infty))$ unit vectors at the same $\R^2$ plane, then we have
%\begin{equation}
%$E^u(\gamma^i_O(\infty))\pitchfork E^s(\gamma^{3-i}_I(-\infty)), \ i=1,2;$
%d\cR[\Psi_{\cS_+^1}\Phi^{\gamma^1_O}_{O_1} E^u(O_1)] \pitchfork  \Phi^{-\gamma^2_I}_{O_2} E^s(O_2),\quad d\cR[\Psi_{\cS_+^2}\Phi^{\gamma^2_O}_{O_2} E^u(O_2)] \pitchfork \Phi^{-\gamma^1_I}_{O_1} E^s(O_1).
%\end{equation}
 The sets $\mathcal{B}_O^i$ and $\mathcal{B}_O^i$, $i=1,2$, have finitely many elements, independent of $\chi$, and moreover, there exists $c>0$, independent of $\chi$, such that
$$\min\{\mathrm{dist}(\mathcal{B}_O^i, \Psi_O^i),\; \;\mathrm{dist}(\mathcal{B}_I^i, \Psi_I^i)\}>c.$$
%\end{enumerate}
%Moreover, the first component $(w_0$ or $G_0)$ of all the vectors $\Phi^{\gamma^i_O(t)}_{O_i} E^u(O_i),\ t\in (-\infty,0]$ and $\Phi^{\gamma^i_I(t)}_{O_i} [\Psi_{\cS_+^i}\Phi^{\gamma^{3-i}_O}_{O_{3-i}} E^u(O_{3-i})],\ t\in \R,\ i=1,2$ are bounded away from zero uniformly for all $t$.
%The transversality is  uniform for large $\chi$.
\end{Prop}

\subsection{Proof of Theorem \ref{ThmMain}}

In this section, we give the proof of Theorem \ref{ThmMain}. 

First, from the result of Theorem \ref{ThmMM}, we choose a point $\sf x$ such that its sojourn times satisfy Definition \ref{DefSojourn} with some $c$ sufficiently small and $C$ sufficiently large. We shall construct a genuine noncollision singularity of the F4BP close to $\sf x$. 

{\it 1. The admissible cube.}

We next introduce the notion of {\it admissible cube.} Let us fix a small  number $\eta>0$ throughout the proof. 
 
\begin{Def} We say that a three-dimensional cube  $\mathsf C$ on $\cS_-$ is \emph{admissible} at step $n$, if 

\begin{enumerate}
\item  $|Y|<\nu_n:= e^{\mu_s\sum_{i=1}^nT_i}$ on $\sf C$, where $T_i$ is as in Definition \ref{DefSojourn} for the $i$-th return;
\item  $T\mathsf C\subset \mathcal C_{\eta}(\bar\bu_1,\bar\bu_2,\bar\bu_3)$.%, so we can parametrize the admissible surface by $\ell_0$ and $G_2$ since the corresponding entries in the vectors $\bu_1$ and $\bu_2$ respectively are nonvanishing;
%\item The $\ell_0$-component contains a interval of length $\dt e^{-T_1}$ around the $\ell_0$ value of $\sx$, and the $G_2$-component contains an interval of length $\dt^2$ around 0.
\item When parametrized by the variables $\ell_0,G_0,G_2$, we have that for each $(G_0,G_2)$, the $\ell_0$-curve intersects $W^s_3$ and has length at least $\eps^2$. For each $(\ell_0,G_2)$, the curve on the $(G_0,\sf g_{01})$-plane has length between $\nu_n/2$ and $\nu_n$, and its two endpoints are on the boundary of the cone $\mathcal C_1(E_3^s(\gamma_I(0)))$. Note that the curve does not intersect $(G_0,\mathsf{g}_{01})=(0,\frac{\pi}{2})$ for $\eta$ small and  its tangent vectors satisfy the cone condition in \eqref{EqdRL}. For each $(\ell_0,G_0)$, the curve on the $(G_2,\mathsf g_{21})$-plane has $G_2$ coordinates in the interval $[-\frac12 ,\frac12] e^{\mu_s\sum_{i=1}^nT_i}.$
\end{enumerate}
\end{Def}
Item (3) is given by the transversality Lemma \ref{LmTrans}. By the sojourn time condition, we always have $e^{\mu_s\sum_{i=1}^nT_i}\ll e^{-\mu_uT_{n+1}}$, so the assumption of Proposition \ref{PropDL2} is always satisfied. 

{\it 2. The local map and the point-deleting procedure. }

We pick an initial admissible cube $\mathsf C$, then apply the local map followed by the renormalization. By Proposition \ref{PropDL2}, the cube $\mathsf C$ is strongly expanded by at least $e^{-T_1v_*/2}$.

We next examine the dynamics of the $Y$-variables. We  first focus on the F3BP, i.e. the $(w_0,\mathsf g_{01})$-components.  The Lagrange fixed point has two eigenvalues $\mu_u(>0>)\mu_s$ when restricted to the $(w_0,\mathsf g_{01})$-plane. By  the cone condition  in \eqref{EqdRL}  and by the $\lambda$-lemma, a curve on $\sf C$ with constant $\ell_0,G_2$ component will approach  $W^u_{3}(\gamma_O(0))$ exponentially fast (by $e^{T \mu_s}$) and get stretched by $e^{T \mu_u}$, where $T$ is the time span defining the local map that can be estimated as $T_1$. We keep only the part of the curve that stays within the cone $\mathcal C_1(E^s_3(\gamma_I(\tau))),\ \tau\geq 0$ and   $\mathcal C_1(E^u_3(\gamma_O(\tau))),\ \tau\leq 0$,  along the orbit. Therefore the remaining segment has length estimated as $e^{\mu_s T_1}$. %If we consider all the four $Y$-variables, the point deleting procedure in the $Y$-components gives us a segment of length estimated as $e^{\mu_s T_1}$ and its $\sf g_{01}$-component contains a neighborhood of 0. 

So we introduce a point-deleting procedure as follows
\begin{enumerate}
\item In the $X$-component, we keep only a subsegment whose sojourn time satisfies the sojourn time inequality. 
\item In the $(G_0,\mathsf{g}_{02})$-component, along the flow, we keep only the part that stays within the cone $\mathcal C_1(E^u_3(\gamma_O(\tau))),\ \tau<0$. 
\item In the $G_2$-component, we keep only the part with $G_2\in [-e^{\mu_s T_1},e^{\mu_s T_1}]$. 
\end{enumerate}

%{\it 3. The $(G_0,\sf g_{02})$ dynamics. }
%{\bf Here I have to define the vectors $\bar\bu_2,\bar\bl_2$ explicitly. }

{\it 3. The perturbation of the nonisoscelesness to the $X$-variables.}

We next estimate how the nonisoscelesness $|Y|\neq 0$ perturbs the $X$-variables. Let us denote by $X'=F(X)$ the equation of motion for the I4BP in the blowup coordinates (see \eqref{4-iso})  and by $\bar X'=F(\bar X)+\mathcal E$ the equation of the F4BP, where the error $\mathcal{E}$ is estimated as $O(Y^2)$ (see equation \eqref{blowup-ws} for the F3BP). Then denoting $\dt X=X-\bar X$, we get $(\dt X)'=DF \dt X+\mathcal E$ and we estimate the divergence of the orbits $$|\dt X(\tau)|\leq e^{M\tau}\int_0^\tau e^{-Ms} \mathcal E(s)ds\leq \bar C e^{M\tau } (Y^2(0))$$ by Gronwall, where $M$ bounds $DF$, which is estimated by $\mu_u$. By  the point deleting process, we see that $Y^2(0)$ is estimated by $e^{2\mu_s(\sum_{i=1}^{n-1} T_i)}$, inductively obtained from the last step if $n> 1$. By our assumption on the sojourn time $\tau\leq T_n$ we thus get an upper bound $|\dt X(\tau)|\leq \bar C e^{-DT_n}$ for $\tau\in[0,T_n]$ and $D=-\frac{Mc+\mu_s}{c}\gg1$ if we choose $c$ small.  The  argument is similar to that in the proof of Proposition \ref{PropLocalMM} (as in Section \ref{SSI4BP}). When $D$ is much larger than $|\mu_u|$ and $|\mu_s|$, we see that perturbation to the $X$-variables is negligible.  %So the orbit for the full four-body problem that survives the point deleting procedure remains in a $o((\log T_n)^{-1})$ neighborhood of the orbit determined by $\mathbb P^n\mathsf x$ in Theorem \ref{ThmMM}.

{\it 4. The global map and the iteration. }

After the point-deleting procedure during the local map, we obtain a cube $\bar{\mathsf C}$ on $\cS_+$ whose tangent lies in the cone $\mathcal C_\eta(\bu_1,E^u_{4}(\gamma_O(0)))$. %Moreover, since the $\mathsf g_{21}$-component  of $\bar\bu_2$ is nonvanishing, we see that the $\mathsf g_{21}$-component of the admissible surface contains a neighborhood of $\mathsf g_{21}=0$ that is independent of $\chi$. In the image of the global map, due to the strong expansion and the fact that $\bu_2$ has nonvanishing $G_2$ and $\mathsf g_{21}$ components, we get that in the image $\mathbb P(\mathsf C)$, we can choose a piece of admissible surface. 
We next apply the global map. We use $d\mathbb G$ in Proposition \ref{PropDG2} to control the difference between an orbit issued from $\bar{\mathsf C}$ and an orbit of the I4BP. 
Estimate \eqref{EqDG>} in Proposition \ref{PropDG2} implies that $\mathbb G(\bar{\mathsf C})$, when projected to the plane span$\{\bar\bu_2,\bar \bu_3\}$,  covers a disk of radius $\frac14e^{\mu_s T_1}$ centered at zero. We perform another point deleting procedure so that we get another admissible cube from $\mathbb G(\bar{\mathsf C})$. We remark that though the global map has no expansion in the $\bar\bu_3$ direction,  the transversality in Lemma \ref{LmTrans} and the cone condition in  \eqref{EqdRL} show that 
it will be expanded by the local map. 

Thus, we can iterate the procedure for infinitely many steps. The point deleting procedure gives rise to a point-deleting procedure on the initial admissible  cube $\mathsf C$, hence after infinitely many steps,  we get a Cantor set of initial conditions.

{\it 5. No binary collisions. }

To show that there is no binary collision, we first invoke the following lemma whose proof is given at the end of Section \ref{SSTrans}. 
\begin{Lm}\label{LmOscG0}
There exist $C$ and $\nu_0$ independent of $\chi$ such that for orbits between the sections $\cS_+^1$ and $\cS_-^2$, if on the section $\cS_+^1$, 
$$|G_0|=\nu_1\in(0,\nu_0], \quad |g_0-g^L_1-\frac{\pi}{2}|,\;|g^L_1-g^L_2-\pi|, \; |G^L_2|\leqslant C|G_0|,$$
and on the section $\cS^2_-$, 
$$|G_0|=\nu_2\in(0,\nu_0],\quad  |g_0-g^R_2-\frac{\pi}{2}|,\;|g^R_1-g^R_2-\pi|, \; |G_2^R|\leqslant C|G_0|,$$
 then we have $2\nu_0>|G_0(t)|>\frac{1}{2}\min\{\nu_1,\nu_2\}$.
\end{Lm}
 
Lemma \ref{LmOscG0} shows that along the piece of orbit defining the global map, there is no double collision since $G_0$ does not vanish. For the local map, it may happen that $G_0$ vanishes at some isolated points whose number does not depend on $\chi$ near triple collision. However, when this happens, by Proposition \ref{PropTrans1}, we have $\psi\neq \pm\pi/2$, so the binary collision actually does not occur.

The proof is now complete. 
\qed
% Now we can repeat the above argument by replacing $Q_1$-$Q_3$-$Q_4$ by $Q_1$-$Q_3$-$Q_4$. The procedure can be iterated for infinitely many steps. In the limit, we will result in a Cantor set on $S_0$ whose points can be iterated for infinitely many steps. 

 \section{The  Hamiltonians and the isosceles three-body problem}\label{SJacobi}
In this section, we first write down the Hamiltonian of the F4BP and study the transform between the left  and right Jacobi coordinates in the first two subsections.  In Section \ref{SSIso}, we write down the equations of motion of the I3BP in blowup coordinates with double collision regularized. In Sections \ref{SSgammaI} and \ref{SSNondeg}, we give the proof of Proposition \ref{ThmGammaI} on the existence of $\gamma_I$ as well as Lemma \ref{LmNondeg}.
 \subsection{The Hamiltonian systems in Jacobi coordinates}
The original Hamiltonian has the form
$$H(P,Q)=\frac{P_1^2}{2m_1}+\frac{P_2^2}{2m_2}+\frac{P_3^2}{2}+\frac{P_4^2}{2}-\frac{m_1m_2}{|Q_1-Q_2|}-\sum_{i=1,2}^{j=3,4}\frac{m_i}{|Q_i-Q_j|}-\frac{1}{|Q_3-Q_4|}.$$
In the (left) Jacobi coordinates, defined in Section \ref{subsection-jacobi}, the Hamiltonian reads
\begin{equation}
\begin{aligned}
&H_{\cdot:\cdot}(\bx,\bp)=\left(\frac{\bp_1^2}{2M_1}+\frac{\bp_0^2}{2M_0} -\frac{m_1}{|\bx_1-\frac{\bx_0}{2}|}-\frac{m_1}{|\bx_1+\frac{\bx_0}{2}|}-\frac{1}{|\bx_0|}\right) \\
&+\left(\frac{\bp_2^2}{2M_2}  -\frac{m_1m_2}{|\bx_2-\frac{2 \bx_1}{m_1+2} |}-\frac{m_2}{|\bx_2+\frac{m_1\bx_1}{m_1+2}-\frac{\bx_0}{2}|}-\frac{m_2}{|\bx_2+\frac{m_1\bx_1}{m_1+2}+\frac{\bx_0}{2}|}\right)
\end{aligned}
\end{equation}
where
$M_0=\frac12,\  M_1=\frac{2m_1}{m_1+2},\  M_2=\frac{m_2(m_1+2)}{m_1+m_2+2}.$

 We will work with two different kinds of Hamiltonian. When the binary $Q_3$-$Q_4$ is close to $Q_1$ (local map), we write the Hamiltonian as
\begin{equation}\label{EqHamLoc}
\begin{aligned}
H_{\cdot:\cdot}(\bx,\bp)&=H_\therefore+H_2\\
H_\therefore&=\frac{\bp_1^2}{2M_1}+\frac{\bp_0^2}{2M_0} -\frac{m_1}{|\bx_1-\frac{\bx_0}{2}|}-\frac{m_1}{|\bx_1+\frac{\bx_0}{2}|}-\frac{1}{|\bx_0|}:=\frac{\bp_1^2}{2M_1}+\frac{\bp_0^2}{2M_0} +V,\\
 H_2&=\frac{\bp_2^2}{2M_2} -\frac{k_2}{|\bx_2|}+U_2, \quad k_2=(m_1+2)m_2,\\
   U_2&=\frac{k_2}{|\bx_2|}-\frac{m_1m_2}{|\bx_2-\frac{2 \bx_1}{m_1+2} |}-\frac{m_2}{|\bx_2+\frac{m_1\bx_1}{m_1+2}-\frac{\bx_0}{2}|}-\frac{m_2}{|\bx_2+\frac{m_1\bx_1}{m_1+2}+\frac{\bx_0}{2}|}.
\end{aligned}
\end{equation}
Note that the subsystem $H_\therefore$ is a three-body problem in Jacobi coordinates and $H_2$ is a Kepler problem perturbed by $U_2$.

When the binary is far away from both $Q_1$ and $Q_2$ (global map), we use the following form of Hamiltonian system, which is a perturbation of three Kepler problems
\begin{equation} \label{EqHamGlob}
\begin{aligned}
H_{\cdot:\cdot}(\bx,\bp)&=\left(\frac{\bp_0^2}{2M_0}-\frac{k_0}{|\bx_0|}\right)+\left(\frac{\bp_1^2}{2M_1}-\frac{k_1}{|\bx_1|}\right)+\left(\frac{\bp_2^2}{2M_2}-\frac{k_2}{|\bx_2|}\right) +U_{01}+U_2\\
k_0&=1,\ k_1=2m_1,\quad\\
U_{01} &=\left(\frac{k_1}{|\bx_1|}-\frac{m_1}{|\bx_1-\frac{\bx_0}{2}|}-\frac{m_1}{|\bx_1+\frac{\bx_0}{2}|}\right).
\end{aligned}
\end{equation}
We denote  the energy $E_i=\frac{|\bp_i|^2}{2M_i}-\frac{k_i}{|\bx_i|}$ for the Kepler problem $(\bx_i,\bp_i)$,  $ i=0,1,2.$

We remark that for the Hamiltonians written  in the right Jacobi coordinates, that is when the binary $Q_3$-$Q_4$ is close to $Q_2$, in the definition of $M_1,k_1,M_2,k_2$, we should permute the subscripts $1$ and $2$. If necessary, we will put a superscript $L$ or $R$ to distinguish the quantities in the left or right Jacobi coordinates.
\subsection{Transition from the left to the right}\label{STransition}
When the binary $Q_3$-$Q_4$ is closer to $Q_i$ than $Q_{3-i}$, we use the Jacobi coordinates \eqref{EqJacobi} to treat $Q_i$-$Q_3$-$Q_4$ as a three-body problem and $Q_{3-i}$ as a perturbation to the three body problem, $i=1,2$.   We now give the linear transform changing the left Jacobi coordinates to the right Jacobi coordinates.
Denote $(\bx_0,\bx_1,\bx_2)^t_L=\cL (q_3,q_1,q_2)^t$ and $(\bx_0,\bx_1,\bx_2)^t_R=\cR (q_3,q_1,q_2)^t$. Then we get $(\bx_0,\bx_1,\bx_2)^t_L=\cL \cR^{-1} (\bx_0,\bx_1,\bx_2)^t_R$ and  $(\bx_0,\bx_1,\bx_2)^t_R=\cR \cL^{-1} (\bx_0,\bx_1,\bx_2)^t_L.$
We have
$$\cL=\left[\begin{array}{ccc}
1&0&0\\
-\frac12&1&0\\
-\frac{1}{m_1+2}&-\frac{m_1}{m_1+2}&1
\end{array}\right],\quad \cR=\left[\begin{array}{ccc}
1&0&0\\
-\frac{1}{m_2+2}&1&-\frac{m_2}{m_2+2}\\
-\frac12&0&1
\end{array}\right].$$
This gives
$$\cL\cR^{-1}=\left[\begin{array}{ccc}
1&0&0\\
0&1&\frac{m_2}{m_2+2}\\
0&-\frac{m_1}{m_1+2}&1-\frac{m_1m_2}{(m_1+2)(m_2+2)}
\end{array}\right],\quad \cR\cL^{-1}=\left[\begin{array}{ccc}
1&0&0\\
0&1-\frac{m_1m_2}{(m_1+2)(m_2+2)}&-\frac{m_2}{m_2+2}\\
0&\frac{m_1}{m_1+2}&1\\
\end{array}\right].$$
We then denote $\frac{m_2}{m_2+2}=\frac{1}{2}M_2^R $ and $\frac{m_1}{m_1+2}=\frac12 M_1^L$.
Explicitly, we have
\begin{equation}\label{EqLR}
\begin{cases}
\bx_1^R&=(1-\frac{M_1^LM_2^R}{4} )\bx_1^L-\frac{M_2^R}{2} \bx_2^L,\\
\bx_2^R&=\frac{ M_1^L}{2}\bx_1^L+ \bx_2^L,\\
\end{cases}\quad
\begin{cases}
\bp_1^R&=\bp_1^L-\frac12 M_1^L\bp_2^L,\\
\bp_2^R&=\frac12 M_2^R\bp_1^L+(1-\frac{1}{4}M_1^LM_2^R )\bp_2^L
\end{cases}
\end{equation}
where the transform of the momenta is obtained by taking the transpose inverse of that of the positions.
\subsection{ The isosceles 3-body problem} \label{SSIso}Without loss of generality, we assume the particle $Q_2$ is infinitely far away and consider the I3BP $Q_1$-$Q_3$-$Q_4$.  The system has  two degrees of freedom and in the blowup coordinates $(r,\psi, v, w)$ (c.f. equation \eqref{EqBlowup}), the  equations of motion read
\begin{equation}\label{triple-iso}
\begin{cases}
r'&=rv\\
v'&=w^2+\frac{1}{2}v^2+\bar V(\psi),\\
\psi'&=w,\\
w'&=-\frac{1}{2} vw-\partial_\psi \bar V(\psi),
\end{cases}
\end{equation}
where \begin{equation}\label{iso-potential}\bar V(\psi)=-\frac{1}{\sqrt{2}\cos\psi}-\frac{4m_1}{\sqrt{2\cos^2\psi+\frac{2(2+m_1)}{m_1}\sin^2\psi}}.\end{equation}
and the energy relation becomes
$rE=\frac{1}{2}v^2+\frac{1}{2}w^2+\bar V(\psi).$
The potential $\bar V$ blows up when $\psi=\pm\pi/2$, corresponding to the double collision between the binary~$Q_3$-$Q_4$. We regularize the double collision  by introducing $\hat w=w\cos \psi $ as well as a new change of time $\frac{d\tau}{d\hat \tau}=\cos\psi$. The equations of motion become
\begin{equation}\label{EqSundman}
\begin{cases}
\frac{dr}{d\hat \tau}&=rv\cos \psi ,\\
\frac{dv}{d\hat \tau}&=\left(2rE-\frac12v^2-\bar V(\psi)) \right)\cos\psi\\
\frac{d\psi}{d\hat \tau}&= \hat w\\
\frac{d\hat w}{d\hat \tau}&=-\frac12 v\hat w\cos\psi -\bar V'(\psi)\cos^2\psi-\sin\psi\cos\psi(2rE-v^2-2\bar V(\psi)),
\end{cases}
\end{equation}
 with the energy relation
 $rE=\frac{\hat w^2}{2\cos^2\psi}+\frac{1}{2} v^2+\bar V(\psi). $ In these new coordinates, the system is an analytic vector field defined on $[0,\infty)\times \R\times [-\pi/2,\pi/2]\times \R$ and the double collision becomes an elastic collision so that the singularity is regularized. Recall that we denote by $\cM_0$  the collision manifold $\{r=0\}$. The dynamics on $\mathcal{M}_0$, as stated in  Theorem \ref{collision-m1} and \ref{collision-m2},  follows by analyzing these equations of motion.
 % Let $P_1$ be the Lagrange fixed point on the right bottom of $\mathcal{M}_0$, that is
  %\[P_1=(0,-v_*,\psi_*,0), \; v_*>0, \psi_*>0.\]
 %From Theorem \ref{collision-m1}, we know that  $P_1$ is a saddle point on $\mathcal{M}_0$ with $1$-dimensional stable and unstable manifolds $\gamma_-$ and $\gamma_O$. The stable manifold $\gamma_-$ comes from the lower arm with $\psi=\frac{\pi}{2}$, and the unstable manifold $\gamma_O$ runs up to the upper arm with $\psi=\frac{\pi}{2}$. So for orbits shadowing these two, $Q_1$ would approach the binary $Q_3$-$Q_4$ from the left side and then leave them from the left side.

% Now we consider $P_1$ as a fixed point in the isosceles system \eqref{triple-iso},  which we denote as $\cM_{iso}$. One checks  (see \cite[Section 3]{D}) that $P_1$ is a saddle point with a $2$-dimensional stable manifold, which we denote as $W_{iso}^s$, and a $1$-dimensional unstable manifold, which lies in $\cM_0$.  Clearly $\gamma_-\subset W_{iso}^s$.
\subsection{Boundary value problem for the I3BP}\label{SSgammaI}
In this section, we show the existence of $\gamma_I.$ Without loss of generality, we consider the I3BP $Q_2$-$Q_3$-$Q_4$.
%\begin{Lm}\label{L0psi}
%On the section $S_-^R=\{r=\frac{1}{\epsilon}\}$, we have that if $L_0\sim1$ and $r_0\sim1$,
%\[\frac{\partial L_0}{\partial \psi}\sim\epsilon^{-1}.\]
%\end{Lm}
%\begin{proof}
%We have the following differential relation,
%\begin{equation*}
%\begin{aligned}
%d\psi&=\frac{\sqrt M_2r_0 dr_2-\sqrt M_0 r_2dr_0}{r^2}=\frac{1}{r}(\sqrt{\frac{M_2}{M_0}}\cos\psi dr_2-\sqrt{\frac{M_0}{M_2}}\sin\psi dr_0),\\
%dr_0&=\frac{\partial r_0}{\partial L_0}dL_0+\frac{\partial r_0}{\partial \ell_0}d \ell_0,\quad
%dr_2=\frac{\partial r_2}{\partial L_2}dL_2+\frac{\partial r_2}{\partial \ell_2}d \ell_2.
%\end{aligned}
%\end{equation*}
%By Lemma \ref{Lm-A1}, we have  $\frac{\partial r_0}{\partial L_0}=\frac{2r_0}{L_0}$. Then by implicit function theorem we obtain that
%$$\frac{\partial L_0}{\partial \psi}=-\frac{\sqrt{M_0}L_0}{\sin2\psi}.$$
%If $L_0\sim1$ and $r_0\sim1$, we have $\psi\sim\epsilon$, which gives the assertion of the lemma.
%\end{proof}
\begin{Lm}\label{energy-partition-infinity} Fix the energy level $E<0$ in the I3BP,  then there exists $C>0$ such that  for any $a\in(-\infty, E-C\epsilon)$ there exists $\sx\in W^s_{\therefore}(O_2)\cap \cS_-^R$ such that
$E_0(\sx)=a$, where $E_0(\sx)$ is the energy of the binary $Q_3$-$Q_4$ at the phase point $\sx$. In particular, for any $a\in(-\infty,E-C\eps )$ there exists an orbit $\gamma_a$ on $W^s_{\therefore}(O_2)$ such that
$E_0(\gamma_a(-\infty))=a$ and $r(-\infty)=\infty$. Moreover, this $\gamma_a$ is locally unique on each energy level $E$.
\end{Lm}
\begin{proof} We shall shoot orbits from a small neighborhood of $O_2$ and trace the flow backward as $\tau\to-\infty$. Equivalently, readers can also consider shooting orbits from a neighborhood of the Lagrange fixed point in the upper half space and trace the flow forward. On the stable manifold  $W^s_{\therefore}(O_2)$ there exist orbits which eject backwardly from $O_2$ with arbitrarily small $r_0$ when crossing the section $\cS_-^R$, thus the momentum $\bp$ can be arbitrarily large speed.  For any $B>0$ large enough, let~$\sx_1\subset W^{s}_{\therefore}(O_2)$ be a  point in a small neighborhood of $O_2$ such that the corresponding orbit crosses the section $\cS^R_-$ with kinetic energy  of $Q_2$ greater than $B$.
It is proved in~\cite{D} that on $W^s_{\therefore}(O_2)$, there exists an orbit which is heteroclinic from the Lagrange fixed point on the upper half space to $O_2$  with fixed $\psi$ value. It is the orbit with homothetic equilateral triangle configuration and the variable $r$ is always bounded.  Let $\sx_2\in W^s_{\therefore}(O_2)$ be a point on this orbit in a small neighborhood of $O_2$.

  Let $\zeta_0(s)$, $0\leqslant s\leqslant1$, be a small smooth curve in $W^s_{\therefore}(O_2)$ with $\zeta_0(0)=\sx_1$ and $\zeta_0(1)=\sx_2$. Let
$s_*$ be the largest number such that the orbit corresponding to the phase point $\zeta(s)$ crosses the section $\cS^R_-$ for all $0\leqslant s\leqslant s_*$.  Then we have $0<s_*<1$.
Let $\zeta(s)\subset W^s_{\therefore}(O_2)$, $s\in[0,s_*]$ be the image of $\zeta_0([0,s_*])$ on the section $\cS_-^R$. Let $\Xi_2(s)$ be the kinetic energy of $Q_2$ at the point $\zeta(s)$. We then have that
$\Xi_2(s_*)=O(\epsilon)$ since otherwise due to the continuous dependence on initial data,  it contradicts the maximality of $s_*$.   Clearly, $\Xi_2(s)$ varies continuously on the curve $\zeta([0,s_*])$. Therefore, $[0,B]\subset\Xi_2([0,s_*])$. Let $E_0(s)$ be the energy of the binary $Q_3$-$Q_4$ at the point $\zeta(s)$. From energy conservation we have that
$E_0(s)+\Xi_2(s)=E+O(\epsilon).$ Therefore there exists $C>0$ such that $(-B+E+C\epsilon, E+C\epsilon]\subset E_0([0,s_*])$.

Finally, for the local uniqueness, we note that the stable manifold $W_\therefore(O_2)$ is an analytic manifold, and so are the energy level set and the section $\cS_-$. When restricted to the section $\cS_-$, we can realize the intersection point $\gamma_a\cap \cS_-$ as the zero of an analytic function, which is necessarily isolated.
\end{proof}
\subsection{Proof of Lemma \ref{LmNondeg}}\label{SSNondeg}
With the above lemma, we now are ready to prove Lemma \ref{LmNondeg}.
\begin{proof}[Proof of Lemma \ref{LmNondeg}]

By Proposition \ref{PropDL1} and Proposition \ref{PropDG1}, a vector (in coordinates $(L_0,\ell_0, L_2,\ell_2)$) of the form \mbox{$(1,O(\eta),O(\eta),O(\eta))$}  is mapped to a vector of the form $O(\epsilon^2\chi,\chi, 1,\chi)$ (in coordinates $(L_0,\ell_0, L_1,\ell_1)$) by $d\mathbb G$. The estimate is uniform on the curve $\xi$. So when projected to $(L_0, \ell_0)$, the image on the section $\cS^R_-$ is almost tangent to $(\epsilon^2,1)$ at every point and has length of order $\chi\cdot \frac{1}{\epsilon\chi}=\epsilon^{-1}\gg1$. Therefore the resulting curve winds around the $\ell_0$-circle many times.

We use the variables $(L_0,\ell_0)$ to describe the section $\cS_{-}^R$ when restricted to the energy level $E$ of the I3BP. The variables $(L_0,\ell_0)$ are defined on a  cylinder with $L_0\in[0, L_*]$ and  $L_*=\sqrt{\frac{M_0k_0^2}{-2E}}$.
 Then by Lemma \ref{energy-partition-infinity}, the intersection of $W_{\therefore}^s(O_2)$ and $\cS_-^R$ contains at least a continuous curve which connects the two parts on the cylinder with $L_0= \sqrt{\frac{M_0k_0^2}{-2E+O(\epsilon)}}$ and $L_0=0$. So it necessarily intersects topologically transversally a homologically nontrivial circle with $0<L_0<L_*-c$ for certain $c>0$. This proves the assertion of the lemma.
\end{proof}

\section{Some preliminary estimates of the  isosceles four-body problem}\label{SBlowupMM}
In this section,  we study the local dynamics of the I4BP. Without loss of generality, we assume that $Q_1$-$Q_3$-$Q_4$ is close to triple collision and $Q_2$ is distance $\chi$ away. 
  In Section \ref{SSI4BP}, we study the I4BP near triple collision and give the proof of Proposition \ref{PropLocalMM} on the dynamics of the local map. In Section \ref{SSRenorm}, we give the proof of Proposition \ref{PropRenorm} and
we prove Proposition \ref{PropDevaneyOrbit}  in Section \ref{SSPartition}.

 \subsection{The I4BP near triple collision}\label{SSI4BP}
 For the I4BP, we need to introduce the perturbation from $Q_2$.
 We convert $(\bx_1,\bp_1, \bx_0,\bp_0)$ into blowup coordinates $(r,v,w, \psi)$ and $(\bx_2,\bp_2)$ into Delaunay coordinates $(L_2,\ell_2)$. In the Hamiltonian \eqref{EqHamLoc}, we have $U_2=O\left(\frac{r}{|x_2|^3}\right),$
where the estimate is obtained by assuming $\bx_0,\bx_1$ are of order $r\ $ and $\bx_2$ is of order $\chi\gg1$. Converting the Keplerian Hamiltonian $\frac{|\bp_2|^2}{2M_2}-\frac{k_2}{|\bx_2|}$ into Delaunay coordinates as $\frac{M_2k_2^2}{L_2^2}$, we obtain $  H_{\cdot:\cdot}=H_\therefore+\frac{M_2k_2^2}{2L_2^2}+U_2. $
 The equations of motion are as follows with double collision regularized as in Section \ref{SSIso} (with the same notation as Section \ref{SSIso})
 %\begin{equation}\label{triple-iso}
%\begin{cases}
%r'&=rv\\
%v'&=w^2+\frac{1}{2}v^2+\bar V(\psi)+r\big(\bx\cdot\partial_{\bx} \hat V_{012}(\bx, \bx_2)\big),\\
%\psi'&=w,\\
%w'&=-\frac{1}{2} vw-\partial_\psi \bar V(\psi)+r^2(\partial_{\bx}\hat V_2\times \bs),
%\end{cases}
%\begin{cases}
%L_2'&=-r^{3/2}\frac{\partial}{\partial\ell_2}\hat V_{012} \\
%\ell_2'&=-\frac{M_2k_2^2}{L_2^3}+r^{3/2}\frac{\partial}{\partial L_2}\hat V_{012}
%\end{cases}
%\end{equation}
%It becomes the following with the double collision regularized
 \begin{equation}\label{4-iso}
\begin{cases}
\frac{dr}{d\hat \tau}&=\cos \psi rv,\\
\frac{dv}{d\hat \tau}&=(2rE-\frac12v^2-\bar V(\psi)) \cos\psi+\cos \psi r\big(\bx\cdot\partial_{\bx} U_2\big),\\
\frac{d\psi}{d\hat \tau}&= \hat w,\\
\frac{d\hat w}{d\hat \tau}&= -\frac12 v\hat w\cos\psi -\bar V'(\psi)\cos^2\psi-\sin\psi\cos\psi(2rE-v^2-2\bar V(\psi))\\
&\quad +\cos^2\psi r(\partial_{\bx_0}U_2\cdot \bx_1-\partial_{\bx_1}U_2\cdot \bx_0),\\
\frac{dL_2}{d\hat{\tau}}&=-\cos\psi r^{3/2}\frac{\partial}{\partial\ell_2}U_2, \\
\frac{d\ell_2}{d\hat \tau}&=\cos\psi r^{3/2}(-\frac{M_2k_2^2}{L_2^3}+\frac{\partial}{\partial L_2}U_2).
\end{cases}
\end{equation}
Compared to \eqref{EqSundman}, the equations of motion here are perturbed by $U_2$. 
Now we are ready to prove Proposition \ref{PropLocalMM}.
\begin{proof}[Proof of Proposition \ref{PropLocalMM}]
We first deal with the I3BP $Q_1$-$Q_3$-$Q_4$.
  Let $O_1$ be the Lagrange fixed point and $\mathscr O$ a small neighborhood of $O_1$. For a fixed energy level $E<0$ of the I3BP, the intersection of $W_\therefore^s(O_1)$ and $\cS_-^1$ is a smooth curve $\Gamma$. Let $\sx_*\in \Gamma$ be the point corresponding to $\gamma_I$. 
 Let $\zeta$ be a smooth curve on $\cS_-^1$ that is transversal to $\Gamma$ at $\sx_*$ and on the correct side. Note that we can find orbit that follows the unstable manifold of $O_1$ as far as we wish by choosing the initial condition on $\zeta$ sufficiently close to $W^s_\therefore(O_1)$.  % Let $\Gamma'_{\sx_*,\delta_3}\subset\Gamma$ be on this side of $W^s_{\therefore}$ and for $\sx\in\Gamma'_{\delta_3}$, $|\sx-\sx_*|<\delta_3$. Then if  $\delta_3$ is small enough, the orbit starting from $\Gamma_{\sx_*,\delta_3}'$ will intersect with the section $\cS^1_+$ inside the set $\mathscr T_{O_1,\delta_2,\delta_2'}$. Clearly, $\delta_3$ can be chosen uniformly for all $\sx_*\in\mathcal{K}$.
%Let define $\mathscr{I}_{\mathcal{K},\delta}$ be such that $\sx\in\mathscr{I}_{\mathcal{K},\delta}$ if $\sx\in \Gamma'_{x_*,\delta_3}$ for some $\sx_*\in\mathcal{K}$,
%$$|\sx-\sx_*|\leqslant\delta\quad \text{and}\quad |\psi-\psi(\sx_*)|\geqslant\frac{\delta}2.$$
% Assume $\delta>0$ is small enough. Then  for $\sx\in\mathscr{I}_{\mathcal{K},\delta}$, we have that   if $\sx(0)=\sx$,  then there exists $\tau_\sx\sim|\log\delta|$ such that
  %\begin{equation}\label{path1}\sx(\tau_\sx)\in \Gamma_{O,\delta_2,\delta_2'}, \quad \sx(\tau)\in N_{P_1,\delta_1},\;\tau\in[C_1, \tau_\sx-C].\end{equation}
  Let $\sx\in \zeta$ be an initial condition whose distance to $\Gamma$ is $\hat\dt\ll 1$. Then following the orbit $\gamma_I$, the orbit $\gamma_\sx$ with initial condition $\sx$
  will follow $\gamma_I$ then $\gamma_O$, and stay within the neighborhood $\mathscr O$ for time of order $\log\hat \dt^{-1}$, but the travel times for the orbit $\gamma_\sx$ from
  $\cS_-^1$ to entering $\mathscr O$ and from leaving $\mathscr O$ to hitting $\cS_+^1$ are both of order $1$.  The set $\mathscr O$ is chosen small so that within $\mathscr O$, we
  have $v$ close to $v_*<0$, which is the value of $v$ at $O_1$.

 From the equation $r'=rv$, we get  $r(\tau)=r(0) e^{\int_0^\tau v(s)ds}.$
 Therefore,  for given $\dt>0$ sufficiently small, by choosing $\hat \dt$ accordingly, we can arrange that when arriving at the section $\cS_+^1$, we have $r\in (\dt^{3/2},\dt)$.
% \begin{equation}\label{r-upper-lower} r\in (\dt^{3/2},\dt)%\leqslant r(0)C_*\delta^{\frac{4v_*}{5}}\text{ and } r(\tau)\geqslant r(0)C_*^{-1}\delta^{\frac{6}{5}v_*},\quad \tau\in[0,\tau_\sx].
 %\end{equation}
 This proves the assertion for the I3BP.

Now we consider the I4BP \eqref{4-iso} by taking into account the perturbation from $Q_2$. From the expression of $U_2$, we estimate the perturbation to the $(r,v,\psi,w)$-equations
 as $r\big(\bx\cdot\partial_{\bx} U_2(\bx, \bx_2)\big)=O(r^2/\chi^3)$ and $r^2(\partial_{\bx_0}U_2\cdot \bx_1-\partial_{\bx_1}U_2\cdot \bx_0)=O(r^2/\chi^3)$. We abbreviate the equations of motion for
 $\sx_1=(r_1,v_1,\psi_1,\hat w_1)$ in the I3BP \eqref{EqSundman} as $\sx_1'=F(\sx_1)$ and that for the I4BP (\ref{triple-iso}) as $\sx_2'=F(\sx_2)+\mathcal{E}(\sx_2,L_2,\ell_2)$ where
 $\mathcal{E}$ includes the error term from $\nabla U_2$ estimated as $r^2/\chi^3$. Taking the difference we obtain the equations for $\dt\sx=\sx_1-\sx_2$
$$\dt\sx'(\tau)=\left(\int_0^1DF(s\sx_1(\tau)+(1-s)\sx_2(\tau))ds\right) \dt\sx(\tau)+\mathcal{E}. $$
So the error of the solution is estimated as $$|\dt\sx(T)|\leq \int_0^T e^{M(\tau-s)}|\mathcal E(\sx(s),L_2(s),\ell_2(s))|\,ds\leq C_\eps e^{MT}/\chi^3. $$
 where $M$ is the upper-bound for $DF$ and $T$ is the travel time for the local map. We obtain $C_\eps$ by integrating $r$ between the two sections ($r$ decays exponentially
  for most of the time). Similarly, we get that the change of $L_2$ is estimated as
$ L_2(T)-L_2(0)=C_\eps /\chi^3. $
For small enough $\dt$ in the statement, the travel time $T$ is mostly spent in a small neighborhood of the Lagrange fixed point, so we can estimate $M$ by the eigenvalues of $O_1$ and $T$ by $C\log d^{-1}$, where $d$ is the distance of the initial condition to $W^s_\therefore(O_1)$. We shall choose $C\chi^{-M}<d< C^{-1} \dt^{1/v_*}$ for some constant $C>1$ so that we have $e^{MT}<\chi$ and $e^{-v_*T}r(0)<\dt$. We then choose $D=\max\{2v_*M,1\}$ in the statement, thus the perturbation of $Q_2$ to I3BP is estimated by $C/\chi^2\leq \dt^2\ll \dt^{3/2}$ so it is negligible.
\end{proof}
%We next linearize \eqref{triple-iso} at the Lagrange fixed point $O$. Since the collision manifold has two dimensions, we use the coordinates $(\psi, w)$ and eliminate $v$ using the energy relation. Linearizing the $\psi'$ and $w'$ equations from \eqref{triple-iso}, we get $\left[\begin{array}{cc}
%0&1\\
%-\bar V''&-\frac12v
%\end{array}\right]$ evaluated at $O$. Thus we get two eigenvalues $\lambda_\pm=-\frac{v}{4}\pm \frac12\sqrt{\frac14 v^2-4\bar V''(\psi)}$. Evaluating the eigenvalues at $O$, we get that $|\lambda_-|<\lambda_+<2|v_*|$. This implies that sufficiently close to $O$, we  have $M$ is close to $\lambda_+$ and the growth of $e^{Mt}$
%We finally estimate the time span in the original time scale $t$.  Since $dt=r^{3/2}d\tau$,  then together with \eqref{r-int}, we  have
 %\[t_\sx=r_0^{3/2}\int_0^{\tau_\sx}e^{\frac{3}{2}\int_0^{\tau}v(s)ds}d\tau,\quad r_0=r(\sx).\]
%From the exponential decay of $r$ within the set $\mathscr O$ we get that there exists $C_3>0$ such that
 %\[\int_0^{\tau_\sx}e^{\frac{3}{2}\int_0^\tau v(s)ds}ds\leqslant C_3,\quad \forall \sx\in\bar\zeta.\]
 %Therefore, in the original time scale $t_\sx\leqslant C_3 r_0^{3/2}.$ Hence we complete the proof of Proposition \ref{PropLocalMM}.

 \subsection{The outcome of the renormalization on the section $\cS_+$}\label{SSRenorm}In this subsection, we give the proof of Proposition  \ref{PropRenorm}.

  \begin{proof}[Proof of Proposition \ref{PropRenorm}]
From the definition of $v$ and $w$, we have that
\[\begin{cases}
M^{-1/2}_0r^{1/2}R_0\cos\psi+M_1^{-1/2}r^{1/2}R_1\sin\psi=v\\
-M_0^{-1/2}r^{1/2}R_0\sin\psi+M_1^{-1/2}r^{1/2}R_1\cos\psi=w.
\end{cases}\]
Then
\[\begin{cases}R_0M_0^{-1/2}=r^{-1/2}(v\cos\psi-w\sin\psi),\\
R_1M_1^{-1/2}=r^{-1/2}(v\sin\psi+w\cos\psi).\end{cases}
\]
Hence we have $$E_0=\frac{|\bp_0|^2}{2M_0}-\frac{k_0}{r_0}=\frac{1}{2}r^{-1}(v\cos\psi-w\sin\psi)^2-\frac{1}{r_0}$$ as energy of  $\bx_0$, and $$E_1=\frac{|\bp_1|^2}{2M_1}-\frac{k_1}{r_1}=\frac{1}{2}r^{-1}(v\sin\psi+w\cos\psi)^2-\frac{2m_1}{r_1}$$ as the energy of $\bx_1$.

For any sufficiently small $\epsilon>0$, when restricted to the section $\cS_+^R$, we have $v=\eps^{-1/2}$. Then, by the energy relation, we have
 $\frac{\epsilon}{2}\leqslant\psi+\frac{\pi}{2}\leqslant\epsilon$, 
$\frac{-3}{\eps r}\leqslant E_0\leqslant \frac{-1}{3\eps r}$, and
$\frac{1}{3\eps r}\leqslant E_1\leqslant \frac{3}{\eps r}.$
 Next we consider the section $\cS_+^R=\{v=\eps^{-1/2}\}$ in the I4BP. Suppose we have $|\bx_2|\geqslant\chi$ and $C^{-1}\leqslant |\bp_2|\leqslant C$ on the section $\cS_+^R$, %Define the set $\mathscr T_{O_1}\subset S_+^1$ in this case as
%\begin{equation}
%???\mathscr T_{O_1}=\{x\in \pi^{-1}(\mathscr T_{O_1}):|\bx_2|\geqslant\chi,\; C^{-1}\leqslant |\bp_2|\leqslant C.\},
%\end{equation}
 where $C>1$ is a fixed  constant. %and $\pi(\cdot)$ denotes the projection from phase space of  the 4-body problem to that of the 3-body problem.
We then apply the renormalization map to the set $\cS_+^R$.  The normalization factor $\lambda$ is of order $\frac{ \epsilon^{-1}}{r}$. After the renormalization, we have that there exists $C_2>1$ such that
$$-C_2\leqslant E_0\leqslant- C_2^{-1},\;C_2^{-1}\leqslant E_1 \leqslant C_2,\;  0<E_2\leqslant C\lambda^{-1}.$$
This gives the relation \eqref{EqRenorm}. 
\end{proof}

%This complete the proof of Proposition \ref{PropRenorm}.
\subsection{Energy partition}\label{SSPartition}
We now prove Proposition \ref{PropDevaneyOrbit}.
\begin{proof}[Proof of Proposition \ref{PropDevaneyOrbit}]
We are in the isosceles setting. Suppose the triple $Q_1$-$Q_3$-$Q_4$ is ejected from the near triple collision with $|\bp_0|$ and $|\bp_1|$ large and $|\bp_2|\sim1$. Suppose the distance between $Q_1$ and $Q_2$ is greater than $\chi$. On the section $\cS^L_+$, we perform a renormalization which rescales $\bp_0$ and $\bp_1$ to  size of order $O(1)$ and $|\bp_2|=O(\beta)$, where $\beta=\lambda^{-1/2}$ and $\lambda$ is the renormalization factor. Then the distance between $Q_1$ and $Q_2$ is at least $\lambda\chi$. Now the triple $Q_1$-$Q_3$-$Q_4$ has total energy $O(\beta^2)\ll 1$. Therefore we have 
\begin{equation}\label{EqEnergyPartition1}E^L_1+E_0+U_{01}+O(\beta^2+(\lambda\chi)^{-1})=0.\end{equation}
We next choose a middle section $\cS^M:=\{ -\beta \bx_1^L=\bx_2^R\}$ and change from left to right Jacobi coordinates on it (noting that $\bx_1^L\parallel \bx_2^R$ with opposite signs). When restricted to $\cS^M$, we estimate from equation \eqref{EqLR} that $|\bx_1^L|\sim\lambda\chi,\ |\bx_2^L|\sim\lambda\chi,\ |\bx_2^R|\sim \sqrt\lambda\chi$, and $|\bx_1^R|\sim\lambda\chi.$ When traveling from $\cS_+^L$ to $\cS^M$, we estimate the change of  $|\bp_2|$ and $E_2$ as $\beta\cdot\chi^{-1}$ by integrating the Hamiltonian equations. By conservation of energy, the change in $E^L_1+E_0+U_{01}$ is on the same order. Evaluating  the energy on $\cS^M$, we have \begin{equation}\label{EqEnergyPartition2}
E^L_1=\frac{1}{2M^L_1}|\bp^L_1|^2+O((\lambda\chi)^{-1}), \quad U_{01}+U_2=O(\beta\cdot\chi^{-1}),\quad E_1^L+E_0+O(\beta^2+\beta\cdot\chi^{-1})=0.\end{equation}
Next we apply \eqref{EqLR} to transition from the left to the right. We get on the section $\cS^M$ that
\begin{equation}\label{EqpLR}0=\left(\frac12 M_1^L+\beta\right) \bx_1^L+ \bx_2^L,\quad \bp_1^R=\bp_1^L+O(\beta),\quad \bp_2^R=\frac12 M_2^R\bp_1^L+O(\beta).\end{equation}
In right Jacobi coordinates,  on  the middle section $\cS^M$, we have $U_{01}+U_2\sim\frac{1}{\sqrt\lambda\chi}$. Therefore, on $\cS^M$,
%$0=H=E_0+E_1^R+E_2^R+O(\frac{1}{\beta\chi})$.
%So we get $$E_1^R+E_2^R=E^L_1+O(\frac{1}{\beta\chi}), \quad \frac{|\bp_1^R|^2}{2M^R_1}+\frac{|\bp_2^R|^2}{2M^R_2}=\frac{|\bp_1^L|^2}{2M^L_1}+O(\frac{1}{\beta\chi}).$$
\[E_1^R=\frac{1}{2M_1^R}|\bp_1^L|^2+O(\beta),\;E_2^R=\frac{1}{2M_2^R}\frac{1}{4}(M_2^R)^2|\bp_1^L|^2+O(\beta),\]\[ E_1^L=\frac{1}{2M_1^L}|\bp_1^L|^2+O(\frac{1}{\lambda\chi})=-E_0+O(\beta).\]
Hence, we have
$\frac{E_0}{E_2^R}=-\frac{4}{M_1^LM_2^R}+O(\beta).$
%\[E_1^R=\frac{M_1^R(k_1^R)^2}{2(L_1^R)^2}=\frac{1}{2M_1^R}|p_1^L|^2,\;E_2^R=\frac{M_2^Rk_2^R}{2(L_2^R)^2}=\frac{1}{2M_2^R}\frac{1}{4}(M_2^R)^2|p_1^L|^2,\; E_1^L=\frac{1}{2M_1^L}|p_1^L|^2=\frac{M_1^L(k_1^L)^2}{2(L_1^L)^2}.\]
%\[|p_1^L|^2=\frac{(M_1^L)^2(k_1^L)^2}{(L_1)^2}.\]
%\[\frac{1}{(L_1^R)^2}=\frac{(M_1^L)^2(k_1^L)^2}{(M_1^Rk_1^R)^2}\frac{1}{(L^L_1)^2},\quad \frac{1}{(L_2^R)^2}=\frac{(M_1^L)^2(k_1^L)^2}{4(k_2^R)^2}\frac{1}{(L^L_1)^2}.\]
Recall that by definition $M_0=\frac{1}{2}$, $k_0=1$, $k_1^L=2m_1$, $k_2^R=2m_2,$ $M_1^L=\frac{2m_1}{m_1+2}$, $M_2^R=\frac{2m_2}{m_2+2}$.
Thus
\[E_1^L=\frac{M_1^L(k_1^L)^2}{2(L_1^L)^2},\;E_2^R=\frac{M_2^R(k_2^R)^2}{2(L_2^R)^2}, \; E_0=-\frac{M_0(k_0)^2}{2(L_0)^2},\]
($\bx_2^R$ and $\bx_1^L$ are in hyperbolic motions and $\bx_0$ is in elliptic motion). Therefore, we get
\begin{equation}\label{LRL1L2}\frac{1}{L_1^R}=\frac{4(m_1+m_2+2)}{(m_1+2)(m_2+2)(m_2+1)}\frac{1}{L_1^L}+O(\beta),\quad \frac{1}{L_2^R}=\frac{m_1^2}{4m^2_2(m_1+2)}\frac{1}{L_1}+O(\beta).\end{equation}
%$$E^R_2=\frac{M_2^R (k_2^R)^2}{2(L_2^R)^2}=\frac{|\bp_2^R|^2}{2 M_2^R}+O(\frac{1}{\beta\chi}),\quad E_0=-\frac{M_0 k_0^2}{2L_0^2}=-E_1^L+O(\beta^2+\frac{1}{\beta\chi}),$$
and  in the limit $\beta\to 0$ we have that
\begin{equation}\label{EqEnergyPartition3}\frac{-E_0^R}{E_2^R}=\frac{4}{M_1^LM_2^R}=\frac{(m_2+2)(m_1+2)}{m_1m_2},\quad \frac{L^2_0}{L_2^2}=\frac{M_1^L M_0k_0^2}{4(k_2^R)^2}=\frac{m_1}{16m_2^2(m_1+2)}. \end{equation}
This is the condition determining $\gamma_I$. Along the orbit  $\gamma$ of the I4BP going from the section $\cS^M$ to $\cS_-^R$ as in the statement, the change of $E_1$ is estimated as $(\lambda\chi)^{-1}\to 0$ as $\dt\to 0$. Evaluating $H_\therefore(\gamma(T))$, we see that $|H_\therefore(\pi_\therefore^2\gamma(T))-H_\therefore(\gamma_I(0))|\to 0 $ as $\dt\to 0$.
 On the other hand,  the stable manifold $W^s_\therefore(O_2)$, the energy level set $H_\therefore^{-1}(E)$ and the section $\cS_-^L$ intersect at a point that depends on $E$ analytically for the I3BP $Q_2$-$Q_3$-$Q_4$. Thus, by assumption (2), we also have $|\pi_\therefore^2\gamma(T)-\gamma_I(0)|\to 0$.
\end{proof}

\section{The $C^1$ dynamics of the I4BP}\label{SDerLl}
In this section, we give the proof of Proposition \ref{PropDL1} (local map)  and the outline of the proof of  Proposition \ref{PropDG1} (global map) and provide some preliminary estimates. We postpone the technical calculations of  Proposition \ref{PropDG1} to Appendix \ref{SSDG}. 
%There are two subsections. %In Section \ref{SSGC0}, we estimate the equations of motion in Delaunay coordinates. From the definition of the sections and the global map, we see that the system can be considered as an $\eps$-dependent small perturbation of three Kepler problems. In Section \ref{SSDG}, we give the proof of Proposition \ref{PropDLl}.

%\subsection{The derivative formula, the equations of motion and estimates of the potential}\label{SSGC0}
\subsection{Derivative of the local map, proof of Proposition \ref{PropDL1}}
In this section, we give the proof of Proposition \ref{PropDL1}. 
\begin{proof}[Proof of Proposition \ref{PropDL1}]
By the Hartman-Grobman Theorem we know that there exists a small neighborhood $\mathscr{O}$ of the hyperbolic fixed point $O_1$ such that the dynamics is $C^1$-conjugate to the linearized dynamics at the fixed point.  Let us consider two  auxiliary sections $A_\pm=\{v_*\pm a\}$ with $a>0$ small enough that they are  inside $\mathscr O$.  Let $\zeta_-$ be the curve where the orbits starting from $\zeta_0$ insect $A_-$ and let $\zeta_+$ be the same for $A_+$.  We now parametrize the points  on  the  curve $\zeta_-$ by their distance $\delta$ to the stable manifold $W_\therefore^s(O_1)$. Without loss of generality, we can assume the curve $\zeta_-$ is straight. Then there exist $C_1, C_2>0$ and $0<\lambda_1<-v_*$\footnote{the numbers $\lambda_1$ and $-v_*$ are the positive eigenvalues of the linearization at $O_1$.} such that for the the linear dynamics at the fixed point $O_1$ and the orbits starting from $\zeta_-$,  we have 
\[r(t)=r_0e^{v_*t},\; v(t)-v_*=-ae^{v_*t}+C_1\delta (e^{-v_*t}-1),\; \psi(t)-\psi_*=e^{-\lambda_1t}(\psi_0-\psi_*)+C_2\delta e^{\lambda_1t}.\] The travel time $\tau_1$ between the two sections $\zeta_-$ and $\zeta_+$ is such that $-ae^{v_*\tau_1}+C_1\delta (e^{-v_*\tau_1}-1)=a$,  then   $\tau_1= \frac{1}{-v_*}\log\frac{a}{C_1}\delta^{-1}+O(\delta)$. Therefore, on the curve $\zeta_+$, we have $\frac{dr}{d\delta}\simeq r\frac{1}{\delta}$ and $\frac{d\psi}{d\delta}\simeq \delta^{-\alpha}$, where $0\leqslant\alpha<1$, depending on $v_*$ and $\lambda_1$. So, on the curve $\zeta_+$,  we have $\frac{dr}{d\psi}=r\frac{1}{\delta^{1-\alpha}}$. On the section $\cS_+^1$, we have $L_0=\epsilon^{1/2}r^{1/2}C(\psi)$.  Since the travel time between $\zeta_+$ and the section $\cS^1_+$ is uniformly bounded,  we then have $\frac{dL_0}{d\psi}\sim \epsilon^{1/2}r^{1/2}\frac{1}{\delta^{1-\alpha}}$.  With our choice of $\zeta_0$, we also have on $\zeta_-$, $|\frac{dL_2}{d\delta}|\ll\delta$ and $\frac{d\ell_2}{d\delta}=O(1)$.  Now the local map acts almost as an identity for the variables $(dL_2,d\ell_2)$.  Therefore, providing $\delta\ll\eta$, the assertion of the lemma follows from the fact that  the renormalization map stretches  $L_0$ and  $L_2$ by a factor $\sqrt{\lambda}=\epsilon^{-1/2}r^{-1/2}$, leaving $\psi$ and $\ell_2$ untouched and $\frac{d\psi}{d\ell_0}\simeq1$. This gives the estimate of the derivative $d\mathcal{R}d\mathbb L$ in Proposition~\ref{PropDL1}. 
\end{proof}
\subsection{Derivative of the global map, proof of Proposition \ref{PropDG1}}
We next show how to estimate the derivative of the global map. The I4BP is described by a set of six Delaunay coordinates
$(L_0,\ell_0,L_1,\ell_1,L_2,\ell_2).$ After the energy reduction, we use the variables $(L_0,\ell_0, L_2,\ell_2)$ for the local map near the triple collision of $Q_1$-$Q_3$-$Q_4$.  But for the global map defined along orbits traveling from  near the triple collision of $Q_1$-$Q_3$-$Q_4$ to the triple collision of $Q_2$-$Q_3$-$Q_4$, we use the variables $(L_0,\ell_0, L_2,\ell_2)$ for the initial condition and $(L_0,\ell_0, L_1,\ell_1)$ for the final condition. 

\subsubsection{The derivative formula}
We first recall a formula from   \cite[Section 7]{X} computing the derivative of the Poincar\'e map defined by cutting the flow $\mathcal V'=\mathcal F(\mathcal V,t)$
 between the sections $S^i$ and $S^f.$ We use $\mathcal V^i$ to denote the values of variables $\mathcal V$ restricted on the {\it initial} section $S^i$, while $\mathcal V^f$ means values of $\mathcal V$ on the {\it final} section $S^f$, and we use $t^i$ to denote the initial time and $t^f$   the final time. We want to compute the derivative
$\cD:=\frac{D \mathcal V^f}{D\mathcal V^i}$ of the Poincar\'{e} map.%along the orbit starting from $(\mathcal V^i_*, \ell_*^i)$ and ending at $(\mathcal V^f_*, \ell_*^f).$
 We have
\begin{equation}
\cD
=\left(\Id-\mathcal F(\ell_4^f)\otimes \frac{D t^f}{D\mathcal V^f}\right)^{-1}\frac{D\mathcal V(t^f)}{D\mathcal V(t^i)}\left(\Id-\mathcal F(t^i)\otimes \frac{D t^i}{D \mathcal V^i}\right).
\label{eq: formald4}
\end{equation}
Here the middle term $\frac{D\mathcal V(t^f)}{D\mathcal V(t^i)}$ is the fundamental solution to the variational equation $(\dt \mathcal V)'=\nabla_{\mathcal V}\mathcal F(\mathcal V, t)\dt \mathcal V$ and the two terms $\left(\Id-\mathcal F(\ell_4^f)\otimes \frac{D t^f}{D\mathcal V^f}\right)^{-1}$ and $\left(\Id-\mathcal F(t^i)\otimes \frac{D t^i}{D \mathcal V^i}\right)$ are called boundary contributions taking into account the issue that different orbits take different times to travel between the two sections. %We refer readers to \cite{X} for the derivation of this formula.
%\subsubsection{The Hamiltonian}

In terms of Delaunay coordinates, %we first write the three-body subsystem as follows
%$$H_\therefore=-\frac{M_0k_0^2}{2L_0^2}-\frac{M_1k_1^2}{2L_1^2} +U_{01},$$
we write the full four-body problem Hamiltonian as
\begin{equation}
\begin{aligned}
H_{\cdot:\cdot}&=-\frac{M_0k_0^2}{2L_0^2}+\frac{M_1k_1^2}{2L_1^2} +\frac{M_2k_2^2}{2L_2^2} +U,\quad U=U_{01}+U_2,
\end{aligned}
\end{equation}
where $U_2$ and $U_{01}$ are the restriction of that in \eqref{EqHamLoc} and \eqref{EqHamGlob} to the isosceles case: %$k_0=1,\; k_1=2m_1,\; k_2=(m_1+2)m_2$,
  $$U_{01}=\left(\frac{k_1}{r_1}-\frac{2m_1}{(r_1^2+r_0^2/4)^{1/2}}\right), \quad U_2=\left(\frac{k_2}{ r_2}-\frac{m_1m_2}{r_2+\frac{2r_1}{m_1+2} } -\frac{2m_2}{\big((r_2-\frac{m_1r_1}{m_1+2})^2+r_0^2/4\big)^{1/2}}\right).$$

Note that  when $ r_1\gg r_0$, we have
$U_{01}=\frac{k_1}{r_1}(1-\frac{1}{(1+r_0^2/4r_1^2)^{1/2}})=\frac{k_1r_0^2}{8r_1^3}+O(\frac{r_0^4}{r_1^5}).$
%\[\begin{split}U_2=&\frac{m_2(m_1+2)}{r_2}-\frac{m_1m_2}{r_2+\frac{2r_1}{m_1+2}}-\frac{2m_2}{r_2+\frac{m_1r_1}{m_2+2}}.\\
%&=\frac{m_1m_2}(\frac{r_2+\frac{2r_1}{m_1+2}-r_2}{r_2(r_2+\frac{2r_1}{m_1+2}}+2m_2(\frac{r_2+\frac{m_1r_1}{m_2+2}-r_2}{r_2(r_2+\frac{m_1r_1}{m_2+2})}\\
%&=m_1m_2\frac{\frac{2r_1}{m_2+2}}{r_2(r_2+\frac{2r_1}{m_1+2}}+2m_2\frac{\frac{m_1r_1}{m_2+2}}{r_2(r_2+\frac{m_1r_1}{m_2+2}}
%\end{split}\]
 \subsubsection{The strategy of the proof }\label{SSStrategyLoc}
We are now ready to give a proof of Proposition \ref{PropDG1}. The strategy is as follows. Suppose we have an orbit going from the left to the right. Suppose it starts from the left section $\cS^L_+=\{v=\epsilon^{-1/2}\}$ and ends at the right section $\cS^R_-=\{r=\epsilon^{-1}\}$. On the left section $\cS^L_+=\{v=\epsilon^{-1/2}\}$ and before the application of the global map, we have done a renormalization. Let us denote by $\lambda$ the renormalization factor and introduce $\beta=\lambda^{-1/2}$.
We define the middle section $\cS^M:=\{-\beta \bx_1^L=\bx_2^R\}$. We also introduce a large parameter $\chi$ giving a lower bound for dist$(Q_1,Q_2)$ after renormalization on $\cS_+^L$, so $\chi$ already contains the $\lambda$ factor (note that the convention for $\chi$ here differs from that in Section \ref{SSPartition}).

When the binary is moving between the section $\cS^L_+$ and $\cS^M$, we eliminate $L_1$ using the total energy conservation, and treat $\ell_1$ as the new time, while using formula \eqref{eq: formald4} to obtain the derivative matrix denoted by $\cD_1=\frac{\partial (L_0,\ell_0,L_2,\ell_2)|_{\cS^M}}{\partial (L_0,\ell_0,L_2,\ell_2)|_{\cS^L_+}}$. Next, restricting to the middle section $\cS^M$, we change the coordinate system from the left to the right. We first convert Delaunay coordinates to Cartesian coordinates. Then we apply \eqref{EqLR} in Section \ref{STransition} to convert the left Jacobi coordinates to the right ones. Finally, we convert the right Cartesian coordinates to the Delaunay coordinates. We obtain $\cD_2=\frac{\partial (L_0,\ell_0,L_1,\ell_1)|_{\cS^M}}{\partial (L_0,\ell_0,L_2,\ell_2)|_{\cS^M}}$.  When moving on the right of the section  $\cS^M$, the calculation is similar to the left case except that we need to permute the subscripts $1$ and $2$. In this case, we eliminate $L_2$ using the total energy conservation and treat $\ell_2$ as the new time,  using formula \eqref{eq: formald4} to obtain the derivative matrix denoted by $\cD_3=\frac{\partial (L_0,\ell_0,L_1,\ell_1)|_{\cS^R_-}}{\partial (L_0,\ell_0,L_1,\ell_1)|_{\cS^M}}$. Then $d\mathbb G$ is obtained by composing $\cD_3\cD_2\cD_1$ and Proposition \ref{PropDG1} follows. 

In the following estimates, we shall use the notations $\sim$ and $\simeq$ in the following way.
\begin{Not}Let $A(\chi)$ and $B(\chi)$ be two functions depending on $\chi$ and we consider the limit as $\chi\to\infty$.
\begin{enumerate}
\item We denote $A(\chi)\sim B(\chi)$, if there exists a constant $C>1$ such that $C^{-1}\leq \frac{A(\chi)}{B(\chi)}\leq C$;
\item We denote $A(\chi)\simeq B(\chi)$, if $\frac{A(\chi)}{B(\chi)}\to 1.$
\end{enumerate}
\end{Not}

Following the strategy outlined above, we give the proof of Proposition \ref{PropDG1} deferred to Appendix \ref{SSDG}. Piercing through the technicalities, the essential point of the proof is to reveal the estimate $\frac{\partial \ell_0}{\partial L_0}\sim \chi$ in Proposition \ref{PropDG1}, which is essential in the proof of Lemma \ref{LmNondeg}. Indeed, this is clear from the Kepler problem, which has the Hamiltonian $H=-\frac{mk^2}{2L^2}$ and equations of motion $\begin{cases}\dot L=0\\
\dot \ell=\frac{mk^2}{L^3}\end{cases}$. Then we can derive the variational equation by differentiating the Hamiltonian equations and the fundamental solutions gives the above estimate. In fact, in the proof of Proposition \ref{PropDL1}, most of the work is devoted to showing that the perturbations do not spoil this estimate.

 In the remaining part of this section, we give the estimate of the Hamiltonian equation as well as the changes in $L_i,\ i=0,1,2. $
\subsubsection{Preliminary estimates of the potentials}
The setting is as follows. Suppose the triple $Q_1$-$Q_3$-$Q_4$ is just ejected from the near triple collision so $|E_0|$ and $|E_1|$ are very large and $E_2$ is of order $1$. We perform a renormalization with the outcome in Proposition \ref{PropRenorm}. In particular, $E_0$ and $E_1$ are of order $1$ but $E_2$ is of order $\beta^2. $ From the Hamiltonian equations, we get that the changes in $E_0,E_1,E_2/\beta^2$ are estimated as $O(\eps^2)$, following the orbit up to the section $\cS^M$ (see Lemma \ref{L01-bound} below). %For convenience, we introduce a number $\beta=\lambda^{-1/2}$ which takes a small value around $\sqrt E_2$ of the left case, which is the same as the definition of $\beta$ in Section \ref{SSPartition}.
However, when we switch from left to right using the transformation \eqref{EqLR}, from the argument in Section \ref{SSPartition}, we see that the energies $E_0, E_1$ and $E_2$ are of the same order.
%Since in the left case the smallness of $E_2$ gives us some good estimates while in the right case $E_1$ (which is of order $1$) does not, we pick the middle section to be\footnote{the  minus sign here is due the fact that $x_1^L$ and $x_2^R$ are of opposite direction.} $$\cS^M=\{-\beta \bx_1^L=\bx_2^R\}$$ so that the orbit stays for shorter time $O(\beta\chi)$ in the right case, which will give us good estimate in the right case.
%When moving from the left section $\cS_+^L$ to the middle section $\cS^M$, we have
%\[r_0\sim 1,\quad r_2\sim \chi,\quad \text{and} \quad r_1\text{ increases from the order }\frac{1}{\epsilon} \text{ to the order }\chi. \]
Applying Proposition \ref{PropRenorm}, we have on the section $\cS^M$ that
$$\bx_1^L\sim \bx_2^L\sim \bx_1^R\sim\chi,\quad \bx_2^R\sim \beta\chi, \quad L_1^L\sim L_1^R\sim L_2^R\sim 1,\quad L_2^L\sim 1/\beta.$$
We always assume $0<\frac{1}{\chi}\ll\beta^2\ll\epsilon^4\ll1.$

In the following calculations, we mainly focus on the left case, that is, the transition between the sections $\cS_+^L$ and $\cS^M$. The right case can be dealt with similarly.%\footnote{ In fact, a large part of the estimates for the right cases can be  obtained automatically by permuting the subscripts 1 and 2 and replacing  $\beta$ by $1$.}.

 \begin{Lm}\label{LmPotential}
Consider the piece of orbit between the sections $\cS_+^L$ and $\cS^M$. Assume there exist some constants $C>1$ and $\beta\ll 1$ such that
\begin{equation}
1/C<L_1<C,\quad 1/C<L_0<C, \quad 1/C<L_2\beta<C ,
\end{equation}
Then we have the following estimate when $\ell_1\gg1$ and $\chi\gg1$
\begin{equation}
\begin{aligned}
\frac{\partial U_{01}}{\partial r_1}&=O\left(\frac{1}{\ell_1^4}\right),\quad \frac{\partial U_{01}}{\partial r_0}=r_0O(\frac{1}{\ell_1^3}), \\
\frac{\partial^2 U_{01}}{\partial r_1^2}&=O(\frac{1}{\ell_1^5}),\quad \frac{\partial^2 U_{01}}{\partial r_0^2}=O(\frac{1}{\ell_1^3}),\quad \frac{\partial^2 U_{01}}{\partial r_0\partial r_1}=O(\frac{1}{\ell_1^4}),\\
\frac{\partial U_2}{\partial r_1}&=O(\frac{1}{\chi^2}),\quad \frac{\partial U_2}{\partial r_0}=r_0O(\frac{1}{\chi^3}),\quad \frac{\partial U_2}{\partial r_2}=O(\frac{1}{\chi^2}),\\
\frac{\partial^2 U_2}{\partial r_1^2}&=O(\frac{1}{\chi^3}),\quad \frac{\partial^2 U_2}{\partial r_0^2}=O(\frac{1}{\chi^3}),\quad \frac{\partial^2 U_2}{\partial r_2^2}=O(\frac{1}{\chi^3}),\\
\frac{\partial^2 U_2}{\partial r_1\partial r_0}&=r_0O(\frac{1}{\chi^4}),\quad \frac{\partial^2 U_2}{\partial r_0\partial r_2}=r_0O(\frac{1}{\chi^4}),\quad \frac{\partial^2 U_2}{\partial r_2\partial r_1}=O(\frac{1}{\chi^3}).\\
\end{aligned}
\end{equation}
\end{Lm}
The proof consists of straightforward calculations, and we skip the details. %After the renormalization, we get that $E_0$ and $E_1$ are of order $1$ and $E_2:=\beta^2\ll1$. When the orbit moves to the right of the section $x_1^L=\chi/2$, we after the coordinate change \eqref{EqLR}, we note that all $E_1$, $E_0 $ and $E_2$ are of order 1 (note $|p_2^L|\sim \beta\ll1$).
The Hamiltonian equation has the following estimates under the assumption of the previous lemma
\begin{equation}\label{EqHamLimit}
\begin{cases}
\dot L_0&=\left(\frac{\partial U_{01}}{\partial r_0} +\frac{\partial U_2}{\partial r_0}\right)\frac{\partial r_0}{\partial \ell_0}=O\left(\frac{1}{\ell_1^3}+\frac{1}{\chi^3}\right)r_0\frac{\partial r_0}{\partial \ell_0}\\%=-\frac{1/2}{(r_1^2+r_0^2/4)^{3/2}}\frac{L_0^4}{(M_0k_0)^2}\sin u_0+\frac{1}{\chi^3},\\
\dot \ell_0&=\frac{M_0k_0^2}{L_0^3}+\left(\frac{\partial U_{01}}{\partial r_0} +\frac{\partial U_2}{\partial r_0}\right)\frac{\partial r_0}{\partial L_0}=\frac{M_0k_0^2}{L_0^3}+O\left(\frac{1}{\ell^3_1}+\frac{1}{\chi^3}\right) r_0\frac{\partial r_0}{\partial L_0}\\%=\frac{M_0k_0^2}{L_0^3}+\frac{r^2_0}{L_0(r_1^2+r_0^2/4)^{3/2}}+\frac{1}{\chi^3},\\
\dot L_1&=\left(\frac{\partial U_{01}}{\partial r_1} +\frac{\partial U_2}{\partial r_1}\right)\frac{\partial r_1}{\partial \ell_1}=O\left(\frac{1}{\ell_1^4}+\frac{1}{\chi^2}\right)\frac{\partial r_1}{\partial \ell_1}\\%=\left(\frac{2}{r_1^3}-\frac{4}{(r_1^2+r_0^2/4)^{3/2}}\right)\left(\frac{L_1^2}{M_1k_1}\right)^2(-\sinh u_1)+\frac{1}{\chi^2},\\
\dot \ell_1&=-\frac{M_1k_1^2}{L_1^3}+\left(\frac{\partial U_{01}}{\partial r_1} +\frac{\partial U_2}{\partial r_1}\right)\frac{\partial r_1}{\partial L_1} =-\frac{M_1k_1^2}{L_1^3}+O\left(\frac{1}{\ell_1^4}+\frac{1}{\chi^2}\right)\frac{\partial r_1}{\partial L_1}\\
%=-\frac{M_1k_1^2}{L_1^3}+\left(\frac{2}{r_1^2}-\frac{4r_1}{(r_1^2+r_0^2/4)^{3/2}}\right)\frac{2 r_1}{L_1}+\frac{1}{\chi},\\
\dot L_2&=-\frac{\partial U_2}{\partial r_2}\frac{\partial r_2}{\partial \ell_2}=O(\frac{1}{\beta^2\chi^2}),\\
\dot \ell_2&=\frac{M_2k_2^2}{L_2^3}+\frac{\partial U_2}{\partial r_2}\frac{\partial r_2}{\partial L_2}=\frac{M_2k_2^2}{L_2^3}+O(\frac{\beta}{\chi}).%=\frac{M_2k_2^2}{L_2^3}+\frac{1}{\chi}.
\end{cases}\end{equation}
Here we use the formula in Appendix \ref{App-delaunay} to obtain
\[\frac{\partial r_0}{\partial L_0}=O(1),\ r_0\frac{\partial r_0}{\partial \ell_0}=O(1),\ \frac{\partial r_1}{\partial L_1}=O(\ell_1),\ \frac{\partial r_1}{\partial \ell_1}=O(1),\ \frac{\partial r_2}{\partial L_2}=O(\beta\chi),\ \frac{\partial r_2}{\partial \ell_2}=O(\beta^{-2}).\]
\begin{Lm} \label{L01-bound}For  the I4BP orbits between the section $\cS_+^L$ and $\cS^M$, there exists $C>0$ such that if,  on the section $\cS_+^L$, we have
\[\frac{1}{C}\leqslant L_1(0)\leqslant C,\quad \frac{1}{C}\leqslant L_0(0)\leqslant C,\quad \frac{1}{C}\leqslant L_2(0)\beta\leqslant C,\]
then, for all time between the two sections, we have
\[\max\{ |L_1-L_1(0)|,\; |L_0-L_0(0)|, \;| L_2\beta-L_2(0)\beta|\}=O(\epsilon^2).\]
\end{Lm}
This can be proved by a simple bootstrap argument with the estimates in \eqref{EqHamLimit}. The argument goes as follows. We start by assuming  the change in $L_i,\ i=0,1,2$ is bounded by a generous constant $C\gg \eps^2$.  This implies that $r_0$ is bounded, and $r_1$ and $\ell_1$ grow at least linearly in time. Then integrating the estimates in \eqref{EqHamLimit} over time of order $\chi$, we get that the change in each $L_i$ is indeed $O(\eps^2)$, which is the value of $(\ell_1^i)^{-2}=\int_{\ell_1^i}^\chi\frac{1}{\ell_1^3}$ evaluated on the section $\cS_+^L$. We refer readers to Lemma 6.6 of \cite{X} for details.

\section{The derivative of the local map in the F4BP}\label{local-gs}
In this section, we estimate the derivative of the local map in the F4BP when the variables $(G_0,g_0, G_1,g_1,G_2,g_2)$ are nonvanishing. Without loss of generality, we consider the local interaction of $Q_1$-$Q_3$-$Q_4$. Reducing the angular momentum conservation, we use the variables $X=(L_0,\ell_0,L_2,\ell_2)$ and $Y=(G_0,\mathsf g_{01}=g_0-g_1-\frac{\pi}{2},G_2,\mathsf g_{21}=g_2-g_1)$. 

By the block diagonal form of the derivative matrix in Lemma \ref{LmDiagonal} in the limit $Y\to 0$, it is enough to work out the two blocks independently in order to prove Proposition \ref{PropDL2}. Indeed, in the variational equation, each off-diagonal term corresponding to the mixed partial derivatives with respect to $X$ and $Y$ is estimated as $O(\nu)$ if we assume $|Y|\leq \nu$. We then apply the DuHamel principle to show that the nonisoscelesness causes a perturbation of $O(e^{MT}\nu)$ to the estimate of $d\mathbb L$ in the isosceles limit, where $T$ is the time defining the local map in the blowup coordinates and $M$ bounds the unperturbed part of the derivative (see the proof of Proposition \ref{PropLocalMM}). Thus this perturbation can be suppressed by bounding $T$ from above by the sojourn time and choosing $\nu$ small accordingly. Therefore for the derivative of the local map, we shall neglect the off-diagonal terms caused by the nonisoscelesness $Y\neq0$ and it is enough for us to work with the $Y$-part of the derivative of the local map in the isosceles limit $|Y|\to 0$, since the $X$-part is given in Proposition \ref{PropDL1}. 

We divide this section into two parts. In Section \ref{SSG0}, we study the derivative of the local map in the F3BP $Q_1$-$Q_3$-$Q_4$. In Section \ref{SSG02}, we study the F4BP. The former embeds into the latter containing the essential information of the $Y$-component of the local map so we study it separately. We also recall the definitions of $E^{u/s}_{3}$ and  $E^{u/s}_{4}$ from Section \ref{SSLoc}. 

%Since the Hamiltonian equation $\dot Y=JD_Y H=0$ along an I4BP orbit,  the the boundary contributions from  the sections $\cS_{+}^L$, $\cS^M$ and $\cS_-^R$ all read as $\mathrm{Id}\pm JD_YH\otimes \nabla_Y \ell_{1,2}=\mathrm{Id}$. Therefore  the boundary terms are trivial.

\subsection{The $C^1$ dynamics of $w_0,g_0$ for F3BP}\label{SSG0}
In this subsection, we restrict ourself to the F3BP $Q_1$-$Q_3$-$Q_4$ with zero total angular momentum, that is, $G_0+G_1=0$.  We study the variational equation for the variables $w_0:=r^{-1/2}G_0$ and $\mathsf{g}_{01}:=g_0-g_1-\frac{\pi}{2}$.
The  equations of motion  for $(w_0,\mathsf{g}_{01})$ take the form
\begin{equation}\label{w0g0eq}
\begin{cases}w_0'=-\frac{1}{2}vw_0-r\frac{\partial U_{01}}{\partial g_0},\\
\mathsf{g}_{01}'=r^{3/2}(\frac{\partial U_{01}}{\partial G_0}-\frac{\partial U_{01}}{\partial G_1}),
\end{cases}\end{equation}
where the potential $U_{01}$ is as in \eqref{EqHamGlob}.

\begin{Lm}\label{LmVar0}
\begin{enumerate}\item The variational equation has the following form with the coefficients matrix evaluated along an I3BP orbit
\begin{equation}\label{variation-w0g0}
\left[\begin{array}{c}
\dt w_0\\
\dt \mathsf{g}_{01}\end{array}\right]'=\left[\begin{array}{cc}
-\frac{1}{2}v-r^{3/2}(\frac{\partial^2 U_{01}}{\partial G_0\partial g_0}-\frac{\partial^2 U_{01}}{\partial g_0\partial G_1}) &-r\frac{\partial^2 U_{01}}{\partial g_0^2}\\
r^2(\frac{\partial^2U_{01}}{\partial G_0^2}-2\frac{\partial^2U_{01}}{\partial G_0\partial G_1}+\frac{\partial^2 U_{01}}{\partial G_1^2})&r^{3/2}(\frac{\partial^2 U_{01}}{\partial G_0\partial g_0}-\frac{\partial^2 U_{01}}{\partial G_1\partial g_0})\end{array}\right]\left[\begin{array}{c}
\dt w_0\\
\dt \mathsf{g}_{01}\end{array}\right].\end{equation}
\item When evaluated at the Lagrange fixed points with $v_*<0$, the coefficient matrix has two eigenvalues $\mu_u>-\frac{v_*}{2}$ and $\mu_s<0$.
\item For any small $\eta>0$, there is $\bar d$ such that for any orbit $\gamma:\ [0,T]\to \cM_\therefore$ of the I3BP with initial condition $\bar d$-close to $\gamma_I(0)$ and $\gamma(T)\in \mathcal S_+^1$, the fundamental solution of \eqref{variation-w0g0} sends any vector $v\in T_{(\gamma(0))}\mathcal S_-^1$ with $v\cap \mathcal C_\eta(E_3^s(\gamma_I(0)))=\{0\}$ to the cone $\mathcal C_\eta(E_3^u(\gamma_O(0)))$, and $v$ expands at a rate of $e^{\mu_u T}$. 
 %The fundamental solution of the variational equation along a near triple collision I3BP orbit $\gamma: \ [0,T]\to \cM$ with $\gamma(0)\in \cS_-^1$ and $\gamma(T)\in \cS_+^1$, which first follows $\gamma_I^1$ and then follows  $\gamma_O^1$  has the form 
%$$e^{T\mu_u}\mathbf{v}_u\otimes \mathbf{v}_s^{\bot}+O(1), $$
%where we have $\mathbf{v}_s=E^s(\gamma_I^1(0))\in T\cS^1_-$, and $\mathbf{v}_u=E^u(\gamma_O^1(0))\in T\cS^1_+$. 

\end{enumerate}
\end{Lm}
\begin{proof}
From \eqref{w0g0eq}, we obtain
the variational equation
$$\left[\begin{array}{c}
\dt w_0\\
\dt \mathsf{g}_{01}\end{array}\right]'=\left[\begin{array}{cc}
-\frac12 v -r\frac{\partial^2U_{01}}{\partial w_0\partial g_0}&-r\frac{\partial^2U_{01}}{\partial g_0\partial\mathsf{g}_{01}}\\
r^{3/2}(\frac{\partial^2 U_{01}}{\partial w_0\partial G_0}-\frac{\partial^2 U_{01}}{\partial w_0\partial G_1})&r^{3/2}(\frac{\partial^2 U_{01}}{\partial G_0\partial \mathsf{g}_{01}}-\frac{\partial^2U_{01}}{\partial G_1\partial \mathsf{g}_{01}})\end{array}\right]\left[\begin{array}{c}
\dt w_0\\
\dt \mathsf{g}_{01}\end{array}\right] $$
Using the conservation of angular momentum $G_0+G_1=0$ and the fact that only the term $\mathsf g_{01}=g_0-g_1$ enters the potential $U_{01}$,   we have 
$\frac{\partial }{\partial w_0}=r^{1/2}(\frac{\partial }{\partial G_0}-\frac{\partial}{\partial G_1}), \ \frac{\partial }{\partial\mathsf{g}_{01}}=\frac{\partial }{\partial g_0}.$ Therefore,
\[\frac{\partial^2 U_{01}}{\partial w_0\partial g_0}=r^{1/2}(\frac{\partial^2  }{\partial G_0\partial g_0}-\frac{\partial^2  }{\partial g_0\partial G_1})U_{01},\quad  \frac{\partial^2 U_{01}}{\partial w_0\partial G_i}=(-1)^ir^{1/2}(\frac{\partial^2  }{\partial G_i^2}-\frac{\partial^2  }{\partial G_0\partial G_1})U_{01},\; i=0,1.\]
This gives us \eqref{variation-w0g0}.  From Lemma \ref{w0g0-deri}, we have an explicit formula  for the quantities in the coefficient matrix, written in the blowup coordinates so that assertion (2) can be verified explicitly. 
% Recall that at the fixed point $(0,v_*,\psi_*,0)$, we have $v_*^2=-2\bar V(\psi_*)$, where the potential  $\bar V(\psi)$ is defined in~\eqref{iso-potential}.
% Then direct calculation shows that the determinant of the  coefficients matrix is negative when evaluated at the Lagrange fixed points with $v_*<0$. Hence, we finishes the proof of the assertion (2).
However, the calculation is a bit involved. In the following, we use  an easier way to prove assertion (2).  Instead of $(w_0,\mathsf{g}_{01})$, we use $(w_0, \theta)$, where $\theta=\theta_0-\theta_1)$. We get the equations of motion from the Hamiltonian equations in terms of polar coordinates.
 $$\begin{cases}
w_0'&=(r^{-1/2}\Theta_0)'=-\frac{1}{2}r^{-3/2}rv\Theta_0+r^{-1/2}\Theta_0'
=-\frac{1}{2} vw_0-r\partial_{\theta_0} V,\\
\theta'&=r^{3/2}\frac{\Theta_0}{M_0 r_0^2}-r^{3/2}\frac{\Theta_1}{M_1r_1^2}=\frac{w_0}{\cos^2\psi}-\frac{w_1}{\sin^2\psi},
\end{cases}
$$where $V$ is as in \eqref{EqHamLoc}.
Taking into account conservation of angular momentum, the variational equation reads
$$\left[\begin{array}{c}
\dt w_0\\
\dt \theta\end{array}\right]'=\left[\begin{array}{cc}
-\frac12 v -r\frac{\partial^2 V}{\partial w_0\partial \theta_0}&-r\frac{\partial^2 V}{\partial \theta_0^2} \\
\frac{1}{\cos^2\psi}+\frac{1}{\sin^2\psi}&0\end{array}\right]\left[\begin{array}{c}
\dt w_0\\
\dt \theta\end{array}\right] $$
Note that $rV$ depends only on $\psi$ and $\theta$, so the term $r\frac{\partial^2 V}{\partial w_0\partial \theta_0}=0$. We have
\begin{equation}
\begin{aligned}
rV&=r\left(-\frac{m_1}{|x_1+\frac{x_0}{2}|}-\frac{m_1}{|x_1-\frac{x_0}{2}|}-\frac{1}{|x_0|}\right)\\
&=-\sum_{i=0,1}\frac{m_1}{\left(\frac{\sin^2\psi}{M_1}+\frac{\cos^2\psi }{4M_0}+(-1)^i \frac{\sin\psi\cos\psi}{(M_1M_0)^{1/2}}\cos(\theta_0-\theta_1)\right)^{1/2}}-\frac{M_0^{1/2}}{\cos\psi}
\end{aligned}
\end{equation}
So we get
$$r\partial^2_{\theta_0} V|_{\theta_0-\theta_1=\pi/2}=-\frac{\frac32m_1M_1^{-1}M_0^{-1}\sin^2\psi\cos^2\psi}{\left(M_1^{-1}\sin^2\psi+\frac{\cos^2\psi }{4M_0}\right)^{5/2}}<0.$$
Restricted to a Lagrangian fixed point in the lower half of $\mathcal M_0$, with $v<0$, we get that the coefficient matrix of the variational equation has trace positive and determinant negative therefore it is also a saddle in the $(w_0,\theta)$ plane. When transformed to the $(w_0,\mathsf{g}_{01})$ coordinates, the fact of being a saddle does not change. At the fixed point, we have $\mu_u+\mu_s=-\frac12v_*$ is the trace of the coefficient matrix and $\mu_u\mu_s<0$. This gives that $\mu_u>-\frac12v_*$ and $\mu_s<0.$

Since the orbit spends most of the time very close to the fixed point, we know that the fundamental solution of the variational equation stretches the unstable direction at the fixed point by a factor $e^{T\mu_u}$ and contracts the stable direction by a factor $e^{T\mu_s}$. Assertion (3) follows from the fact that $E_3^s(\gamma^1_I(0))\in T\cS_-^1$ is the push-backward of the stable direction along $\gamma^1_I$, $E^u_3(\gamma^1_O(0))\in T\cS_+^1$ is the push-forward of the unstable direction along $\gamma^1_O$, and the orbit stays very close to $\gamma^1_I$, then very close to $\gamma^1_O$ for the whole passage. 
\end{proof}
\subsection{The $C^1$ dynamics of $G_0,g_0-g_1,G_2,g_1-g_2$ for the F4BP}\label{SSG02}
In this subsection we study the F4BP with total angular momentum zero, that is, $G=G_0+G_1+G_2=0$. We introduce a  transform $(G_0,g_0,G_1,g_1,G_2,g_2)\mapsto(G_0,\mathsf{g}_{01}, G_2,\mathsf{g}_{12},G,g_1)$ where $\mathsf{g}_{01}=g_0-g_1-\frac{\pi}{2}$, $\mathsf{g}_{12}=g_1-g_2-\pi$. Due to the total  angular momentum conservation and the fact that the Hamiltonian depends on $g_i$ only through  the terms $g_i-g_j$, $i,j\in\{0,1,2\}$, $i\neq j$, we just need to study the four variables $ (w_0,\mathsf{g}_{01}, w_2,\mathsf{g}_{12})$, where $w_0=r^{-1/2}G_0$ and $w_2=r^{-1/2}G_2$. 
The equations of motion are
$$\begin{cases}
w_0'&=-\frac{1}{2}vw_0-r\frac{\partial (U_{01}+U_2)}{\partial g_0}\\
\mathsf g'_{01}&=r^{3/2}\Big(\frac{\partial (U_{01}+U_2)}{\partial G_0}-\frac{\partial (U_{01}+U_2)}{\partial G_1}\Big)\\
w'_2&=-\frac{1}{2}vw_2-r\frac{\partial U_2}{\partial g_2}\\
\mathsf g'_{12}&=r^{3/2}(\frac{\partial (U_{01}+U_2)}{\partial G_1}-\frac{\partial U_2}{\partial G_2})
\end{cases}.$$
\begin{Lm}\label{LmFundLoc} Suppose the distance $|Q_1-Q_2|\geq \chi\gg 1$, then
\begin{enumerate}
\item the coefficient matrix of the  variational equation for the variables $(w_0,\mathsf{g}_{01},w_2,\mathsf{g}_{12})$ along an I4BP orbit has the  form
$$\left[\begin{array}{cccc}
-\frac{v}{2}-r^{3/2}(\frac{\partial^2U_{01}}{\partial g_0\partial G_0}-\frac{\partial^2 U_{01}}{\partial g_0\partial G_1})&-r\frac{\partial^2U_{01}}{\partial g_0^2}&r^{3/2}\frac{\partial^2U_{01}}{\partial g_0\partial G_1}&0\\
r^2(\frac{\partial^2U_{01}}{\partial G_0^2}-2\frac{\partial^2 U_{01}}{\partial G_0\partial G_1}+\frac{\partial^2 U_{01}}{\partial G_1^2})&r^{3/2}(\frac{\partial^2 U_{01}}{\partial G_0\partial g_0}-\frac{\partial^2 U_{01}}{\partial G_1\partial g_0})&r^{2}(-\frac{\partial^2U_{01}}{\partial G_0\partial G_1}+\frac{\partial^2 U_{01}}{\partial G_1^2})&0\\
0&0&-\frac{v}{2}&0\\
r^2(\frac{\partial^2 U_{01}}{\partial G_0\partial G_1}-\frac{\partial^2 U_{01}}{\partial G_1^2})&r^{3/2}\frac{\partial^2U_{01}}{\partial G_1\partial g_0}&-r^{2}\frac{\partial^2U_{01}}{\partial G_1^2}&0
\end{array}\right]+O(\frac{r}{\chi});
$$
\item  For any small $\eta>0$, there exists $\bar d$ such that for any orbit $\gamma:\ [0,T]\to \cM $ of the I4BP with initial condition $\bar d$-close to $\gamma_I(0)$ and $\gamma(T)\in \mathcal S_+^1$, the fundamental solution sends any 2-plane $P\subset T_{(\gamma(0))}\mathcal S_-^1$ with $P\cap \mathcal C_\eta(E_{4}^s(\gamma_I(0)))=\{0\}$ to the cone $\mathcal C_\eta(E_{4}^u(\gamma_O(0)))$, and each vector in $P$ gets expanded by a rate of at least $e^{-\frac{v_*}{2}}$. 
\item The projection to the $(w_0,\mathsf g_{01})$-component  of the pushforward of $\mathbf e_1$ to $T_{\gamma_I(0)}\cS_-^1$ $($respectively  the pushforward of $\mathbf e_4$ to $T_{\gamma_O(0)}\cS_+^1)$ is $E^s_3(\gamma_I(0))$ $($respectively $E^u_3(\gamma_O(0)))$.
\end{enumerate}
\end{Lm}
\begin{proof}
Note that since $|Q_2-Q_1|\geqslant\chi\gg1$, all the second derivatives of $U_2$ are of order $O(\frac{1}{\chi})$. Then the asserted form of the coefficient matrix of the variational equation for $(w_0,\mathsf{g}_{01},w_2,\mathsf{g}_{12})$ can be obtained in the same  way as in   Lemma \ref{LmVar0} by taking  into account the total angular momentum conservation. 
Let us denote the resulting matrix as
$$\mathbb{A}(\tau)=\left[\begin{array}{cc}A_{11}(\tau)&A_{12}(\tau)\\
A_{21}(\tau)&A_{22}(\tau)\end{array}\right]+O(\frac{r}{\chi}),$$ where $A_{ij}$, $ij=1,2$ are 2 by 2 matrices.  Then from the estimates in Lemma \ref{w0g0-deri}, we have 
\[A_{11}=\left[\begin{array}{cc}O(1) & O(1)  \\O(1)  & O(1) \end{array}\right], A_{12}=\left[\begin{array}{cc}a(=O(1)) & 0 \\b(=O(1)) & 0\end{array}\right], A_{21}=\left[\begin{array}{cc}0 & 0 \\-b & a \end{array}\right], A_{22}=\left[\begin{array}{cc}O(1) & 0 \\O(1) & 0\end{array}\right].\]
Note that $A_{11}$ is the same as the coefficient matrix in \eqref{variation-w0g0}. When evaluated at the fixed point $O_1$, the matrix $\mathbb{A}(\tau)$ has eigenvalues $\mu_u>-\frac{v_*}{2}>0>\mu_s$, where $\mu_u$ and $\mu_s$ are the eigenvalues of the matrix $A_{11}$, corresponding to the eigenvectors of the form respectively \begin{equation}\label{EqEigenvector}\mathbf{e}_1=(\mathbf{u}_u,0,\kappa_1), \quad \mathbf{e}_2=(*,*,1,*),\quad \mathbf{e}_3=(0,0,0,1),\quad \mathbf{e}_4=(\mathbf{u}_s,0, \kappa_2),\end{equation}% Therefore $\mathbb{A}$ has two stable directions and two unstable directions. The stable directions  are
%$$\mathbf{e}_s^1=(\mathbf{u}_s,0, \kappa_1),\quad \mathbf{e}_s^2=(0,0,0,1), $$
where $(0,\kappa_1)=A_{21}\cdot \mathbf{u}_u/\mu_u$ and $\mathbf{u}_u\in\mathbb{R}^2$ is the eigenvector of the matrix $A_{11}$ corresponding to $\mu_u$, %and the unstable directions are
%\[\mathbf{e}_u^1=(\mathbf{u}_u,0,\kappa_2), \quad \mathbf{e}_u^2=(*,*,*,*), \]
  $(0,\kappa_2)=A_{21}\cdot \mathbf{u}_s/\mu_s$ and $\mathbf{u}_s\in\mathbb{R}^2$ is the eigenvector of the matrix $A_{11}$ corresponding to $\mu_s$. %Note that $A_{11}$ is exactly the coefficient matrix of \eqref{variation-w0g0}. 

Item (2) follows from an argument similar to item (2) of Lemma \ref{LmVar0}. % Since the orbit spends most of the time inside a very small neighborhood of the fixed point and the time spending outside this neighborhood has a uniform upper bound, we approximate the orbit $\gamma$ by $\gamma_I$ and $\gamma_O$ outside the neighborhood. Denote by $P$ the two dimensional plane spanned by $\mathbf{e}_1, \ \mathbf{e}_3$ based at the Lagrange fixed point. Pushforward it along the orbit $\gamma_I$ backward to the section $\cS_-^1$, we obtain a plane denoted by  $P'$. For any vector $v\in T\cS_-^1$ that is not in the cone $\mathcal C_\eta(P')$ where $\eta$ is chosen to be a small number but is independent of the dist$(\gamma(0),\gamma_I(0))$, when pushforward by the solution of the variational equation along $\gamma$ and staying for a long time close to the fixed point, the image of the vector will be almost parallel to $\mathbf{e}_u^1$, hence when further pushforward along $\gamma_O$ to the section $\cS_+^1$, we get almost the vector $\bar \bu$.  This gives the structure in part (2) of the statement. 

%we know that the fundamental solution of the corresponding 
%variational equation expands along the direction $\mathbf{e}_u^1$ by a  the  factor $r_*^{-C_1}$, along the direction $\mathbf{e}_u^2$ by a factor $r_*^{-C_2}$ and contract along the direction $\mathbf{e}_s^1$ by a factor $r_*^{C_3}$, leaving the direction $\mathbf{e}_s^2$ almost untouched. Here the constants $C_i>0$, $i=1,2,3$ depending only on the eigenvalues.  Hence we have the assertion (2) of the lemma by denoting $\mathbf{v}_s^i$, $i=1,2$, the push-back of the the vectors $\mathbf{e}_s^i$, $i=1,2$ along $\gamma_I^1$ to the section $\cS^1_-$ and $\mathbf{v}_u^i$, $i=1,2$, the push-forward of the the vectors $\mathbf{e}_u^i$, $i=1,2$ along $\gamma_O^1$ to the section $\cS^1_+$. 

We finally prove item (3). We use  Picard iteration to integrate the  variational equation for the variables $(w_0,\mathsf{g}_{01},w_2,\mathsf{g}_{12})$.  Let us denote the $n$-th iterative integral $$\mathbb{A}_n(\tau)=\int_0^\tau\mathbb{A}(s)\cdot\mathbb{A}_{n-1}(s)ds:=\left[\begin{array}{cc}\mathbb{A}_n^{11}(\tau) & \mathbb{A}_n^{12}(\tau) \\\mathbb{A}_{n}^{21}(\tau) & \mathbb{A}_{n}^{22}(\tau)\end{array}\right],\quad \mathbb{A}_1(\tau)=\int_0^\tau \mathbb{A}(s)ds,\quad n\in\mathbb{N}.$$
 From    the fact that $A_{12}A_{21}=0$,  $A_{22}A_{21}=0$, we argue by induction, knowing that the zero entries are preserved for each  $\mathbb{A}_n(\tau)$, $n\in\mathbb{N}$, and 
\[\begin{split}\mathbb{A}_n(\tau)&=\int_0^\tau\left[\begin{array}{cc}A_{11}(s) & A_{12}(s) \\A_{21}(s) & A_{22}(s)\end{array}\right]\cdot\left[\begin{array}{cc}\mathbb{A}_{n-1}^{11}(s) & \mathbb{A}_{n-1}^{12}(s) \\\mathbb{A}_{n-1}^{21}(s) & \mathbb{A}_{n-1}^{22}(s)\end{array}\right]ds\\
&=\left[\begin{array}{cc}\int_0^\tau A_{11}(s)\cdot\mathbb{A}_{n-1}^{11}(s)ds & \int_0^{\tau}A_{11}(s)\cdot\mathbb{A}_{n-1}^{12}(s)ds +\int_0^\tau A_{12}(s)\mathbb{A}_{n-1}^{22}(s)ds\\\int_0^\tau A_{21}(s)\cdot\mathbb{A}_{n-1}^{11}(s)ds & \int_0^\tau A_{21}(s)\cdot\mathbb{A}_{n-1}^{12}(s)ds+\int_0^\tau A_{22}(s)\mathbb{A}_{n-1}^{22}(s)ds\end{array}\right].
\end{split}\]
Therefore, by induction,   we obtain that  the fundamental solution of the variational equation for the variables $(w_0,\mathsf{g}_{01},G_2,\mathsf{g}_{12})$ takes the form
\begin{equation}\label{fund-form1}\left[\begin{array}{cc}E^{11}(\tau)&  E^{12}(\tau)\\
\int_0^\tau A_{21}(s)E^{11}(s)ds& E^{22}(\tau)+\mathbb{I}\end{array}\right]+O(\frac{e^\tau}{\chi}),\end{equation}
where  $E^{11}(\tau)$ be the fundamental solution of the equation $\dot{Y}=A_{11}Y$, which is  the variational equation \eqref{variation-w0g0} for the variables $(w_0,\mathsf{g}_{01})$ and $E^{12}(\tau)$, $E^{22}(\tau)$ are of the form $\left[\begin{array}{cc}*&0\\
*&0\end{array}\right]$. 
This special form of the solution to the variational equation gives item (3) by definition of $E^{u/s}_{3}$ and  $E^{u/s}_{4}$. 
\end{proof}

\section{The derivative of the global map in the F4BP}\label{SG0g0}
In this section, we study the $C^1$ estimate of the global map in  the F4BP by allowing the variables $G_i,g_i$, $i=0,1,2$, to be non-vanishing. As we have explained in the beginning of the last section, we shall first work out the estimate of $d\mathbb G$ in the isosceles limit $Y\to0$, then handle the nonisoscelesness $Y\neq 0$ by the DuHamel principle. Unlike the local map case, where the off-diagonal terms caused by the nonisoscelesness can be neglected, the off-diagonal terms are estimated here as $O(|Y|\chi)$ and should be taken into account, since $\chi$ grows to infinity when we iterate.  The delicate part in the proof of Proposition \ref{PropDG2} is to prove the preservation of a three dimensional cone by the global map, although the map is nearly the identity in one direction of the cone. We have to show that this nearly identity direction is not spoiled by the $O(|Y|\chi)$ off-diagonal terms. For this purpose, some detailed information of the off-diagonal terms are needed. 
%Recall that we always restrict ourself to the case with zero total angular momentum,
%$G_0+G_1+G_2=0.$

First, in Section \ref{SSDiag}, we prove Lemma \ref{LmDiagonal} on the diagonal form of the derivative matrix. In Section \ref{SSGlob12}, we study the $C^1$-dynamics of the variables $G_i$, $g_i$, $i=0,1,2$ in the isosceles limit. In Section \ref{SSAMReduction}, we show how to reduce the variables $G_i$, $g_i$, $i=0,1,2$ to $(G_0,\mathsf g_{01}, G_2,\mathsf g_{12})$ and $(G_0,\mathsf g_{02}, G_1,\mathsf g_{12}).$ In Section \ref{non-isosceles-sec}, we study the derivative of the global map in the F4BP, taking into account the off-diagonal terms caused by the nonisoscelesness, and we complete the proof of Proposition \ref{PropDG2}. In Section \ref{SSTrans}, we give the proof of the transversality Lemma \ref{LmTrans} and Lemma \ref{LmOscG0}. %We  study only the variables $G_i,g_i$, $i=0,1,2$, while the others have already been treated in Section \ref{SDerLl}.    %Note that the above equality is invariant under the renormalization map defined in Definition \ref{renormalize-map}.

%As in Section \ref{SDerLl}, we study the variational equations evaluated along the I4BP orbits between the sections $\cS^L_{+}=\{ v^L=\eps^{-1/2}\}$ and $\cS^R_-=\{r^R=\frac{1}{\epsilon}\}$ and between the two sections, we introduce a middle section $\cS^M=\{-\beta \bx_1^L=\bx_2^R\}$.  We perform the change of coordinates, from the left Jacobi coordinates to the right one, on the section $\cS^M$.
\subsection{The block diagonalization}\label{SSDiag} We here prove Lemma \ref{LmDiagonal}.
\begin{proof}[Proof of Lemma \ref{LmDiagonal}]
By symmetry, we know that the subset $\{G_0=G_1=G_2=g_0-\pi/2=g_1=g_2=0\}$ is an invariant sub-manifold for the full 4-body problem. This is exactly the isosceles system we considered before. On this subset, we have the Hamiltonian equation $$\dot G_0=\dot G_1=\dot G_2=\dot g_0=\dot g_1=\dot g_2=0.$$
To show that the derivative in the statement is block-diagonalized, it is enough to show that the variational equation is block-diagonalized. We use formula \eqref{EqVar}. Let us consider the $JD^2_XH$ first.  Take the term $\frac{\partial^2 H}{\partial L_1\partial G_0}$ as an example. Since $\frac{\partial H}{\partial G_0}=0$ when evaluated on $\{G_0=G_1=G_2=g_0-\pi/2=g_1=g_2=0\}$, we still have $0$ when we take the derivative again with respect to $L_1$. Therefore, when restricted to the isosceles subsystem, the mixed derivative  $\frac{\partial^2 H}{\partial L_1\partial G_0}=0$.  In the same way, we show that the  other mixed derivatives $\frac{\partial^2 H}{\partial \xi \partial \eta}=0$, $\xi\in\{L_i,\ell_i,\;i=0,1,2\}$ and $\eta\in\{G_i,g_i,i=0,1,2\}$.   This proves the statement of the lemma.
\end{proof}

\subsection{The $C^1$-dynamics of the variables $G_i$, $g_i$, $i=0,1,2$ in the isosceles limit}\label{SSGlob12}

In this section, we compute the derivative of the global map for the variables $(G_0,g_0, G_1,g_1,G_2,g_2)$. We shall follow the general strategy outlined in Section \ref{SSStrategyLoc}. We break the orbit from $\cS_+^R$ to $\cS_-^L$ into two pieces by the middle section $\cS^M$. For each piece, we apply equation \eqref{eq: formald4} to compute the derivative, then we also compute the derivative of the coordinate change from left to right on the middle section $\cS^M$.  Since we have the Hamiltonian equations $\dot Y_i=J D_{Y_i}H=0$ along an orbit of the I4BP, where  $J=\left[\begin{array}{cc}0&-1\\
1&0\end{array}\right]$ and $Y_i=(G_i,g_i),\ i=0,1,2$, we get that the boundary contributions in \eqref{eq: formald4} are all the identity, so in the following, we only work on the fundamental solutions to the variational equation and the left to right transition. We remark that the contribution from the latter is essential for the $(G_1,g_1,G_2,g_2)$-component.

According to the strategy outlined in Section \ref{SSStrategyLoc}, as well as the discussion at the beginning of the current section, we split the proof into two steps:
(1) The variational equation,
(2) the transition from the left to the right.

{\bf Step 1,  the variational equation. }

We first consider the variational equation evaluated on orbits between the sections $\cS_+^L$ and $\cS^M$. Denoting $Y_0=(G_0,g_0)$, $Y_1=(G_1,g_1)$, $Y_2=(G_2,g_2)$ and $Y=(Y_1,Y_2)$, the Hamiltonian equation can be written as $\dot Y_0=JD_{Y_0}Y_0$, $\dot Y_1=JD_{Y_1}H,\ \dot Y_2=JD_{Y_2}H$ 
and the variational equation can be written as
\begin{equation}\label{VarEq-012}
\left[\begin{array}{c}
\dot{\dt Y}_0\\
\dot{\dt Y}_1\\
\dot{\dt Y}_2
\end{array}\right]=\left[\begin{array}{ccc}
JD_{Y_0}^2H&JD_{Y_0Y_1}^2H&JD_{Y_0Y_2}H\\
JD_{Y_1Y_0}^2H&JD_{Y_1}^2 H &JD_{Y_2Y_1}^2 H\\
JD_{Y_2Y_0}^2H&JD_{Y_1Y_2}^2 H&JD_{Y_2}^2 H
\end{array}\right]\left[\begin{array}{c}
\dt Y_0\\
\dt Y_1\\
\dt Y_2
\end{array}\right].
\end{equation}
%
%We notice the fact $L_2\sim \beta^{-1}$ so that by Lemma \ref{LmD2G} we have $$\frac{\partial x_2}{\partial G_2}=\frac{L_2}{M_2k_2}(\ell_2,0)=O(\beta\chi),\quad
%\frac{\partial^2 x_2}{\partial G_2^2}=\frac{1}{M_2k_2}(0, 2+\ell_2)=O(\beta^2\chi).$$
%Next recall the potentials \begin{equation}
%\begin{aligned}
%&k_0=1,\ k_1=2m_1,\ k_2=(m_1+2)m_2\\
%&U_{01}=\left(\frac{k_1}{|x_1|}-\frac{m_1}{|x_1-\frac{x_0}{2}|}-\frac{m_1}{|x_1+\frac{x_0}{2}|}\right),\\
%&U_2=\left(\frac{k_2}{|x_2|}-\frac{m_1m_2}{|x_2-\frac{2 x_1}{m_1+2} |}-\frac{m_2}{|x_2+\frac{m_1x_1}{m_1+2}-\frac{x_0}{2}|}-\frac{m_2}{|x_2+\frac{m_1x_1}{m_1+2}+\frac{x_0}{2}|}\right).
%\end{aligned}
%\end{equation}
From the estimates \eqref{U2G1s} and \eqref{u2-estimates},  we have the following estimate of the variational equation,
\begin{equation}\label{jdyh1}
\begin{aligned}
JD^2_{Y_1}H&=\left[\begin{array}{cc}
-\frac{\partial^2 (U_{01}+U_2)}{\partial G_1\partial g_1}&-\frac{\partial^2(U_{01}+U_2)}{\partial g^2_1}\\
\frac{\partial^2 (U_{01}+U_2)}{\partial G_1^2}&\frac{\partial^2(U_{01}+ U_2)}{\partial G_1\partial g_1}
\end{array}\right]=(\frac{1}{r_1^3}+\frac{r_1\chi}{(r_1+\chi)^3}) O\left[\begin{array}{cc} 1&1\\
1&1\end{array}\right],\\
J D^2_{Y_2} H&=\left[\begin{array}{cc}
-\frac{\partial^2 U_2}{\partial G_2\partial g_2}&-\frac{\partial^2 U_2}{\partial g^2_2}\\
\frac{\partial^2 U_2}{\partial G_2^2}&\frac{\partial^2 U_2}{\partial G_2\partial g_2}
\end{array}\right]=O\left[\begin{array}{cc}
\frac{\beta}{\chi}&\frac{1}{\chi}\\
\frac{\beta^2}{\chi}&\frac{\beta}{\chi}
\end{array}\right],\\
J D^2_{Y_2Y_1} H&=\left[\begin{array}{cc}
-\frac{\partial^2 U_2}{\partial G_2\partial g_1}&-\frac{\partial^2 U_2}{\partial g_2\partial g_1}\\
\frac{\partial^2 U_2}{\partial G_2\partial G_1}&\frac{\partial^2 U_2}{\partial G_1\partial g_2}
\end{array}\right]=O\left[\begin{array}{cc}
\frac{\beta}{\chi}&\frac{1}{\chi}\\
\frac{\beta}{\chi}&\frac{1}{\chi}
\end{array}\right],\quad\\
 J D^2_{Y_1Y_2} H&=\left[\begin{array}{cc}
-\frac{\partial^2 U_2}{\partial G_1\partial g_2}&-\frac{\partial^2 U_2}{\partial g_1\partial g_2}\\
\frac{\partial^2 U_2}{\partial G_1\partial G_2}&\frac{\partial^2 U_2}{\partial g_1\partial G_2}
\end{array}\right]=O\left[\begin{array}{cc}
\frac{1}{\chi}&\frac{1}{\chi}\\
\frac{\beta}{\chi}&\frac{\beta}{\chi}
\end{array}\right],\\
\end{aligned}
\end{equation}
and
\begin{equation}\label{jdyh2}
\begin{aligned}
JD_{Y_0}^2H&=\left[\begin{array}{cc}
-\frac{\partial^2 (U_{01}+U_2)}{\partial G_0\partial g_0}&-\frac{\partial^2(U_{01}+U_2)}{\partial g^2_0}\\
\frac{\partial^2 (U_{01}+U_2)}{\partial G_0^2}&\frac{\partial^2(U_{01}+ U_2)}{\partial G_0\partial g_0}
\end{array}\right]=(\frac{1}{r_1^3}+\frac{1}{\chi^3})O\left[\begin{array}{cc} 1&1\\
1&1\end{array}\right],\\
 JD_{Y_0Y_1}^2H&=\left[\begin{array}{cc}
-\frac{\partial^2 (U_{01}+U_2)}{\partial G_1\partial g_0}&-\frac{\partial^2(U_{01}+U_2)}{\partial g_0\partial g_1}\\
\frac{\partial^2 (U_{01}+U_2)}{\partial G_0\partial G_1}&\frac{\partial^2(U_{01}+ U_2)}{\partial G_0\partial g_1}
\end{array}\right]=(\frac{1}{r_1^3}+\frac{1}{\chi^3})O\left[\begin{array}{cc} 1&1\\
1&1\end{array}\right],\\
 J D^2_{Y_0Y_2} H&=\left[\begin{array}{cc}
-\frac{\partial^2 U_2}{\partial G_2\partial g_0}&-\frac{\partial^2 U_2}{\partial g_2\partial g_0}\\
\frac{\partial^2 U_2}{\partial G_2\partial G_0}&\frac{\partial^2 U_2}{\partial G_0\partial g_2}
\end{array}\right]=O\left[\begin{array}{cc}
\frac{1}{\chi^3}&\frac{1}{\chi^3}\\
\frac{1}{\chi^3}&\frac{1}{\chi^3}
\end{array}\right],
\end{aligned}
\end{equation}
Note that the estimates in \eqref{jdyh2} are so small that their contributions to the solution to the variational equation are negligible. 
We first consider the variational  equation 
\begin{equation}\label{VarEq-12}
\left[\begin{array}{c}
\dot{\dt Y}_1\\
\dot{\dt Y}_2
\end{array}\right]=\left[\begin{array}{cc}
JD_{Y_1}^2 H &JD_{Y_2Y_1}^2 H\\
JD_{Y_1Y_2}^2 H&JD_{Y_2}^2 H
\end{array}\right]\left[\begin{array}{c}
\dt Y_1\\
\dt Y_2
\end{array}\right].
\end{equation}
Between the sections $\cS_+^L$ and $\cS^M$ the quantity $r_1$ increases almost linearly with respect to time from the order $\epsilon^{-1}$ to the order $\chi$. Then we integrate the  equation \eqref{VarEq-12} over a time interval of order $\chi$, to get the following estimate of the fundamental solution
\begin{equation}\label{EqFundG1g1G2g2}\cD_G^L=\mathrm{Id}+O\left[\begin{array}{cccc}
1&1 &\beta&1\\
1& 1 &\beta &1\\
1&1&\beta&1\\
\beta&\beta&\beta^2&\beta
\end{array}\right].\end{equation}
For the 2 by 2 sub-matrix in the upper left corner, we have the following explicit form
\begin{equation}\label{explicit1}\frac{2m_1m_2}{m_1+2}A\cdot \left[\begin{array}{cc}\frac{1}{L_1}&1 \\
- \frac{1}{L_1^2}&- \frac{1}{L_1}\end{array}\right]+O(\epsilon^2),\end{equation}
where $L_1$ is from the initial condition on the section $\cS_+^L$ and  $A=\int_{\gamma(t)} \frac{r_1(t)r_2(t)}{(r_2(t)+\frac{2}{m_1+2}r_1(t))^3}-\frac{r_1(t)r_2(t)}{((r_2(t)-\frac{m_1r_1(t)}{m_1+2})^2+\frac{r_0^2(t)}{4})^{3/2}}dt$ is a negative  constant given by integrating the integral along an I4BP orbit $\gamma(t)$ \footnote{Note that for the orbit $\gamma(t)$, $r_0(t)\sim1$,  $r_1(t)\sim\epsilon^{-1}+\frac{k_1^L}{L_1^L}t$ and $r_2(t)\sim\chi+\beta t$. With these, we have that $A=\frac{m_1}{L_1^L}\Big(\frac{1}{(1+\frac{4m_1}{m_1+2}\frac{1}{L_1^L})^2}-\frac{1}{(1-\frac{2m^2_1}{m_1+2}\frac{1}{L_1^L})^2}\Big)+O(\epsilon^2).$ Nonetheless, we do not need this specified value of $A$. } between the sections $\cS^L_+$ and $\cS^M$. 
Indeed, to get this explicit form, we need only consider a linear ODE system with the coefficient matrix $\left[\begin{array}{cc}
-\frac{\partial^2 U_2}{\partial G_1\partial g_1}&-\frac{\partial^2 U_2}{\partial g^2_1}\\
\frac{\partial^2 U_2}{\partial G_1^2}&\frac{\partial^2  U_2}{\partial G_1\partial g_1}
\end{array}\right].$ Using the formulas from \eqref{U2G1s}, we have that the matrix takes the from
$ab(t)\left[\begin{array}{cc} \frac{1}{L_1}&1\\
- \frac{1}{L_1^2} &- \frac{1}{L_1} \end{array}\right]$. If we make a time change $d\tau=ab(t) dt$, then in the new time $\tau$, the variational equation has constant coefficients and the coefficient matrix has rank 1. Then the variational equation can be integrated explicitly as above. Moreover, by a similar calculation, we find that the $(1,4)$ entry is also $-L_1$ times the $(2,4)$ entry, as is the case in \eqref{explicit1}, so that to the leading order term, the first two rows of the matrix $\mathcal D^L_G-\mathrm{Id}$ are parallel. 

%\[\begin{split}A&=\int_{\gamma(t)} \frac{r_1(t)r_2(t)}{(r_2(t)+\frac{2}{m_1+2}r_1(t))^3}-\frac{r_1(t)r_2(t)}{((r_2(t)-\frac{m_1r_1(t)}{m_1+2})^2+\frac{r_0^2(t)}{4})^{3/2}}dt\\
%&=\frac{m_1}{L_1^L}\Big(\frac{1}{(1+\frac{4m_1}{m_1+2}\frac{1}{L_1^L})^2}-\frac{1}{(1+\frac{2m^2_1}{m_1+2}\frac{1}{L_1^L})^2}\Big)+O(\epsilon^2).
%\end{split}\]
%Note that all the terms involved $U_2$ have a common factor $m_2$. We have also the following expression,
%\begin{equation}\label{m2-dep}\cD_G^L=\mathrm{Id}+m_1O(\epsilon^2)+m_2O(1).\end{equation}

For the right case, that is between the sections $\cS^M$ and $\cS_-^R$,  the estimate for the coefficient matrix of the variational equation in \eqref{jdyh1} remains valid by replacing $\beta$ with $1$ and switching  the positions of the subindexes $1$ and $2$.   Note that the quantity $r_2$ decreases almost linearly with respect to time from the order $\beta\chi$ to the order $\epsilon^{-1}$.  Now we integrate the variational equation \eqref{VarEq-12} over a time interval of order $\beta\chi$. We then get the following estimate for its fundamental solution
\begin{equation}\label{D1-right}\cD_G^R=\mathrm{Id}+O\left[\begin{array}{cccc}\beta&\beta&\beta&\beta\\
\beta&\beta&\beta&\beta\\
\beta&\beta&\epsilon^2&\epsilon^2\\
\beta&\beta&\epsilon^2&\epsilon^2\end{array}\right].\end{equation}
With the estimates in \eqref{jdyh2}, using Picard iteration, we know the fundamental solution of the full variational equation \eqref{VarEq-012}  between the sections $\cS^L_+$ and $\cS^M$ takes the form
\begin{equation}\label{D-G-L}
\left[\begin{array}{cc}\mathbb{I}_{2\times2}+O_{2\times2}(\epsilon^2)& O_{2\times6}(\epsilon^2)\\
O_{6\times2}(\epsilon^2)&\mathcal D_{G}^L\end{array}\right],\; %\Delta_L=O_{2\times6}(\epsilon^2),\; \bar\Delta_L=O_{6\times2}(\epsilon^2),
\end{equation}
and between the sections $\cS^M$ and $\cS^R_-$ takes the form 
\begin{equation}\label{D-G-R}
\left[\begin{array}{cc}\mathbb{I}_{2\times2}+O_{2\times2}(\epsilon^2)& O_{2\times6}(\epsilon^2)\\
O_{6\times2}(\epsilon^2)&\mathcal D_{G}^R\end{array}\right],\; %\Delta_R=O_{2\times6}(\epsilon^2),\; \bar\Delta_R=O_{6\times2}(\epsilon^2).
\end{equation}
{\bf Step 2, Transition from the left to the right.}
\begin{Lm}\label{LmTransitionGg}
We have the following derivative estimate on the section $\cS^M=\{-\beta x_1^L=x_2^R\}$
$$\frac{\partial(G_0,g_0,G_1,g_1,G_2,g_2)^{R}}{\partial(G_0,g_0,G_1,g_1,G_2,g_2)^{L}}\Big|_{\cS^M}=\left[\begin{array}{cc}\mathbb{I}_{2\times2}&0\\
0&\mathbf{v}\otimes dG_1^R+O(1)\end{array}\right],$$
where $\mathbf{v}=(1,c_1,-1,c_2)^T$ with 
$c_1=-\frac{|\frac{\partial x_1^R}{\partial G_1}|}{|x_1^R|}=-\frac{1}{L_1^R}$, $c_2=\frac{|\frac{\partial x_2^R}{\partial G_2}|}{|x_2^R|}=\frac{1}{L_2^R}$ and
$dG_1^R:=\frac{\partial G_1^R}{\partial Y^L}\simeq\chi\cdot(\frac{1}{L_1^L},1,\beta,-1)$.
\end{Lm}
The physical meaning of this derivative matrix is clear. First, there is no coordinate change in the $(G_0,g_0)$-components, so we get an identity $\mathbb I_{2\times 2}$. The $O(\chi)$ rank-1 tensor part is caused by the coordinate change as follows. For instance, a $\Delta$-change in the angles $g_1^L,g_2^L$ will cause a change in position of order $\chi\Delta$ when arriving at the middle section, which changes the angular momentum by the same order. We refer readers to Appendix \ref{appendix-D} for detailed proofs. 

Define
\begin{equation}\label{ul-def}
\bar\bu_2= \cD_G^R\cdot\left[\begin{array}{cccc}
1\\
c_1\\
-1\\
c_2
\end{array}\right]=\left[\begin{array}{cccc}
1\\
c_1\\
-1\\
c_2
\end{array}\right]+O(\eps^2+\beta),\quad 
\bar\bl_2=\frac{1}{\chi} d G_1^R \cdot\cD_G^L=\frac{1}{\chi} d G_1^R+O(\eps^2+\beta),\end{equation}
where the last = follows from the fact that the first two rows of $\cD^L_G-\mathrm{Id}$ are parallel to the leading order term. Hence the derivative of the global map for the variables $G_i,g_i$, $i=0,1,2$, between the sections $\cS^L_+$ and $\cS^R_-$ takes the following form,
\begin{equation}\label{gs-tensor1}
\frac{\partial(G_0,g_0,G_1,g_1,G_2,g_2)^{R}\big|_{\cS^R_-}}{\partial(G_0,g_0,G_1,g_1,G_2,g_2)^{L}\big|_{\cS^L_+}}=\left[\begin{array}{cc}\mathbb{I}_{2\times2}&0\\
0&0\end{array}\right]+\chi\bar{\mathbf{u}}_2\otimes\bar{\mathbf{l}}_2.
\end{equation}
\subsection{The angular momentum reduction}\label{SSAMReduction}
Note that the  calculations in the previous subsection treat $(G_i,g_i)$, $i=0,1,2$ as independent variables. In this subsection, we exploit the total angular momentum conservation $G_0+G_1+G_2=0$ and the fact that the Hamiltonian depends only on the relative angles $g_i-g_j,\ i,j=0,1,2$ but not on any individual angle $g_i$. 

On the left section $\cS^L_+$ (near the triple collision of $Q_1$-$Q_3$-$Q_4$), we introduce the coordinates $\bar{Y}_L:=(G_0,\mathsf g_{01}, G_2,\mathsf g_{21}, G,g_1)$, and on the right section $\cS_-^R$ (near the triple collision of $Q_2$-$Q_3$-$Q_4$), we use the coordinates $\bar{Y}_R=(G_0,\mathsf g_{02},G_1,\mathsf g_{12}, G, g_2)$. We also introduce $Y_L=(G_0,g_0, G_2,g_2,G_1,g_1)^L$ and  $Y_R=(G_0,g_0, G_1,g_1,G_2,g_2)^R$. We thus find the matrix $\frac{\partial \bar{Y}_L}{\partial Y_L}=\left[\begin{array}
{cc|cccc}
1&0&0&0&0&0\\
0&1&0&0&0&-1\\
\hline
0&0&1&0&0&0\\
0&0&0&1&0&-1\\
1&0&1&0&1&0\\
0&0&0&0&0&1
\end{array}\right]:=\left[\begin{array}
{c|c}
\mathbb{I}_{2\times2}&A\\
\hline
B&C
\end{array}\right]$ whose inverse has the form
$\frac{\partial Y_R}{\partial \bar{Y}_R}=\left[\begin{array}
{cc|cccc}
1&0&0&0&0&0\\
0&1&0&0&0&1\\
\hline
0&0&1&0&0&0\\
0&0&0&1&0&1\\
-1&0&-1&0&1&0\\
0&0&0&0&0&1
\end{array}\right]:=\left[\begin{array}
{c|c}
\mathbb{I}_{2\times2}&A'\\
\hline
B'&C'
\end{array}\right],$
where we have $$A'=-A,\ B'=-B,\ C'-\mathbb{I}_{4\times4}=-(C-\mathbb{I}_{4\times4}).$$

Using the derivatives of the global map as in \eqref{gs-tensor1},
we compose
$$\frac{\partial \bar{Y}_R}{\partial {Y}_R} \frac{\partial Y_R}{\partial {Y}_L}\frac{\partial Y_L}{\partial \bar{Y}_L}=\left[\begin{array}
{c|c}
\mathbb{I}_{2\times2}&A\\
\hline
B&C
\end{array}\right]\left[\begin{array}
{c|c}
\mathbb{I}_{2\times2}&0\\
\hline
0&X
\end{array}\right]\left[\begin{array}
{c|c}
\mathbb{I}_{2\times2}&A'\\
\hline
B'&C'
\end{array}\right]=\left[\begin{array}
{c|c}
\mathbb{I}_{2\times2}+AXB'&A'+AXC'\\
\hline
B+CXB'&BA'+CXC'
\end{array}\right]$$
where $X=\chi \bar\bu_2\otimes \bar\bl_2$. 
In the resulting matrix, we eliminate the last two rows and last two columns corresponding to $G$ and $g_i,\ i=1,2$. Note that this elimination will trivialize the matrices $A',B',BA'$. Thus, we get \begin{equation}\label{EqG}\left[\begin{array}
{c|c}
\mathbb{I}_{2\times2}+AXB'&AXC'\\
\hline
CXB'&CXC'
\end{array}\right]=\left[\begin{array}
{c|c}
\mathbb{I}_{2\times2}&0\\
\hline
0&0
\end{array}\right]+\chi \bar\bu_2'\otimes \bar\bl_2',\end{equation}
where $\bar\bu_2'=(A\bar\bu_2, C\bar\bu_2)=(0,-c_2,1,c_1-c_2)$, and $\bar\bl_2'=(\bar\bl_2 B', \bar\bl_2 C')=(-\beta,0,\frac{1}{L_1^L}-\beta,1)$. 
In the resulting matrix, we eliminate the last two rows and last two columns corresponding to $G$ and $g_i,\ i=1,2$.
\subsection{The non-isosceles case}\label{non-isosceles-sec}
We next include the perturbations coming from the non-isoscelesness and complete the proof of Proposition \ref{PropDG2}.
\begin{proof}[Proof of Proposition \ref{PropDG2}]
Between the section $\cS^L_+$ and $\cS^M$, let us consider the non-isosceles orbits such that 
\begin{equation}\label{small-gs}|G_0|,|G_1|,|G_2|,|g_0-g_1-\frac{\pi}{2}|,|g_0-g_2+\frac{\pi}{2}|,|g_1-g_2-\pi|\leqslant \nu, \quad 0<\nu\ll1.
\end{equation}
From \eqref{D-G-L} we know that the above conditions are satisfied if they are valid on the section~$\cS^L_+$. 
We compute the  variational equation of the variables $(X,Y)$ where $X=(L_0,\ell_0,L_2,\ell_2)$ and  $Y=(G_0,g_0,G_1,g_1,G_2,g_2)$ along such an orbit. Recall that we reduce the variable $L_1$ using the conservation of energy and treat the variable $\ell_1\simeq r_1$ as the new time. With the estimates in Appendix \ref{mix-derivatives},   we then have in the leading term 
\begin{equation}\label{VarEq-XY}
\left[\begin{array}{c}\dot{\dt X}\\
\dot{\dt Y}\end{array}\right]=\left[\begin{array}{cc}K_{11}& K_{12}\\
K_{21}&K_{22}\end{array}\right]\left[\begin{array}{c}\dt X\\
\dt Y\end{array}\right],
\end{equation}
where $K_{11}$ is the coefficient matrix \eqref{coefficients-M1} of the variational equation for $X$, and $K_{22}$ is the coefficient matrix in \eqref{VarEq-012}.

{\bf Step 1.} {\it From the left section $\mathcal S_-^L$ to the middle section $\mathcal S^M$.}

 We have the estimate
\begin{equation}\label{estimate-K12} K_{12}\simeq\nu\left[\begin{array}{cccc|cc}\ & ( \frac{1}{r_1^3}+\frac{1}{\chi^3})_{2\times 4} &\ &\ & (\frac{1}{\chi^3})_{2\times 2} & \ \\
\hline
%\frac{1}{r_1^3}+\frac{1}{\chi^3} & \frac{1}{r_1^3}+\frac{1}{\chi^3} & \frac{1}{r_1^3}+\frac{1}{\chi^3} & \frac{1}{r_1^3}+\frac{1}{\chi^3} & \frac{1}{\chi^3} & \frac{1}{\chi^3}\\
\frac{1}{\chi^3} & \frac{1}{\chi^3} & \frac{r_1}{\chi^3} & \frac{r_1}{\chi^3} & \frac{
1}{\chi^2} & \frac{r_1}{\beta^{2}\chi^3}\\
\frac{1}{\chi^3} & \frac{1}{\chi^3} & \frac{\beta r_1}{\chi^2} & \frac{\beta r_1}{\chi^2} & \frac{\beta r_1}{\chi} & \frac{\beta r_1}{\chi^2} 
\end{array}\right]\simeq K_{21}^\prime.
\end{equation}
Let us denote by $\mathbb{X}(\ell_1)$ the fundamental solution of the equation $\dot{\dt X}=K_{11}\dt X$, and by $\mathbb{Y}(\ell_1)$ the fundamental solution of the equation $\dot{\dt Y}=K_{22}\dt Y$. The leading terms of $\mathbb{X}(\ell_1)$ and $\mathbb{Y}(\ell_1)$ are in \eqref{D-X-L} and \eqref{D-G-L}, respectively. We have the estimate $\mathbb{X}(\ell_1)=\ell_1 \bar \bu_1\otimes \bar\bl_1+O((\eps+\beta)\ell_1)$ by \eqref{D-X-L}  (see Proposition \ref{PropDG1}), and $\mathbb{Y}=O(1)$. We next show how to estimate the fundamental solution to \eqref{VarEq-XY}. From the first row of \eqref{VarEq-XY}, we get $\dot{\dt X}=K_{11} \dt X+K_{12}\dt Y$. Since we have $\frac{\dt X(0)}{\dt Y(0)}=0$, we get by the DuHamel principle
$\frac{\dt X(t)}{\dt Y(0)}=\int_0 ^t \mathbb X(t-s)K_{12}(s)\frac{\dt Y(s)}{\dt Y(0)}\,ds.$ %\quad \frac{\dt X(t)}{\dt X(0)}=\mathbb X(t)+\int_0 ^t \mathbb X(t-s)K_{12}(s)\frac{\dt Y(s)}{\dt X(0)}\,ds. $$
In this integral, the term $\frac{\dt Y(s)}{\dt Y(0)}$ is unknown, but we can replace it by $\mathbb Y(s)$ to get an estimate of $\frac{\dt X(t)}{\dt Y(0)}$ and finally verify the validity of this substitution for the purpose of doing estimates. For this purpose, we have to verify that in $ \frac{\dt Y(t)}{\dt Y(0)}=\mathbb Y(t)+\int_0 ^t \mathbb Y(t-s)K_{21}(s)\frac{\dt X(s)}{\dt Y(0)}\,ds$, the integral term is dominated by the leading term. 
Similarly, in the estimate of 
$\frac{\dt Y(t)}{\dt X(0)}=\int_0 ^t \mathbb Y(t-s)K_{21}(s)\frac{\dt X(s)}{\dt X(0)}\,ds$
we replace $\frac{\dt X(s)}{\dt X(0)}$ by $\mathbb X(s)$. Now comes an important difference in the estimate of the two terms $\frac{\dt X(t)}{\dt Y(0)}$ and $\frac{\dt Y(t)}{\dt X(0)}.$ We take $t=O(\chi)$ as the time going from the left section to the middle section. Then $\mathbb{X}(\chi-s)=(\chi-s) \bar \bu_1\otimes \bar\bl_1+$ lower order terms. In $K_{12}(s)$, the distance $r_1$ grows linearly with respect to $s$, so we get the estimate $\frac{\dt X(t)}{\dt Y(0)}=O(\nu\eps^2\chi)$ after the integration. However,  in the estimate of $\frac{\dt Y(t)}{\dt X(0)}$, we get that $\mathbb Y(t-s)$ is bounded, and the integral is convergent, so we get the estimate $\frac{\dt Y(t)}{\dt X(0)}=O(\nu)$. 

Finally, we have the following estimate for the fundamental solution of \eqref{VarEq-XY},
\begin{equation}\label{DGleft}\left[\begin{array}{cc}(1+O(\nu))\mathcal{D}_2^L&\mathcal{K}^L_{12}\\
O_{6\times 4}(\nu)&(1+O(\nu))\mathcal{D}_G^L\end{array}\right],
\end{equation}
where $\mathcal{D}_2^L$ is in \eqref{D-X-L}, $\mathcal{D}_G^L$ is in \eqref{EqFundG1g1G2g2}, and $\mathcal{K}^L_{12}=\nu \eps^2\chi \bar \bu_1\otimes \hat\bl_1+ O(\nu(\eps^3+\eps^2\beta)\chi)$
with $\hat\bl_1=\bar\bl_1\int_0^\chi  K_{12}(s)\mathbb Y(s)\,ds.$

%+ O\left[\begin{array}{cccccc}\epsilon^2 &\epsilon^2 &\epsilon^2 & \epsilon^2 & \beta &\beta \\
%\epsilon^2\chi & \epsilon^2\chi & \epsilon^2\chi & \epsilon^2\chi &\beta^3& \beta\\
%\frac{1}{\chi} & \frac{1}{\chi} & \beta & \beta& \beta & \beta \\
%\epsilon^2\beta^3\chi & \epsilon^2\beta^3\chi &\epsilon^2\beta^3\chi& \epsilon^2\beta^3\chi & \beta^4
 %&\beta^2\end{array}\right]=O(\nu\eps^2\chi).\]

{\bf Step 2.} {\it From the middle section to the right.}
 
Consider  orbits between the section $\cS^M$ and $\cS^R_-$ satisfying the condition \eqref{small-gs}. From \eqref{D-G-R} we know that this condition is true if it is true on the section $\cS^M$. Recall that now $X=(L_0,\ell_0,L_1,\ell_1)$. We reduce the variable $L_2$ using the conservation of energy, and treat the variable $\ell_2\simeq r_2$ as the new time; $\ell_2$ decreases from $O(\beta\chi)$ on  $\cS^M$ to $O(\epsilon^{-1})$ on  $\cS^R_-$.   We still denote the variational equation as \eqref{VarEq-XY}, and now the estimate for the entries in $K_{12}$ and $K_{21}$ are given by
\begin{equation}\label{estimate-K12-R} K_{12}\simeq\nu\left[\begin{array}{cccc|cc}\frac{1}{r_2^3}+\frac{1}{\chi^3} & \frac{1}{r_2^3}+\frac{1}{\chi^3} & \frac{1}{\chi^3} & \frac{1}{\chi^3}& \frac{1}{r_2^3}+\frac{1}{\chi^3} & \frac{1}{r_2^3}+\frac{1}{\chi^3} \\
\frac{1}{r_2^3}+\frac{1}{\chi^3} & \frac{1}{r_2^3}+\frac{1}{\chi^3} &\frac{1}{\chi^3} & \frac{1}{\chi^3} & \frac{1}{r_2^3}+\frac{1}{\chi^3} & \frac{1}{r_2^3}+\frac{1}{\chi^3}\\
\hline
\frac{1}{\chi^3} & \frac{1}{\chi^3} & \frac{r_2}{\chi^3} & \frac{r_2}{\chi^3} & \frac{
1}{\chi^2} & \frac{r_2}{\chi^3}\\
\frac{1}{\chi^3} & \frac{1}{\chi^3} & \frac{r_2^2}{\chi^3} & \frac{r_2^2}{\chi^3} & \frac{r_2^2}{\chi^3} & \frac{ r_2^2}{\chi^3} \end{array}\right]\simeq K_{21}^\prime.
\end{equation}
By the same argument as before, we have the following form for the fundamental solution of the corresponding variational equation,
\begin{equation}\label{DGright}\left[\begin{array}{cc}(1+\nu)\mathcal{D}_3^R&O(\nu\eps)\\
{\mathcal{K}}_{21}^R&(1+O(\nu))\mathcal{D}_G^R\end{array}\right],\end{equation}
where $\mathcal{D}_3^R$ is in \eqref{D-X-R}, $\mathcal{D}_G^R$ is in \eqref{D1-right},  and $\mathcal{K}_{21}^R=O(\nu\eps^2\beta\chi).$ Note that the estimates of the off-diagonal blocks in \eqref{DGright} are quite different from that in \eqref{DGleft}. In particular, the $\nu\chi$-dependent block in \eqref{DGleft} is the (1,2)-block while that in \eqref{DGright} is the (2,1)-block. The difference is caused by the fact that $r_1$ is increasing linearly in $t$ while $r_2$ is decreasing linearly in $t$. 

% \[\mathcal{K}_{21}^R\simeq\int_{\beta\chi}^{\frac{1}{\epsilon}}K_{21}(\ell_2)\mathbb{X}(\ell_2)d\ell_2\simeq\nu\left[\begin{array}{cccc}\epsilon^4\beta\chi&\epsilon^2&\epsilon^4\beta\chi&\frac{1}{\chi}\\
%\epsilon^2\beta\chi&\epsilon^2&\epsilon^2\beta\chi&\frac{1}{\chi}\\
%\beta^2&\frac{1}{\chi}&\beta^3&\frac{1}{\chi}\\
%\beta^2&\frac{1}{\chi}&\beta^3&\frac{1}{\chi}\\%
%\epsilon^4\beta\chi&\epsilon^2&\epsilon^4\beta\chi&\frac{1}{\chi}\\
%\epsilon^2\beta\chi&\epsilon^2&\epsilon^2\beta\chi&\frac{1}{\chi}
 %\end{array}\right]=O(\nu\eps^2\beta\chi).\]
 
 {\bf Step 3.} {\it The transition from left to right.}
 
We next consider  the transition from the left to right on the section $\cS^M$, assuming 
\begin{equation}\label{small-g-L}|G^L_1|,|G^L_2|,|g^L_1-g^L_2-\pi|=O(\nu),\end{equation}
and after the transition,
\begin{equation}\label{small-g-R}|G^R_1|,|G^R_2|,|g^R_1-g^R_2-\pi|=O(\nu_0),
\end{equation}
then we have 
\begin{equation}\label{non-iso-transit}\frac{\partial (X^R,Y^R)}{\partial(X^L,Y^L)}\Big|_{\cS^M}=\left[\begin{array}{c|cc}\mathcal{M}_{4\times4}&\Big(\begin{array}{c}0\\O_{2\times2}(\frac{\nu}{\chi^2})\end{array}\Big)&\Big(\begin{array}{c}0\\ 
\nu_0\chi (O(1)_{2\times1})\otimes dG_1^R\end{array}\Big)\\
\hline
0&\mathbb{I}_{2\times2}&0\\
\Big(\begin{array}{ccc}O_{4\times2}(\frac{\nu}{\chi^2})&O_{4\times2}(\nu\chi)\end{array}\Big)&O_{4\times2}(\frac{\nu}{\chi^2})&\chi\mathbf{v}\otimes dG_1^R\end{array}\right]
\end{equation}
where $\mathcal{M}$ is the transition matrix for the $X^{L,R}$  given \eqref{middle-d-l} and $\mathbf{v}\otimes dG_1^R$ is in Lemma \ref{LmTransitionGg}. The proof of the above form is presented in Appendix \ref{appendix-D}.

 {\bf Step 4.} {\it Completing the proof of Proposition \ref{PropDG2}.}

Composing the resulting matrices in the above Step 1-Step 3, we get the following estimate of the derivatives of the global map
$\frac{\partial (X^R,Y^R)}{\partial(X^L,Y^L)}=\chi \bar\bu_1\otimes \bar\bl_1+ \chi \bar{\mathbf{u}}_2\otimes\bar{\mathbf{l}}_2+O(\nu\chi)$
where $\mathbf{u}_1=(0,1,0,1)$, $\mathbf{l}_1=(1,0,0,0)$ and $\bar\bu_2$, $\bar\bl_2$ are defined in \eqref{ul-def}.

This gives the tensor structure decomposition as stated in Proposition \ref{PropDG2}. We next show the cone preservation property \eqref{EqDG}. For this purpose, we pick a vector $ \bu\in E^u_{4}(\gamma_O(0))$ that is orthogonal to $\bar\bl_2$. We then apply the three matrices composing $d\mathbb G$ one by one, $d\mathbb G=\eqref{DGright} \eqref{non-iso-transit}\eqref{DGleft}$. Note that the first four components of $\bu$ vanish, corresponding to the $X$-variables. After applying \eqref{DGleft}, we get a resulting vector of the form $(\nu \eps^2\chi \bar \bu_1,0_{1\times 6})+(0_{1\times 4},(1+O(\nu))\mathcal D_G^L \bu)$. The image of $(\nu \eps^2\chi \bar \bu_1,0_{1\times 6})$ after applying \eqref{non-iso-transit} and \eqref{DGright} is a vector in $\mathcal C_\eta(\bar\bu_1)$, so it remains to consider $(0_{1\times 4},(1+O(\nu))\mathcal D_G^L \bu)$. From the multiplication of \eqref{non-iso-transit} and \eqref{DGleft} we find that $\bar\bl_2=dG_1^R(1+O(\nu))\mathcal D_G^L$, so from the choice of $\bu$, we see that in the application of \eqref{non-iso-transit} to $(0_{1\times 4},(1+O(\nu))\mathcal D_G^L \bu)$, the $O(\chi)$ and $O(\nu_0\chi)$ terms in \eqref{non-iso-transit} disappear and we obtain a vector of the form $(O(\nu/\chi^2)_{1\times 4}, \bu'_{1\times 2},O(\nu/\chi^2)_{1\times 4})$, where  $\bu'$ is the  projection of $(1+O(\nu))\mathcal D_G^L \bu$  to the $(G_0,\mathsf g_{01})$-component. After the further application of \eqref{DGright}, we get a vector of the form $\bar \bu_3:=(O(\nu\eps)_{1\times 4}, \bu''_{1\times 2},O(\nu)_{1\times 4})$ where $\bu''$ is the projection of $(1+O(\nu))\mathcal{D}_G^R (\bu',O(\nu/\chi^2)_{1\times 4})$  to the $(G_0,\mathsf g_{01})$-component. This is also the  projection of $\bu$  to the $(G_0,\mathsf g_{01})$-component up to  a small error. This completes the proof of the proposition.

\end{proof}
\subsection{Proof of the transversality Lemma \ref{LmTrans} and Lemma \ref{LmOscG0}}\label{SSTrans}
We next give the proof of the transversality Lemma \ref{LmTrans}
\begin{proof}[Proof of Lemma \ref{LmTrans}]
From Step 4 of the last subsection, we find that $\bar \bu_3$ is defined as follows. We first select a vector in $E_{4}^u(\gamma_O(0))$ that is perpendicular to $\bar\bl_2$, which is uniquely determined up to a nonzero scalar multiple. Indeed, the projection of $\bar\bl_2$ to $E_{4}^u(\gamma_O(0))$ is not a point, since we have $\bar\bl_2\cdot \mathbf e_1\neq 0$ (c.f. $\bar\bl$ \eqref{EqG} for $\bar\bl_2$ and \eqref{EqEigenvector} for $\mathbf e_1$). 
For this vector, we project to the $(G_0,\mathsf{g}_{02})$ component to get $\bar\bu_3$. This projection does not kill the vector either since $\bar\bl_2$ has the first two entries vanishing (see \eqref{EqG}) and $\mathbf e_1$ has the first two entries nonvanishing. 

Note that the transversality condition is formulated in the isosceles limit. The orbits $\gamma_I$ and $\gamma_O$ are analytic functions of the masses $m_1$ and $m_2$ ($\gamma_I$ depends on both masses and $\gamma_O$ depends on only one). So are the tangent spaces $E^{u/s}_{4}$. Both the vectors $\bar{\mathbf u}_2$ and $\bar{\mathbf u}_3$ are analytic functions of the masses by their definitions. We thus see that the transversality condition is equivalent to the nonvanishing of the determinant, which is an analytic function of the masses $m_1,m_2$. Then the statement follows from Lemma \ref{LmUlam}. 
\end{proof}
We complete this section with the proof of Lemma \ref{LmOscG0}.

\begin{proof}[Proof of Lemma \ref{LmOscG0}] We again cut the orbit into two pieces, between  $\cS_+^L$ and $\cS^M$ and between $\cS^M$ and $\cS^R_-$. 
 For the piece of orbit $\gamma(t)$ between $\cS^L_+$ and $\cS^M$ with $\gamma(0)\in\cS_+^L$, assume we have 
 \begin{equation}\label{assume-bound}|G_0(t)|, \; |G_1(t)|\;|G_2(t)|, \; |g_0-g_1-\frac{\pi}{2}|,\; |g_1-g_2-\pi|\leqslant C_*|G_0(0)|\ll1,\end{equation}
 where $C_*\ll\epsilon^{-1}$ is a constant to be determined later. 
 Then the equation of motion for $G_0$ takes the form
 \begin{equation}\label{G0-1/2}\dot G_0=-\frac{\partial U_{01}+U_2}{\partial g_0}=O\Big((\frac{1}{r_1^3}+\frac{1}{\chi^3})C_*|G_0(0)|\Big).\end{equation}
 Therefore we have $|G_0(t)-G_0(0)|<\frac{1}{2}|G_0(0)|.$ We then consider the equations of motion for $G_2$ and $\mathsf{g}_{12}=g_1-g_2-\pi$. 
 With the estimate in Appendix \ref{mix-derivatives} and the conservation of total angular momentum $G_0+G_1+G_2=0$, we have 
 $$\begin{cases}\dot{\mathsf{g}}_{12}=c_{11}(t)G_2+c_{12}(t)\mathsf{g}_{12}+b_1(t)|G_0(0)|,\\
\dot G_2=c_{21}(t)G_2+c_{22}(t)\mathsf{g}_{12}+b_2(t)|G_0(0)|,\end{cases}$$
  where $|c_{ij}(t)|, \;|b_i(t)|\leqslant \frac{tC_1}{\chi^2}$ with $C_1$ depending only on $m_1$, $m_2$. 
  Therefore, $|G_2(t)|, \;|\mathsf{g}_{12}|\leqslant C_2|G_0(0)|,$ where $C_2$  depends only on  $m_1$, $m_2$. 
   Plugging  the obtained estimates for $\mathsf{g}_{12}$, $G_2$, $G_0$ into the equation of motion for $\mathsf{g}_{01}=g_0-g_1-\frac{\pi}{2}$, we have $|g_{01}(t)|\leqslant C_3|G_0(0)|,$
   with $C_3$ depending only on  $m_1$, $m_2$.  Therefore we prove \eqref{G0-1/2} for this piece of the orbit if we choose in the assumption~\eqref{assume-bound} $C_*>\max\{C_2,C_3\}$. 
   
   Similarly, we show that \eqref{G0-1/2} remains true for the piece of orbit between the sections $\cS^M$ and $\cS^R_-$ if we consider $\cS^R_-$ as the initial section and take the equations of motions backward in time.
   Thus, we complete the proof of Lemma~\ref{LmOscG0}.
\end{proof}
%
%      \begin{figure}[ht]
%\begin{subfigure}{0.48\textwidth}
%\includegraphics[height=4.5cm,width=7cm]{gamma_I_backward1}
%\end{subfigure}
%\begin{subfigure}{0.48\textwidth}
%\includegraphics[height=4.5cm,width=7cm]{gamma_I_forward1}
%\end{subfigure}
%\begin{subfigure}{0.48\textwidth}
%\includegraphics[height=4.5cm,width=7cm]{gamma_I_backward2}
%\end{subfigure}
%\begin{subfigure}{0.48\textwidth}
%\includegraphics[height=4.5cm,width=7cm]{gamma_I_forward2}
%\end{subfigure}
%\begin{subfigure}{0.48\textwidth}
%\includegraphics[height=4.5cm,width=7cm]{gamma_I_backward3}
%\end{subfigure}
%\begin{subfigure}{0.48\textwidth}
%\includegraphics[height=4.5cm,width=7cm]{gamma_I_forward3}
%\end{subfigure}
%
%\caption{The left: backward solution along $\gamma_I^2$; the right: forward solution along $\gamma_I^2$}
%\end{figure}
%

\appendix
\section{The triple collision blow up in the F3BP}\label{app-blowup}
In this section,  we perform a triple collision blow up to study the dynamics of the 3-body problem in a neighborhood of the isosceles case.
We will focus on the blow-up of the triple $Q_1$-$Q_3$-$Q_4$. The same results hold for the other triple $Q_2$-$Q_3$-$Q_4$.
Let us introduce
\[\bx=(M_0^{1/2}\bx_0,M_1^{1/2}\bx_1)\in\mathbb{R}^4, \quad \bp=(M_0^{-1/2}\bp_0,M_1^{-1/2}\bp_0)\in\mathbb{R}^4,\]
Where the reduced masses $M_i$ are defined in (5.1).
Then we can rewrite the Hamiltonian $H_\therefore$ in  \eqref{EqHamLoc} as $H_\therefore(\bx,\bp)=\frac{1}{2} |\bp|^2+V (\bx).$
The Hamiltonian equation is $\begin{cases}\dot{\bx}=\bp,\\ \dot{\bp}=-\partial_{\bx}H_\therefore=-\partial_\bx V (\bx)\end{cases}.$
 We define the {\it blow-up   coordinates near triple collision}  as:
\begin{equation}\label{EqBlowup1}
\begin{aligned}
\begin{cases}
r&=|\bx|,\\
\mathbf{s}&=\bx/r,\\
v&=r^{1/2}(\bs\cdot\bp)\\%=r^{1/2}(R_0r_0+R_1r_1),\\
\bw&=r^{1/2}(\bp-(\bs\cdot \bp)\bs).
\end{cases}
\end{aligned}
\end{equation}
Note from the definition that
$|\bs|=1,\quad \bs\cdot\bw=0.$ We also scale the time variable of the system by $dt=r^{3/2}d\tau,$ and use $'$ to denote $\frac{d}{d\tau}$.
The equations of motion read as
\begin{equation}
\begin{aligned}
\begin{cases}
r'&=rv,\\
\bs'&=\bw,\\
v'&=|\bw|^2+\frac{1}{2}v^2+\bar V(\bs),\\
\bw'&=-\frac{1}{2}v\bw-|\bw|^2\bs-\nabla \bar V(\bs)-\bar V(\bs)\bs.
\end{cases}
\end{aligned}
\end{equation}
Here $\bar V(\bs)=r V (\bx)$, $\nabla \bar V(\bs)=r^2\nabla_{\bx} V (\bx)$. The corresponding energy relation is 
\begin{equation}\label{energy1}rE=\frac{1}{2}|\bw|^2+\frac{1}{2} v^2+\bar V(\bs), \end{equation}
where $E$ is the constant value of the Hamiltonian.

Recalling the Jacobi-polar coordinates defined in subsection \ref{def-jacobi-polar}, we  introduce
\[\bs=(\cos\psi(\cos\theta_0,\sin\theta_0),\sin\psi(\cos\theta_1,\sin\theta_1)),\]
where $\psi=\arctan\frac{M_1^{1/2}r_1}{M_0^{1/2}r_0}$, and  the orthonormal basis of the orthogonal complement of $\bs$ in $\mathbb{R}^4$,
\[\begin{cases}
\be_w=(-\sin\psi(\cos\theta_0,\sin\theta_0),\cos\psi(\cos\theta_1,\sin\theta_1)),\\
\be_0=(-\sin\theta_0,\cos\theta_0,0,0),\\
\be_1=(0,0,-\sin\theta_1,\cos\theta_1).
\end{cases}
\]
Then we have $\bw=w\be_w+w_0\be_0+w_1\be_1,$
where
\[\begin{cases}
w=\bw\cdot \be_w=r^{1/2}\bp\cdot\be_w=r^{1/2}(-\frac{R_0\sin\psi}{M_0^{1/2}}+\frac{R_1\cos\psi}{M_1^{1/2}}),\\
w_0=\cos\psi\bw\cdot\be_0= \cos\psi r^{1/2}\bp\cdot\be_0=r^{-1/2}G_0,\\
w_1=\sin\psi\bw\cdot\be_1= \sin\psi r^{1/2}\bp\cdot\be_1=r^{-1/2}G_1.
\end{cases}\]
Clearly, we have $|\bw|^2=w^2+\frac{w_0^2}{\cos^2\psi}+\frac{w^2_1}{\sin^2\psi}$. Direct calculation gives us the following equations of motion for the full 3-body problem,
\begin{equation}\label{blowup-ws}
\begin{cases}
r'=rv,\\
v'=w^2+\frac{w_0^2}{\cos^2\psi}+\frac{w^2_1}{\sin^2\psi}+\frac{1}{2}v^2+\bar V(s),\\
\psi'=w,\\
w'=-\frac{1}{2}wv-\bar\nabla V(\bs)\cdot\be_w-w_0^2\frac{\sin\psi}{\cos^3\psi}+w_1^2\frac{\cos\psi}{\sin^3\psi},\\
w_0'=-\frac{1}{2}vw_0-r\partial_{g_0}U_{01},\\
w_1'=-\frac{1}{2}vw_1-r\partial_{g_1}U_{01},\\
g_0'=r^{3/2} \frac{\partial U_{01}}{\partial G_0}=\frac{r^{3/2}m_1}{2}\big(-\frac{\langle \bx_1,\partial_{G_0} \bx_0\rangle}{|\bx_1-\frac{\bx_0}{2}|^3}+\frac{\langle \bx_1,\partial_{G_0} \bx_0\rangle}{|\bx_1+\frac{\bx_0}{2}|^3}\big),\\
g_1'=\frac{r^{3/2}}{2} \frac{\partial U_{01}}{\partial G_1}=r^{3/2}m_1\big(\frac{\langle \partial_{G_1}\bx_1, \bx_1\rangle}{|\bx_1|^3}+\frac{\langle \partial_{G_1}\bx_1, \bx_1-\frac{\bx_0}{2}\rangle}{|\bx_1-\frac{\bx_0}{2}|^3}+\frac{\langle \partial_{G_1}\bx_1, \bx_1+\frac{\bx_0}{2}\rangle}{|\bx_1+\frac{\bx_0}{2}|^3}\big).\end{cases}
\end{equation}

\section{Delaunay coordinates}\label{App-delaunay}
In this appendix, we give the relations between Cartesian coordinates and Delaunay coordinates. Consider the Kepler problem $H(P,Q)=\frac{|P|^2}{2m}-\frac{k}{|Q|},$ $(P,Q)\in \R^2\times \R^2$.
\subsection{Elliptic motion}
We have the following relations which explain the physical and geometrical meaning of the Delaunay coordinates.
$a=\frac{L^2}{mk}$ is the semimajor axis, $b=\frac{LG}{mk}$ is the semiminor axis, $E=-\frac{k}{2a}$ is the energy, $G=Q\times P$ is the angular momentum, and $e=\sqrt{1-\left(\frac{G}{L}\right)^2}$ is the eccentricity. Moreover, $g$ is the argument of apapsis and $\ell$ is the mean anomaly.  We can relate $\ell$ to the polar angle $\psi$ through
the equations
$\tan\frac \psi 2 = \sqrt{\frac{1+e}{1-e}}\cdot\tan\frac u 2,\quad u-e\sin u=\ell.$
We also have Kepler's law $\frac{a^3}{T^2}=\frac{k/m}{(2\pi)^2}$ which relates the semimajor axis
$a$ and the period $T$ of the elliptic motion.

Denoting the body's position by $Q=(q_1, q_2)$ and its momentum by $P=(p_1,p_2)$ we have the following formulas
in the case $g=0$
\begin{equation}
\label{DelEll}
\begin{aligned}
q_1=\frac{L^2}{mk}\left(\cos u-\sqrt{1-\frac{G^2}{L^2}}\right), \quad &
q_2=\frac{LG}{mk} \sin u,
\\
p_1=-\frac{mk}{L}\frac{\sin u}{1-\sqrt{1-\frac{G^2}{L^2}} \cos u}, \quad &
p_2=\frac{mk}{L^2}\frac{G\cos u}{1-\sqrt{1-\frac{G^2}{L^2}}\cos u},
\end{aligned}
\end{equation}
where $u$ and $l$ are related by $u-e\sin u=\ell$. This is an ellipse whose major axis lies on the horizontal axis, with one focus at the origin and the other focus on the negative horizontal axis. The body moves counterclockwise if $G>0$. 
%
%\[e_{G}=\frac{-\frac{G}{L^2}}{\sqrt{1-(\frac{G}{L})^2}}=\frac{-\frac{G}{L^2}}{e}, e_{GG}=-\frac{1}{L^2e}+\frac{\frac{G}{L^2}}{e^{2}}e_G.\]
%So we have $e_G|_{G=0}=0$, $e_{GG}|_{G=0}=-\frac{1}{L^2}$.
%\[u_{G}-e_{G}\sin u-e\cos u u_{G}=0.\]
%\[u_{GG}-e_{GG}\sin u-e_{G}\cos u u_G-e_{G}\cos u u_{G}+e\sin u (u_G)^2-e\cos u u_{GG}=0.\]
%So we have $u_G|_{G=0}=0$, $(u_{GG}+\frac{1}{L^2}\sin u- u_{GG}\cos u)|_{G=0}=0$.
%
%We have the partial derivative calculations
%\begin{equation}
%\begin{aligned}
%u_G(1-e\cos u)&=e_G\sin u,\quad u_G|_{G=0}=e_G|_{G=0}=0,\\
%u_{GG}(1-e\cos u)|_{G=0}&=e_{GG}\sin u|_{G=0},\quad
%e_{GG}|_{G=0}=-\frac{1}{L^2}.
%\end{aligned}
%\end{equation}
\subsection{Hyperbolic motion}
For Kepler hyperbolic motion we have similar expressions for the geometric and physical quantities
\[a=\frac{L^2}{mk}, \ b=\frac{LG}{mk},\ E=\frac{k}{2a},\ G=Q\times P,\ e=\sqrt{1+\left(\frac{G}{L}\right)^2}\]
as well as the parametrization of the orbit
\begin{equation}
\begin{aligned}
q_1=\frac{L^2}{mk}
\left(\cosh u-\sqrt{1+\frac{G^2}{L^2}}\right), \quad
& q_2=\frac{LG}{mk} \sinh u, \\
p_1=-\frac{mk}{L}\frac{\sinh u}{1-\sqrt{1+\frac{G^2}{L^2}} \cosh u}, \quad
& p_2=-\frac{mk}{L^2}\frac{G\cosh u}{1-\sqrt{1+\frac{G^2}{L^2}}\cosh u}.
\end{aligned}\label{eq: delaunay4}
\end{equation}
where $u$ and $\ell$ are related by
\begin{equation}u-e\sinh u=\ell, \text{ and } e=\sqrt{1+\left(\frac{G}{L}\right)^2}.
\label{eq: hypul}
\end{equation}
%Differentiating this relation, we have
%\[u_{G}-e_{G}\sinh u-e\cosh u u_{G}=0, u_{GG}-e_{GG}\sinh u-e_{G}\sinh u u_{G}-e_{G}\cosh u u_G-e\sinh u (u_G)^2-e\cosh u u_{GG}=0.\]
%So we have
%\begin{equation}
%e_G|_{G=0}=0, e_{GG}=\frac{1}{L^2},\
%u_G(1-\cosh u)|_{G=0}=0,\quad (u_{GG}-\frac{1}{L^2}\sinh u-u_{GG}\cosh u)|_{G=0}=0.\end{equation}
This hyperbola is symmetric with respect to the $x$-axis, open to the right, and the particle moves counterclockwise on it when $u$ increases ($l$ decreases) in the case when minus the angular momentum $G=Q\times P>0$.
The angle $g$ is defined to be the angle measured from the positive $x$-axis to the symmetric axis. There are two such angles that differ by $\pi$, depending on the orientation of the symmetric axis. %This $\pi$ difference disappears in  the hamiltonian equation after taking derivative so it does not matter which angle we choose. Here $g$ does not appear in \eqref{DelEll} and \eqref{eq: delaunay4} because the argument of apoapsis or perigee is chosen to be zero or $\pi$. 
In the general case, we need to rotate the $(q_1,q_2)$ and $(p_1,p_2)$ using the matrix $\mathrm{Rot}(g)=\left[\begin{array}{cc}\cos g & -\sin g \\\sin g & \cos g\end{array}\right]$.
To parameterize the orbits in the main body of the paper, we use the following convention. In Jacobi coordinates \eqref{EqJacobi}, for the binary, we introduce $(\bx_0,\bp_0)=\mathrm{Rot}(g_0)(Q,P)$, where $(Q,P)$ is given in \eqref{DelEll} (we shall give subscript $0$ to the Delaunay variables). For $(\bx_1,\bp_1)$ and $(\bx_2,\bp_2)$, we note that $\bx_1$ points to the left and $\bx_2$ points to the right, so we introduce the parameterizatoin $(\bx_1,\bp_1)=\mathrm{Rot}(g_1+\pi)(Q,P)$ for $(Q,P)$ in \eqref{eq: delaunay4} with subscript $1$ for the Delaunay variables, and $(\bx_2,\bp_2)=\mathrm{Rot}(g_2)(Q,P)$ for $(Q,P)$ in \eqref{eq: delaunay4} with subscript $2$ for the Delaunay variables. Then the isosceles case corresponds to $G_0=G_1=G_2=g_0-\pi/2=g_1=g_2=0. $

We have the following lemma describing the derivatives of the polar variables  with respect to the Delaunay variables.
\begin{Lm}\label{Lm-A1} For  elliptic motion, we have
$ r_0=\frac{1}{M_0 k_0} L_0^2(1-\cos u_0), \quad u_0-\sin u_0=\ell_0,$
$$\frac{\partial r_0}{\partial\ell_0}=\frac{\partial r_0}{\partial u_0}\frac{d u_0}{d\ell_0}=\frac{L_0^2}{M_0k_0}\frac{\sin u_0}{1-\cos u_0}=\left(\frac{L_0^2}{M_0k_0}\right)^2\frac{\sin u_0}{r_0},\quad \frac{\partial^2 r_0}{\partial \ell_0^2}=O(\frac{1}{r_0^2}).$$
For hyperbolic motion, we have
$ r_1=\frac{1}{M_1k_1}L_1^2(\cosh u_1-1),\quad u_1-\sinh u_1=\ell_1.$
$$\frac{\partial r_1}{\partial\ell_1}=\frac{\partial r_1}{\partial u_1}\frac{d u_1}{d\ell_1}=\frac{L_1^2}{M_1k_1}\frac{\sinh u_1}{1-\cosh u_0}=\left(\frac{L_1^2}{M_1k_1}\right)^2\frac{-\sinh u_1}{r_1}.$$
If $\bx_0$ is elliptic and $\bx_1$ is hyperbolic we have
$$v=r^{-1/2}\bx\cdot \bp=r^{-1/2}(-L_0 e_0\sin u_0-L_1 e_1\sinh u_1).$$
%\begin{enumerate}
%\item In the left case we have $$\frac{\partial r_2}{\partial \ell_2}=L_2^2=\beta^{-2},\quad \frac{\partial r_2}{\partial L_2}=2L_2\ell_2=\beta\chi,\quad \frac{\partial r_1}{\partial \ell_1}=L_1^2=O(1),\quad \frac{\partial r_1}{\partial L_1}=2L_1\ell_1.$$
%\item In the right case we have
 %$$\frac{\partial r_2}{\partial \ell_2}=L_2^2=O(1),\quad \frac{\partial r_2}{\partial L_2}=2L_2\ell_2,\quad \frac{\partial r_1}{\partial \ell_1}=L_1^2=O(1),\quad \frac{\partial r_1}{\partial L_1}=2L_1\ell_1=O(\chi).$$
%\end{enumerate}
\end{Lm}
We cite the next lemma from \cite{DX} to simplify the calculation of the derivative with respect to $\ell$.
\begin{Lm}[Lemma A.2 of \cite{DX}] \label{Lm: simplify}
Let $u$ be the function of $\ell, G$ and $L$ given by \eqref{eq: hypul} and let $\sigma=sign(u)$ when $|u|$ is large.
Then we can approximate $u$
by $\ln (-\sigma\ell/e)$ in the following sense.
\[\sigma u-\ln\frac{-\sigma\ell}{e}=O(\ln|\ell|/\ell),\quad \frac{\partial u}{\partial\ell}=\sigma 1/\ell+O(1/\ell^2),\]
\[\left(\frac{\partial }{\partial L},\frac{\partial }{\partial G}\right)\left(u+\sigma\ln e\right)=O(1/|\ell|),\quad \left(\frac{\partial }{\partial L},\frac{\partial }{\partial G}\right)^2\left(u+\sigma\ln e\right)=O(1/|\ell|),\]
The estimates above are uniform as long as $|G|\leq C,$ $1/C\leq L \leq C,$ $\ell>\ell_0$ for some constant $C>1$
and the implied constants in $O(\cdot)$ depend only on $C$ and $\ell_0.$
%when $u\gtrless 0$ and $|\ell|\to \infty$ for some constant $C$.\label{Lm: simplify}\\
\end{Lm}
For derivatives of the Cartesian variables with respect to the Delaunay variables, we have the following relations when we restrict to the isosceles scenario.
\begin{Lm}\label{LmD2G}In the isosceles limit $G_1=G_0=g_1=g_0-\pi/2=0$, where  $\bx_1$ is in hyperbolic motion horizontally, and $\bx_0$ is in  elliptic motion vertically,  we have
\[\begin{aligned}
&\bx_1=-\left(\frac{L_1^2}{M_1k_1}
\left(\cosh u_1-e_1\right),\frac{L_1G_1}{M_1k_1}\sinh u_1\right),\quad \bx_0=\left(\frac{L_0G_0}{M_0k_0}\sin u_0,-\frac{L_0^2}{M_0k_0}
\left(\cos u_0-e_0\right)\right),\\
&\frac{\partial \bx_0}{\partial G_0}=\frac{L_0}{M_0k_0}(\sin u_0,0),\quad
\frac{\partial^2 \bx_0}{\partial G_0^2}=\frac{1}{M_0k_0}(0, -\cos u_0-2),\\
&\frac{\partial \bx_1}{\partial G_1}=-\frac{L_1}{M_1k_1}(0,\sinh u_1),\quad
\frac{\partial^2\bx_1}{\partial G_1^2}=\frac{1}{M_1k_1}( 2+\cosh u_1,0),\\
&\frac{\partial \bp_i}{\partial G_i}\bot \bp_i,\;\frac{\partial \bp_i}{\partial g_i}\bot \bp_i,\;\frac{\partial \bp_i}{\partial L_i}\parallel \bp_i,\;\frac{\partial \bx_i}{\partial G_i}\bot \bx_i,\; \frac{\partial \bx_i}{\partial g_i}\bot \bx_i,\;\frac{\partial \bx_i}{\partial L_i}\bot \bx_i,\;\frac{\partial \bx_i}{\partial L_i}\parallel \bx_i,\;\frac{\partial \bx_i}{\partial \ell_i}\bot \bx_i,\; i=0,1.\\
\end{aligned}
\]
%In particular,
%\begin{equation}\label{deri-GeH}\frac{\partial \bx_0}{\partial G_0}=\frac{-1}{M_0k_0}(\bx_0\cdot \bp_0,0),\;\frac{\partial^2 \bx_0}{\partial G_0^2}=\frac{1}{M_0k_0}(0,1)-\frac{1}{M_0^2k_0^2}\bx_0\cdot (\bp_0\cdot \bp_0).\end{equation}
Moreover, the formula for $\bx_2$ and its derivatives are obtained from that of $\bx_1$ by changing the subscript and rotating by $\pi.$
\end{Lm}
%\begin{Rk}The relations \eqref{deri-GeH} hold true even with $\bx_0$ in hyperbolic motion. \end{Rk}
Lemmas \ref{Lm-A1} and \ref{LmD2G} follow by direct calculation from the formulas in \eqref{DelEll} and \eqref{eq: delaunay4}.

\section{The derivative of the global map in the I4BP}\label{SSDG}

We now begin the proof of  Proposition \ref{PropDL1} following the strategy outlined in Section \ref{SSStrategyLoc}. It  consists of the following steps. We shall divide the orbit into the left piece (from $\cS_+^L$ to $\cS^M$) and the right piece (from  $\cS^M$ to $\cS_-^R$).
\begin{itemize}
\item Step 1: estimating the fundamental solutions to the variational equation in both the left and right pieces.
\item Step 2: the boundary terms $\mathrm{Id}+  JD_XH\otimes \nabla_X \ell_1$ for both the left and right pieces. 
\item Step 3: the coordinate change from left to right.
\end{itemize}
%We first set up the problem.
%We pick a middle section by setting $x_1^L=\frac\chi2$.

%To the left of the middle section, we divide the Hamiltonian equations by the $\dot\ell_1$ equation and derive the corresponding variational equation for $\dt(L_0,\ell_0,L_2,\ell_2)^L$. We will integrate the variational equation from the left section to the middle section. Next, to the right of the middle section, we need to divide all the Hamiltonian equations by the $\dot \ell_2$ equation and derive the variational equation.  So the estimate of the variational equations is the same for the left and the right provided we permute the subscripts $1$ and $2$ and note that $|\ell_1|$ increases and $|\ell_2|$ decreases. The boundary contributions are used to take into account of the fact that different orbits arrive at the same section at different time. The formula $Id+  JD_XH\otimes \nabla_X \ell_1$ was derived in Section 7 of \cite{X}. The transition from the left coordinates to the right is given in Section \ref{STransition} in terms of Cartesian coordinates. Here we need to do the calculation in the Delaunay coordinates.

%The estimate will be done by assuming $r_2\sim \chi,\ r_1\sim |\ell_1|,\ L_1\sim L_0\sim 1,\ L_2\sim \beta^{-1}.$

{\bf Step 1: the variational equation.}

From Lemma \ref{L01-bound}, we learn that \begin{equation}0<\frac{1}{C'}\leqslant |\dot\ell_1|\leqslant C'\end{equation} for some suitable $C'>1$.  So  we  treat $\ell_1$ as the new time variable. Then we  remove the equation for variable $ L_1$ from the system of equations~\eqref{EqHamLimit}  and divide all the equations for the remaining variables $(L_0,\ell_0, L_2,\ell_2)$ by $\dot\ell_1$. Let $X=(L_0,\ell_0,L_2,\ell_2)$, let $D_X$ be the partial derivative $(\frac{\partial }{\partial L_0},\frac{\partial }{\partial \ell_0},\frac{\partial }{\partial L_2},\frac{\partial }{\partial \ell_2} )$, and let $\nabla_X$ be the covariant derivative $\nabla_X=D_X+(D_X L_1) \partial_{L_1}$, which takes into account of the dependence on $X$ through $L_1$.
The Hamiltonian equation now has the form
$\frac{dX}{d\ell_1} =\frac{1}{\dot \ell_1}J D_X H$ and the variational equation has the form
\begin{equation}\label{EqVar}
\begin{aligned}
\frac{d}{d\ell_1} (\dt X)&=\nabla_X\left(\frac{1}{\dot \ell_1}J D_X H\right)\dt X\\
 &= \left(\frac{1}{\dot \ell_1}J \nabla_X(D_X H)-\frac{1}{(\dot \ell_1)^2}J D_X H\otimes \nabla_X\dot \ell_1\right)\dt X \\
 &=\left(\frac{1}{\dot \ell_1}J (D^2_X H+D_{L_1}D_X H\otimes D_X L_1)-\frac{1}{(\dot \ell_1)^2}J D_X H\otimes \nabla_X\dot \ell_1\right)\dt X
\end{aligned}
\end{equation}
From the Hamiltonian equation \eqref{EqHamLimit}, we have  \begin{equation}\label{H1}JD_XH=O\left(\frac{1}{\ell_1^3},1,\frac{1}{(\beta\chi)^2},\beta^3\right).\end{equation}

We next let $ A(a,b,c)= \frac{\partial U}{\partial a}\frac{\partial^2 a}{\partial b\partial c}+\frac{\partial^2 U}{\partial a^2}\frac{\partial a}{\partial b}\frac{\partial a}{\partial c }$, and $B(a,b)= \frac{\partial^2 U_2}{\partial r_0\partial r_2}\frac{\partial r_0}{\partial a}\frac{\partial r_2}{\partial b}$. Then we get
$$JD^2_{X}H\simeq\left[
\begin{array}{cccc}
-A(r_0, L_0,\ell_0)& -A(r_0,\ell_0,\ell_0)
&-B(\ell_0, L_2)&-B(\ell_0, \ell_2)\\
1+A(r_0, L_0,L_0)& A(r_0,L_0,\ell_0)&B(L_0, L_2)&B(L_0, \ell_2)\\
-B( L_0,\ell_2)&-B(\ell_0,\ell_2)&-A(r_2, L_2,\ell_2)& -A(r_2,\ell_2,\ell_2)
\\
B(L_0,L_2)&B(\ell_0,L_2)&\beta^4+A(r_2, L_2,L_2)& A(r_2,L_2,\ell_2)
\end{array}\right].
$$
 Using Lemma \ref{LmPotential},  this gives the estimate
\begin{equation}\label{H2}JD^2_{X}H=O\left[
\begin{array}{cccccc}
\frac{1}{\ell_1^3}& \frac{1}{r_0\ell_1^3}&\frac{\beta}{\chi^2}&\frac{1}{\chi^3\beta^2}\\
1& \frac{1}{\ell_1^3}&\frac{\beta}{\chi^2}&\frac{1}{\chi^3\beta^2}\\

\frac{1}{\chi^3\beta^2}&\frac{1}{\chi^3\beta^2}&\frac{1}{\chi^2\beta}
&\frac{1}{\chi^3\beta^4}\\
\frac{\beta}{\chi^2}&\frac{\beta}{\chi^2}&\beta^4&\frac{1}{\chi^2\beta}
\end{array}\right].
\end{equation}
Note that the $(1,2)$ entry has a singularity when $r_0=0$, corresponding to double collision between the binary $Q_3$-$Q_4$. We will later show that it is integrable with respect to the time~$\ell_1$.

The following lemma gives the estimates of the partial derivatives of $L_1$ with respect to the other variables.
\begin{Lm}\label{LmDL1}
We have
 $\frac{\partial L_1}{\partial \ell_1}=O\left(\frac{1}{\ell_1^4}+\frac{1}{\chi^2}\right),\ D_XL_1=O\left(1,\frac{r_0}{\ell_1^3},\beta^3,\frac{1}{(\beta\chi)^2}\right). $
 \end{Lm}
 We postpone the proof of this lemma to the end of the section.
Since we have 
\[
JD_X(D_{L_1}H)=\frac{\partial r_1}{\partial L_1}\left(\frac{\partial^2 U}{\partial r_0\partial r_1}\frac{\partial r_0}{\partial L_0},\frac{\partial^2 U}{\partial r_0\partial r_1} \frac{\partial r_0}{\partial \ell_0},\frac{\partial^2 U}{\partial r_1\partial r_2} \frac{\partial r_2}{\partial L_2},\frac{\partial^2 U}{\partial r_1\partial r_2}\frac{\partial r_2}{\partial \ell_2} \right),\]
with the estimates in Lemma \ref{LmPotential}, we have
$JD_X(D_{L_1}H)=O \left(\frac{1}{\ell_1^3},\frac{1}{\ell_1^3}, \frac{\beta}{\chi^2},\frac{1}{\beta^2\chi^3}\right).$
Using Lemma \ref{LmDL1}, we have
\begin{equation}\label{H3}
\nabla_X\dot\ell_1=\left(\frac{3M_1k_1^2}{L_1^4}+A(r_1,L_1,L_1)\right)D_XL_1+JD_X(D_{L_1}H)=O\left(1,\frac{1}{\ell_1^3},\beta^3,\frac{1}{(\beta\chi)^2}\right)
\end{equation}

%There are also two boundary terms given by $(id+d\ell_2\otimes \mathcal F)$.

%From the condition $x_2^L=$const, we obtain $2L_2^L\ell_2^LdL_2^L+(L_2^L)^{2}d\ell_2^L=0 $

Putting the estimates \eqref{H1}, \eqref{H2} and  \eqref{H3}  together according to formula \eqref{EqVar}, we get
\begin{equation}\label{coefficients-M1}\nabla_X\left(\frac{1}{\dot \ell_1}J D_X H\right)=O\left[
\begin{array}{cccccc}
\frac{1}{\ell_1^3}& \frac{1}{r_0\ell_1^3}&\frac{\beta}{\chi^2}+\frac{\beta^3}{\ell_1^3}&\frac{1}{\chi^2\beta}\\
1& \frac{1}{\ell_1^3}&\beta^3&\frac{1}{\chi^2\beta^2}\\
\frac{1}{\beta^2\chi^2}&\frac{1}{\beta^2\chi^3}&\frac{1}{\chi^2\beta}
&\frac{1}{\chi^3\beta^4}\\
\beta^3&\frac{\beta^3}{\ell_1^3}+\frac{1}{\beta^2\chi^2}&\beta^4&\frac{1}{\chi^2\beta}
\end{array}\right].
\end{equation}
We next explain  how to estimate the fundamental solution of the variational equation. The only complication appears at the $(2,1)$ entry, i.e. the term $\frac{\partial\dot L_0}{\partial \ell_0}\simeq\frac{1}{r_0\ell_1^3}$. Note that $r_0$ is allowed to be zero. To integrate the variational equation to obtain the estimate in the statement, we need to show that $1/r_0(\ell_1)$ is integrable with respect to $\ell_1$. Indeed, from the formulas in Lemma~\ref{Lm-A1}, we have that $r_0=0$ if $\ell_0=u_0=0$. Near $\ell_0=0$, we have $u_0-\sin u_0=\frac{u_0^3}{6}+O(u_0^5)=\ell_0$. This gives $u_0=(6\ell_0)^{1/3}+o(\ell_0^{1/3})$. So we have
$$r_0=\frac{1}{M_0 k_0} L_0^2\left(\frac{u_0^2}{2}+o(u_0^2)\right)=\frac{L_0^2}{M_0 k_0}\frac{(6\ell_0)^{2/3}}{2}+o(\ell_0^{2/3}).$$
This shows that $\int_0^1\frac{1}{r_0(\ell_0)}\,d\ell_0$ is bounded. Taking the change of variable $\ell_0=\ell_0(\ell_1)$, since for orbits between the sections $S^+$ and $S^M$, we have  $\frac{1}{C''}\leqslant| \frac{d\ell_1}{d\ell_0}|\leqslant C''$ for suitable $C''\geqslant1$, we know that $\int_{\ell_1(0)}^{\ell_1(1)}\frac{1}{r_0(\ell_1)}\,d\ell_1$ is also bounded.  Therefore, when integrating with respect to the time $\ell_1$, the quantity $\frac{1}{r_0\ell_1^3}$ is equivalent to $\frac{1}{\ell_1^3}$.

Finally, we integrate the variational equation $\frac{d}{d\ell_1} (\dt X)=\nabla_X\left(\frac{1}{\dot \ell_1}J D_X H\right)\dt X$ to obtain an estimate of its fundamental solution. We integrate over a ``time" (that is, $\ell_1$) interval of length of order $\chi$, since on the sections $\cS^L_+$ and $\cS^M$ we have $\ell_1\sim 1/\eps$ and $\ell_1\sim\chi$ respectively. 
Consider a linear ODE given by $\frac{dY}{dt}=\Lambda(t)Y$ with the initial condition $Y(0)$ as an identity matrix. Using Picard iteration, we write the solution  as
\begin{equation}\label{EqPicard}Y(t)=\text{Id}+\int_0^t\Lambda(s)\cdot Y(s)ds=\text{Id}+\int_0^t\Lambda(s)ds+\int_0^t\Lambda(s)(\int_0^s\Lambda(s')ds')ds+\cdots.
 \end{equation}
It turns out that the Picard iteration stablizes quickly and two steps are enough.
We have the following estimate for the fundamental solution of the corresponding variational equation,
\begin{equation}\label{D-X-L}
\cD^L_2= \Id+O\left[
\begin{array}{cccc}
 \epsilon  & \epsilon ^2 & \beta ^3 & \frac{1}{\beta  \chi } \\
 \chi  & \epsilon ^2 \chi  & \beta ^3 \chi  & \frac{1}{\beta } \\
 \frac{1}{\beta ^2 \chi } & \frac{\epsilon ^2}{\beta ^2 \chi } & \frac{1}{\beta  \chi } & \frac{1}{\beta ^3 \chi ^2} \\
 \beta ^3 \chi  & \beta ^3 \epsilon ^2 \chi  & \beta ^4 \chi  & \beta ^2 \\
\end{array}
\right],
\end{equation}
For the orbits between the middle section $\cS^M$ and the right section $\cS_-^R$, the estimates in \eqref{coefficients-M1} of the coefficient matrix for the corresponding variational equation remain valid if we replace $\beta$ with $1$ and change the subindex $1$ to $2$. Then integrating  the variational equation over the time $\ell_2$ on an interval of order $O(\beta\chi)$, decreasing from $\beta\chi$ on the section $\cS^M$ to $O(\epsilon^{-1})$ on the section $\cS_-^R$,  we have the following estimates for its fundamental solution,
\begin{equation}\label{D-X-R}
\cD^R_3=\Id+O\left[
\begin{array}{cccc}
 \beta  \epsilon ^2 \chi  & \epsilon ^2 & \beta  \epsilon ^2 \chi  & \frac{\beta }{\chi } \\
 \beta  \chi  & \epsilon  & \beta  \chi  & \beta ^2 \\
 \frac{\beta }{\chi } & \frac{\epsilon }{\chi ^2} & \frac{\beta }{\chi } & \frac{\beta }{\chi ^2} \\
 \beta  \chi  & \epsilon  & \beta  \chi  & \beta ^2 \\
\end{array}
\right].
\end{equation}
Here the order $\beta\epsilon^2\chi$ terms in the first row come from quantities of the form $\frac{1}{\ell_2^3}(\beta\chi-\ell_2)$ which appear  in the second iteration when we multiply the (1,2) entry and the (2,1) entry. 

{\bf Step 2, the boundary contributions.}

The boundary terms $(\Id+\nabla_X\ell_{1,2}\otimes JD_XH)$ are calculated on the sections $\cS_+^L$, $\cS_-^R$ and $\cS^M$.  For the middle section $\cS^M$, we need to compute the boundary contribution before and after  changing  Jacobi coordinates from left  to  right.
We first consider the left side of the  section $\cS^M=\{-\beta \bx_1^L=\bx_2^R\}$. Recall that $X=(L_0,\ell_0,L_2^L,\ell_2^L)$. For the Hamiltonian equation, we have $JD_XH=O\left(\frac{1}{\chi^3},1,\frac{1}{(\beta\chi)^2} ,\beta^3\right).$ Using \eqref{EqLR} we get $-\beta \bx_1^L=\bx_2^R=\frac12M_1^L\bx_1^L+ \bx_2^L$ on $\cS^M$. Substituting  the formulas from Lemma \ref{LmD2G} and differentiating the above relation, we have
$$\frac{1}{M^L_1k^L_1}\left(\frac12M_1^L+\beta\right)(2L_1^L\ell_1^LdL_1^L+(L_1^L)^2d\ell_1^{L})=-\frac{1}{M^L_2k^L_2}(2L^L_2\ell^L_2dL^L_2+(L^L_2)^2d\ell^L_2).$$
From the implicit function theorem, we get $D_X\ell_1^L=O(0,0,\beta\chi,\frac{1}{\beta^2}),\; \partial_{L_1^L}\ell_1^L\simeq\chi.$
Combining with Lemma \ref{LmDL1}, we have  $\nabla_X\ell_1^L=O\left(\chi,\frac{1}{\chi^2}, \beta\chi,\frac{1}{\beta^2}\right).$
So, from the left side of the section $\cS^M$, we have the boundary contribution,
\begin{equation}\cD^L_M=\Id+O\left(\frac{1}{\chi^3},1,\frac{1}{(\beta\chi)^2} ,\beta^3\right)\otimes O\left(\chi,\frac{1}{\chi^2}, \beta\chi,\frac{1}{\beta^2}\right) \end{equation}

Next, working from the right side of the  section $\cS^M$, we have $X=(L_0,\ell_0,L^R_1,\ell^R_1)$. For the Hamiltonian vector field, we have
$JD_XH=O\left(\frac{1}{(\beta\chi)^3},1,\frac{1}{\chi^2} ,1\right).$
From  \eqref{EqLR}, we obtain $-\beta (\bx_1^R+\frac12M_2^R\bx_2^R)=\bx_2^R$  on $\cS^M$. Once again we substitute the formulas from  Lemma \ref{LmD2G} and differentiate the above relation. We obtain,
$$\frac{1}{M^R_1k^R_1}\beta(2L_1^R\ell_1^RdL_1^R+(L_1^R)^{2}d\ell_1^R)=\frac{1}{M^R_2k^R_2}\left(1+\frac\beta2M_2^R\right)\Big(2L_2^R\ell_2^RdL_2^R+(L_2^R)^{2}d\ell_2^R\Big).$$
Using the  implicit function theorem, we have
$D_X\ell_2^R=O(0,0,\beta\chi,\beta),\; \partial_{L_2^R}\ell_2^R\sim \beta\chi,$ and with Lemma \ref{LmDL1}\footnote{We need to replace $\beta$ with 1, $L_1^L$ with $L_2^R$ and $\ell_1$ with $\ell_2^R$}, we obtain
$\nabla_X\ell_2^R=O\left(\beta\chi,\frac{1}{\beta^2\chi^2},\beta \chi,\beta\right).$
Therefore, from the right side of the section $\cS^M$, we have the boundary contribution,
\[\cD^R_M=\Id+O\left(\frac{1}{(\beta\chi)^3},1,\frac{1}{\chi^2} ,1\right)\otimes O\left(\beta\chi,\frac{1}{(\beta\chi)^2},\beta \chi,\beta\right),\]

On the section $\cS_+^L$, after renormalization,  we have the following by Proposition \ref{PropRenorm},
 \[\frac{1}{C}\leqslant L_1,L_0,\ell_0,\beta L_2\leqslant C,\; r_0<C, \text{ and }\ell_1\sim\frac{1}{\epsilon},\; r_2\sim \chi.\] For the Hamiltonian vector field we have $JD_X H= O(\eps^3,1,\frac{1}{\beta^2\chi^2},\beta^3)$. Next, using the formulas from  Lemma \ref{Lm-A1}, we have that
\[
\begin{aligned}
0&=dv=-\frac{1}{2}\frac{dr}{r} v+r^{-1/2}d(L_0\sin u_0+L_1\sinh u_1)\\
&=-\frac{1}{2}\left(\frac{\sqrt M_0 r_0}{r} dr_0+\frac{\sqrt M_1 r_1}{r} dr_1\right)\frac{v}{r}+r^{-1/2}d(L_0\sin u_0+L_1\sinh u_1)\\
&\simeq \eps^{3/2} (r_0\frac{\partial r_0}{\partial L_0}dL_0+r_0\frac{\partial r_0}{\partial \ell_0}d\ell_0)+\epsilon^{1/2}(\frac{\partial r_1}{\partial L_1}dL_1+\frac{\partial r_1}{\partial \ell_1}d\ell_1)\\
&+\eps^{1/2}(\sin u_0 dL_0+L_0\frac{\cos u_0}{1-\cos u_0}d\ell_0+\sinh u_1 dL_1+L_1\frac{\cosh u_1}{1-\cosh u_1}d\ell_1)\\
&=(\eps^{1/2}\sin u_0+O(\eps^{3/2}) dL_0+(\eps^{1/2}L_0\frac{\cos u_0}{1-\cos u_0}+O(\eps^{3/2})) d\ell_0+\eps^{-1/2} dL_1+O(\eps^{1/2})d\ell_1
\end{aligned}
\]
Note that $r_0=\frac{1}{M_0k_0}L_0^2(1-\cos u_0)$, we then have $L_0\frac{\cos u_0}{1-\cos u_0}\sim r_0^{-1}\sim1.$
Therefore, $\partial_{\ell_0}\ell_1\sim1$.
This gives
$D_X\ell_1=O(1,1,0,0),\ D_{L_1}\ell_1\sim\epsilon^{-1},\ \text{and}\ \nabla_X\ell_1=O\left(\eps^{-1},1,\eps^{-1}\beta^3,\frac{\eps^{-1}}{(\beta\chi)^2}\right).$
So we get the boundary contribution from the section $\cS_+^L$,
\[\cD^L_+=\Id+O(\eps^3,1,\frac{1}{\beta^2\chi^2},\beta^3)\otimes O\left(\eps^{-1},1,\eps^{-1}\beta^3,\frac{\eps^{-1}}{(\beta\chi)^2}\right).\]
On the section $\cS_-^R=\{r=\epsilon^{-1}\}$,
\[0=dr=\frac{1}{r}(M_0r_0dr_0+M_2r_2dr_2)=\frac{M_0r_0}{r}(\frac{dr_0}{d L_0}d L_0+\frac{dr_0}{d\ell_0}d\ell_0)+\frac{M_2r_2}{r}(\frac{d r_2}{d L_2}dL_2+\frac{dr_2}{d\ell_2}d\ell_2).\]
% \[\frac{1}{C}\leqslant L_1,L_0, L_2\leqslant C,\; r_0=O(1), \text{ and }\ell_2\sim\frac{1}{\epsilon},\; r_1\sim \chi.\]
Then using again Lemma \ref{LmDL1} and \ref{Lm-A1}, we have
\[D_X\ell_2=O(\epsilon, \epsilon, 0,0),\quad D_{L_2}\ell_2\sim\epsilon^{-1}\quad \text{and}\quad \nabla_X\ell_2=O(\epsilon^{-1},\epsilon^2,\epsilon^{-1},\frac{\epsilon^{-1}}{\chi^2}).\]
For the Hamiltonian vector field, we have
$J\nabla_X H=O(\epsilon^3, 1,\frac{1}{\chi^2},1).$
Hence, we get the boundary contribution from the section $S^R_-$,
\[
\cD^R_-=\Id+O(\eps^3,1,\frac{1}{\chi^2},1)\otimes O\left(\eps^{-1},\eps^{2},\eps^{-1},\frac{\eps^{-1}}{\chi^2}\right).
\]

{\bf Step 3, transition from the left to the right. }

At the middle section $\cS^M$, we change the left Jacobi coordinates to the right Jacobi coordinates. By Lemma \ref{LmdL}, we have that on $\cS^M$,  with $X=(L_0,\ell_0,L_2^L,\ell_2^L)$,
\[D_XL_1^{R}=O(0,0,\beta,\frac{1}{\beta^2\chi}),\quad D_{L_1^L}L_1^R\simeq1,\;D_{\ell_1^L}L_1^R\simeq\frac{1}{\chi}, \]\[ D_X\ell_1^R=O(0,0,\beta\chi,\beta^{-2}),\;\partial_{L_1^L}\ell_1^R\simeq\chi, \; D_{\ell_1^L}\ell_1^R\simeq1. \]
With Lemma \ref{LmDL1}, we have
$\nabla_XL_1^R=O(1,\frac{1}{\chi^3},\beta,\frac{1}{\beta^2\chi})$, and $\nabla_X\ell_1^R=O(\chi,\frac{1}{\chi^2},\beta\chi,\beta^{-2}).$
Note that the transition from left to right does not effect the binary $\bx_0$. We then  obtain the  following transition matrix on the section $\cS^M$, 
\begin{equation}\label{middle-d-l}\cM=\frac{\partial(L_0,\ell_0,L_1,\ell_1)^R|_{\cS^M}}{\partial(L_0,\ell_0,L_2,\ell_2)^L|_{\cS^M}}=O\left[\begin{array}{cccc}
1&0&0&0\\
0&1&0&0\\
1&\frac{1}{\chi^3}&\beta&\frac{1}{\beta^2\chi}\\
\chi&\frac{1}{\chi^2}&\beta\chi&\beta^{-2}
\end{array}\right].
\end{equation}

%Note that in the integration, we use the fact $\int_0^1\frac{1}{r_0(\ell_0)}\,d\ell_0$ is bounded to handle the singularity $1/r_0$.

After completing the three steps, using the formula \eqref{eq: formald4}, we multiply all the matrices together and obtain the derivative matrix $$\frac{\partial(L^R_0,\ell^R_0,L^R_1,\ell^R_1)\big|_{S_-^R}}{\partial(L^L_0,\ell^L_0,L^L_2,\ell^L_2)\big|_{S_+^L}}=(\cD_-^R)^{-1}\cD_3^R\cD_M^R\cM (\cD_M^L)^{-1}\cD_2^L\cD_+^L=O\left[
\begin{array}{cccc}
 \epsilon ^2 \chi  & \epsilon ^4 \chi  & \beta  \epsilon ^2 \chi  & \frac{\epsilon ^2}{\beta ^2} \\
 \chi  & \epsilon ^2 \chi  & \beta  \chi  & \frac{1}{\beta ^2} \\
 1 & \epsilon ^2 & \beta  & \frac{1}{\beta ^2 \chi } \\
 \chi  & \epsilon ^2 \chi  & \beta  \chi  & \frac{1}{\beta ^2} \\
\end{array}
\right]. $$
Thus we complete the proof of Proposition \ref{PropDG1}. \qed

We end this section with the proof of Lemma \ref{LmDL1}.

\begin{proof}[Proof of Lemma \ref{LmDL1}]
Differentiating the Hamiltonian we get
\[
\begin{aligned}
0&=\frac{M_0k_0^2}{L_0^3}dL_0-\frac{M_1k_1^2}{L_1^3}dL_1 -\frac{M_2k_2^2}{L_2^3}dL_2 +dU_{01}+dU_2\\
&=\frac{2M_0k_0^2}{L_0^3}dL_0-\frac{2M_1k_1^2}{L_1^3}dL_1 -\frac{2M_2k_2^2}{L_2^3}dL_2\\
&+\left(\frac{\partial U_{01}}{\partial r_1} +\frac{\partial U_2}{\partial r_1}\right) \left(\frac{\partial r_1}{\partial L_1}d L_1+ \frac{\partial r_1}{\partial \ell_1}d \ell_1\right)+\left(\frac{\partial U_{01}}{\partial r_0} +\frac{\partial U_2}{\partial r_0}\right) \left(\frac{\partial r_0}{\partial L_0}d L_0+ \frac{\partial r_0}{\partial \ell_0}d \ell_0\right)\\
&+\frac{\partial U_2}{\partial r_2} \left(\frac{\partial r_2}{\partial L_2}d L_2+ \frac{\partial r_2}{\partial \ell_2}d \ell_2\right)
\end{aligned}
\]
The estimate in the statement is then obtained by the implicit function theorem.
\end{proof}
\section{The second order derivatives with respect to $G$ and $g$}\label{App-2-derviatives}
In this section, we compute the second derivatives of the potentials with respect to $G_i,g_i$, $i=0,1,2$ and evaluate them in the isosceles limit.
%Recall the potentials
%\[
%\begin{aligned}
%&k_0=1,\ k_1=2m_1,\ k_2=(m_1+2)m_2\\
%&U_{01}=\left(\frac{k_1}{|\bx_1|}-\frac{m_1}{|\bx_1-\frac{\bx_0}{2}|}-\frac{m_1}{|\bx_1+\frac{\bx_0}{2}|}\right),\\
%& U_2=\left(\frac{k_2}{|\bx_2|}-\frac{m_1m_2}{|\bx_2-\frac{2 \bx_1}{m_1+2} |}-\frac{m_2}{|\bx_2+\frac{m_1\bx_1}{m_1+2}-\frac{\bx_0}{2}|}-\frac{m_2}{|\bx_2+\frac{m_1\bx_1}{m_1+2}+\frac{\bx_0}{2}|}\right).
%\end{aligned}
%\]

\subsection{The second order derivatives of $U_{01}$ with respect to $G_i,g_i$, i=0,1}
Note that  $g_i$ and $G_i$ are contained only in $x_i$, $i=0,1$. So we have the following calculation.
\[\begin{aligned}
\frac{1}{m_1}\frac{\partial U_{01}}{\partial g_0}&=-\frac{\langle \bx_1,\partial_{g_0} \frac{\bx_0}{2}\rangle}{|\bx_1-\frac{\bx_0}{2}|^3}+\frac{\langle \bx_1,\partial_{g_0}\frac{\bx_0}{2}\rangle}{|\bx_1+\frac{\bx_0}{2}|^3}\\
\frac{1}{m_1}\frac{\partial U_{01}}{\partial G_0}&=-\frac{\langle \bx_1-\frac{\bx_0}{2},\partial_{G_0} \frac{\bx_0}{2}\rangle}{|\bx_1-\frac{\bx_0}{2}|^3}+\frac{\langle \bx_1+\frac{\bx_0}{2},\partial_{G_0}\frac{\bx_0}{2}\rangle}{|\bx_1+\frac{\bx_0}{2}|^3}\\
\frac{1}{m_1}\frac{\partial^2 U_{01}}{\partial g_0^2}&=-\frac{\langle \bx_1,\partial_{g_0}^2\frac{\bx_0}{2}\rangle}{|\bx_1-\frac{\bx_0}{2}|^3}+\frac{\langle \bx_1,\partial_{g_0}^2\frac{\bx_0}{2}\rangle}{|\bx_1+\frac{\bx_0}{2}|^3}-\frac{3\langle \bx_1,\partial_{g_0}\frac{\bx_0}{2}\rangle^2}{ |\bx_1-\frac{\bx_0}{2}|^5}-\frac{3\langle \bx_1,\partial_{g_0}\frac{\bx_0}{2}\rangle^2}{|\bx_1+\frac{\bx_0}{2}|^5}\\
\frac{1}{m_1}\frac{\partial^2 U_{01}}{\partial G_0^2}&=-\frac{\langle-\partial_{G_0}\frac{\bx_0}{2},\partial_{G_0}\frac{\bx_0}{2}\rangle+\langle \bx_1-\frac{\bx_0}{2},\partial_{G_0}^2\frac{\bx_0}{2}\rangle}{|\bx_1-\frac{\bx_0}{2}|^3}+\frac{\langle\partial_{G_0}\frac{\bx_0}{2},\partial_{G_0}\frac{\bx_0}{2}\rangle+\langle \bx_1+\frac{\bx_0}{2},\partial_{G_0}^2\frac{\bx_0}{2}\rangle}{|\bx_1+\frac{\bx_0}{2}|^3}\\
&-\frac{3\langle \bx_1-\frac{\bx_0}{2},\partial_{G_0}\frac{\bx_0}{2}\rangle^2}{ |\bx_1-\frac{\bx_0}{2}|^5}-\frac{3\langle \bx_1+\frac{\bx_0}{2},\partial_{G_0}\frac{\bx_0}{2}\rangle^2}{|\bx_1+\frac{\bx_0}{2}|^5}\\
\frac{1}{m_1}\frac{\partial^2 U_{01}}{\partial g_0\partial G_0}&=-\frac{\langle \bx_1,\partial_{g_0 G_0}^2\frac{\bx_0}{2}\rangle}{|\bx_1-\frac{\bx_0}{2}|^3}+\frac{\langle \bx_1,\partial_{g_0G_0}^2\frac{\bx_0}{2}\rangle}{|\bx_1+\frac{\bx_0}{2}|^3}
\end{aligned}\]
\[\begin{aligned}
&-\frac{3\langle \bx_1,\partial_{g_0}\frac{\bx_0}{2}\rangle\langle \bx_1-\frac{\bx_0}{2},\partial_{G_0}\frac{\bx_0}{2}\rangle}{ |\bx_1-\frac{\bx_0}{2}|^5}-\frac{3\langle \bx_1,\partial_{g_0}\frac{\bx_0}{2}\rangle\langle \bx_1+\frac{\bx_0}{2},\partial_{G_0}\frac{\bx_0}{2}\rangle}{|\bx_1+\frac{\bx_0}{2}|^5}\\
\frac{1}{m_1}\frac{\partial U_{01}}{\partial g_1}&=\frac{\langle -\frac{\bx_0}{2},\partial_{g_1}\bx_1\rangle}{|\bx_1-\frac{\bx_0}{2}|^3}+\frac{\langle \frac{\bx_0}{2},\partial_{g_1}\bx_1\rangle}{|\bx_1+\frac{\bx_0}{2}|^3}\\
\frac{1}{m_1}\frac{\partial U_{01}}{\partial G_1}&=-\frac{2\langle\bx_1,\partial_{G_1}\bx_1\rangle}{|\bx_1|^3}+\frac{\langle \bx_1-\frac{x_0}{2},\partial_{G_1}\bx_1\rangle}{|\bx_1-\frac{\bx_0}{2}|^3}+\frac{\langle \bx_1+\frac{\bx_0}{2},\partial_{G_1}\bx_1\rangle}{|\bx_1+\frac{\bx_0}{2}|^3}\\
\frac{1}{m_1}\frac{\partial^2 U_{01}}{\partial g_1^2}&=\frac{\langle -\frac{\bx_0}{2},\partial_{g_1}^2\bx_1\rangle}{|\bx_1-\frac{\bx_0}{2}|^3}+\frac{\langle \frac{\bx_0}{2},\partial_{g_1}^2\bx_1\rangle}{|\bx_1+\frac{\bx_0}{2}|^3}
-3\frac{\langle -\frac{\bx_0}{2},\partial_{g_1}\bx_1\rangle^2}{|\bx_1-\frac{\bx_0}{2}|^5}-3\frac{\langle \frac{\bx_0}{2},\partial_{g_1}\bx_1\rangle^2}{|\bx_1+\frac{\bx_0}{2}|^{5}}\\
\frac{1}{m_1}\frac{\partial^2 U_{01}}{\partial G_1^2}&=-2\frac{\langle \partial_{G_1}\bx_1,\partial_{G_1}\bx_1\rangle+\langle\bx_1,\partial_{G_1}^2\bx_1\rangle}{|\bx_1|^3}+6\frac{\langle\bx_1,\partial_{G_1}\bx_1\rangle^2}{|\bx_1|^5}+\frac{\langle \partial_{G_1}\bx_1, \partial_{G_1}\bx_1\rangle+\langle \bx_1-\frac{\bx_0}{2},\partial_{G_1}^2\bx_1\rangle}{|\bx_1-\frac{x_0}{2}|^3}\\
&-3\frac{\langle \bx_1-\frac{\bx_0}{2},\partial_{G_1}\bx_1\rangle^2}{|\bx_1-\frac{\bx_0}{2}|^5}+\frac{\langle \partial_{G_1}\bx_1, \partial_{G_1}\bx_1\rangle+\langle \bx_1+\frac{\bx_0}{2},\partial_{G_1}^2\bx_1\rangle}{|\bx_1+\frac{\bx_0}{2}|^3}-3\frac{\langle \bx_1+\frac{\bx_0}{2},\partial_{G_1}\bx_1\rangle^2}{|\bx_1+\frac{\bx_0}{2}|^5}\\
\frac{1}{m_1}\frac{\partial^2 U_{01}}{\partial g_1\partial G_1}&=\frac{\langle-\frac{\bx_0}{2},\partial_{g_1G_1}\bx_1\rangle}{|\bx_1-\frac{\bx_0}{2}|^3}+\frac{\langle\frac{\bx_0}{2},\partial_{g_1G_1}\bx_1\rangle}{|\bx_1+\frac{\bx_0}{2}|^3}-3\frac{\langle-\frac{\bx_0}{2},\partial_{g_1}\bx_1\rangle\langle\bx_1-\frac{\bx_0}{2},\partial_{G_1}\bx_1\rangle}{|\bx_1-\frac{\bx_0}{2}|^5}\\
&-3\frac{\langle\frac{\bx_0}{2},\partial_{g_1}\bx_1\rangle\langle\bx_1+\frac{\bx_0}{2},\partial_{G_1}\bx_1\rangle}{|\bx_1+\frac{\bx_0}{2}|^5}
\end{aligned}
\]
\[\begin{aligned}
\frac{1}{m_1}\frac{\partial^2U_{01}}{\partial g_0\partial g_1}&=\frac{\langle-\partial_{g_0}\frac{\bx_0}{2},\partial_{g_1}\bx_1\rangle}{|\bx_1-\frac{\bx_0}{2}|^3}-3\frac{\langle -\frac{\bx_0}{2},\partial_{g_1}\bx_1\rangle\langle \bx_1,-\partial_{g_0}\frac{\bx_0}{2}\rangle}{|\bx_1-\frac{x_0}{2}|^5}+\frac{\langle\partial_{g_0}\frac{\bx_0}{2},\partial_{g_1}\bx_1\rangle}{|\bx_1+\frac{\bx_0}{2}|^3}-3\frac{\langle \frac{\bx_0}{2},\partial_{g_1}\bx_1\rangle\langle \bx_1,\partial_{g_0}\frac{\bx_0}{2}\rangle}{|\bx_1+\frac{x_0}{2}|^5}\\
\frac{1}{m_1}\frac{\partial^2U_{01}}{\partial g_0\partial G_1}&=\frac{\langle-\partial_{g_0}\frac{\bx_0}{2},\partial_{G_1}\bx_1\rangle}{|\bx_1-\frac{\bx_0}{2}|^3}-3\frac{\langle -\frac{\bx_0}{2},\partial_{G_1}\bx_1\rangle\langle \bx_1,-\partial_{g_0}\frac{\bx_0}{2}\rangle}{|\bx_1-\frac{x_0}{2}|^5}+\frac{\langle\partial_{g_0}\frac{\bx_0}{2},\partial_{g_1}\bx_1\rangle}{|\bx_1+\frac{\bx_0}{2}|^3}-3\frac{\langle \frac{\bx_0}{2},\partial_{G_1}\bx_1\rangle\langle \bx_1,\partial_{g_0}\frac{\bx_0}{2}\rangle}{|\bx_1+\frac{x_0}{2}|^5}\\
\frac{1}{m_1}\frac{\partial^2U_{01}}{\partial G_0g_1}&=\frac{\langle-\partial_{G_0}\frac{\bx_0}{2},\partial_{g_1}\bx_1\rangle}{|\bx_1-\frac{\bx_0}{2}|^3}-3\frac{\langle -\frac{\bx_0}{2},\partial_{g_1}\bx_1\rangle\langle \bx_1,-\partial_{G_0}\frac{\bx_0}{2}\rangle}{|\bx_1-\frac{x_0}{2}|^5}+\frac{\langle\partial_{g_0}\frac{\bx_0}{2},\partial_{g_1}\bx_1\rangle}{|\bx_1+\frac{\bx_0}{2}|^3}-3\frac{\langle \frac{\bx_0}{2},\partial_{g_1}\bx_1\rangle\langle \bx_1,\partial_{G_0}\frac{\bx_0}{2}\rangle}{|\bx_1+\frac{\bx_0}{2}|^5}\\
\frac{1}{m_1}\frac{\partial^2 U_{01}}{\partial G_0\partial G_1}&=\frac{\langle-\partial_{G_0}\frac{\bx_0}{2}, \partial_{G_1}\bx_1\rangle}{|\bx_1-\frac{\bx_0}{2}|^3}-3\frac{\langle\bx_1-\frac{\bx_0}{2},\partial_{G_1}\bx_1\rangle\langle \bx_1-\frac{\bx_0}{2},-\partial_{G_0}\frac{\bx_0}{2}\rangle}{|\bx_1-\frac{\bx_0}{2}|^5}\\
&+\frac{\langle\partial_{G_0}\frac{\bx_0}{2}, \partial_{G_1}\bx_1\rangle}{|\bx_1+\frac{\bx_0}{2}|^3}-3\frac{\langle\bx_1+\frac{\bx_0}{2},\partial_{G_1}\bx_1\rangle\langle \bx_1+\frac{\bx_0}{2},\partial_{G_0}\frac{\bx_0}{2}\rangle}{|\bx_1+\frac{\bx_0}{2}|^5}
\end{aligned}
\]
By Lemma \ref{LmD2G},  when evaluating at the isosceles limit, the above second order derivatives can be  simplified to the following
\begin{equation}\label{w0g0-deri-1}
\begin{aligned}
\frac{1}{m_1}\frac{\partial^2 U_{01}}{\partial G_0\partial g_0}&=-\frac{3\langle \bx_1,\partial_{g_0}\bx_0\rangle\langle\bx_1,\partial_{G_0}\bx_0\rangle}{2(r_1^2+r_0^2/4)^{5/2}},\quad \frac{1}{m_1}\frac{\partial^2 U_{01}}{\partial g_0^2}=\frac{-3(r_1r_0)^2}{2(r_1^2+r_0^2/4)^{5/2}},\\
\frac{1}{m_1}\frac{\partial^2 U_{01}}{\partial G_0^2}&=\frac{ \langle \bx_0, \frac{\partial^2 \bx_0}{\partial G_0^2}\rangle+\langle \frac{\partial \bx_0}{\partial G_0}, \frac{\partial \bx_0}{\partial G_0}\rangle}{2(r_1^2+r_0^2/4)^{3/2}}-\frac{3\langle \bx_1, \frac{\partial \bx_0}{\partial G_0}\rangle^2}{2(r_1^2+r_0^2/4)^{5/2}},\\
\frac{1}{m_1}\frac{\partial^2U_{01}}{\partial g_0\partial G_1}&=-\frac{3}{2}\frac{\langle \bx_0,\partial_{G_1}\bx_1\rangle\langle \bx_1,\partial_{g_0}\bx_0\rangle}{|\bx_1-\frac{\bx_0}{2}|^5},\quad \frac{1}{m_1}\frac{\partial^2U}{\partial G_0\partial G_1}=-\frac{3}{2}\frac{\langle \bx_0,\partial_{G_1}\bx_1\rangle\langle \bx_1,\partial_{G_0}\bx_0\rangle}{|\bx_1-\frac{\bx_0}{2}|^5},\\
\frac{1}{m_1}\frac{\partial^2U_{01}}{\partial G_1^2}&=-2(\langle\partial_{G_1}\bx_1,\partial_{G_1}\bx_1\rangle+\langle \bx_1,\partial_{G_1}^2\bx_1\rangle)(\frac{1}{|\bx_1|^3}-\frac{1}{|\bx_1-\frac{\bx_0}{2}|^3})-\frac{3}{2}\frac{\langle \bx_0,\partial_{G_1}\bx_1\rangle^2}{|\bx_1-\frac{\bx_0}{2}|^5}.
\end{aligned}
\end{equation}
Note $\bp_0=(0, R_0), \bx_0=(0,x_0)=M_0^{-1/2}r(0,\cos\psi)$, $\bx_1=(x_1,0)=M_1^{-1/2}r(\sin\psi,0)$, 
$R_0=\sqrt{M_0}r^{-1/2}(v\cos\psi-w\sin\psi)$ and $R_1=\sqrt{M_1}r^{-1/2}(v\sin\psi+w\cos\psi)$. Then  by Lemma \ref{LmD2G}, we have $\partial_{g_0}\bx_0=-M_0^{-1/2}r(\cos\psi,0)$,
%\[x_0\cdot p_0=r^{1/2}(v\cos^2\psi-w\sin\psi\cos\psi).\]
\[\begin{aligned}\frac{\partial \bx_0}{\partial G_0}&= {sign}(x_0)\frac{1}{M_0k_0}\big(\langle \bx_0,\bp_0\rangle,0\big)= {sign}(x_0)\frac{1}{M_0k_0}(r^{1/2}(v\cos^2\psi-\frac{1}{2}w\sin2\psi),0).\\
\frac{\partial^2\bx_0}{\partial G_0^2}&=\frac{1}{M_0k_0}(0,- {sign}(x_0))-\frac{|\bp_0|^2}{M_0^2k_0^2}\bx_0=\frac{1}{M_0k_0}(0,- {sign}(x_0)-\frac{1}{M_0k_0}M_0^{1/2}\cos\psi(v\cos\psi-w\sin\psi)^2)
\end{aligned}
\]
Similarly, we have the  formulas  for $\bx_1$,
\[\begin{aligned}&\frac{\partial\bx_1}{\partial G_1}= {sign}(x_1)\frac{1}{M_1k_1}(0,\langle \bx_1,\bp_1\rangle)= {sign}(x_1)\frac{1}{M_1k_1}(0,r^{1/2}(v\sin^2\psi+\frac{1}{2}w\sin 2\psi)),\\
&\frac{\partial^2\bx_1}{\partial G_1^2}=\frac{1}{M_1k_1}(- {sign}(x_1),0)-\frac{|\bp_1|^2}{M_1^2k_1^2}\bx_1=\frac{1}{M_1k_1}(- {sign}(x_1)-\frac{1}{M_1k_1}M_1^{1/2}\sin\psi(v\sin\psi+w\cos\psi)^2,0).\end{aligned}\]
 Hence we have the following statement.
\begin{Lm}\label{w0g0-deri} When restricted to the isosceles limit, in the blowup coordinates, we have
\[\begin{aligned}\frac{1}{m_1}\frac{\partial^2 U_{01}}{\partial g_0^2}&=-\frac{1}{r}\frac{3(M_0M_1)^{-1}(\sin\psi \cos\psi)^2}{2(\frac{1}{M_1}\sin^2\psi+\frac{\cos^2\psi}{4M_0})^{5/2}},\\
\frac{1}{m_1}\frac{\partial^2 U_{01}}{\partial G_0^2}&=\frac{-\frac{\cos\psi}{M_0^{3/2}k_0}}{2r^2(\frac{1}{M_1}\sin^2\psi+\frac{\cos^2\psi}{4M_0})^{3/2}}
 -\frac{3}{2}\frac{\big(\frac{1}{M_1^{1/2}M_0k_0}\sin\psi(v\cos^2\psi-\frac{1}{2}w\sin2\psi)\big)^2}{r^2(\frac{1}{M_1}\sin^2\psi+\frac{\cos^2\psi}{4M_0})^{5/2}},\\
\frac{1}{m_1}\frac{\partial^2 U_{01}}{\partial G_0\partial g_0}&=\frac{3}{2}\frac{\frac{1}{(M_0)^{3/2}M_1k_0}\sin^2\psi\cos\psi(v\cos^2\psi-\frac{1}{2}w\sin2\psi)}{r^{3/2}(M_1^{-1}\sin^2\psi+\frac{\cos^2\psi}{4M_0})^{5/2}},\\
\frac{1}{m_1}\frac{\partial^2U_{01}}{\partial g_0\partial G_1}&
=\frac{3}{2}\frac{\frac{1}{M_1^{3/2}M_0k_1}(\frac{1}{2}\sin2\psi)^2(v|\sin\psi|+sign(\psi)w\cos\psi)}{r^{3/2}(M_1^{-1}\sin^2\psi+\frac{\cos^2\psi}{4M_0})^{5/2}},\\
\frac{1}{m_1}\frac{\partial^2 U_{01}}{\partial G_0\partial G_1}&
=-\frac{3}{2}\frac{\frac{1}{M_0^{3/2}M_1^{3/2}k_0k_1}\sin^2\psi\cos^2\psi(v|\sin\psi|+sign(\psi)w\cos\psi)(v\cos\psi-w\sin\psi)}{r^2(M_1^{-1}\sin^2\psi+\frac{\cos^2\psi}{4M_0})^{5/2}},\\
\frac{1}{m_1}\frac{\partial^2 U_{01}}{\partial G_1^2}&=\frac{\frac{2}{M_1^{3/2}k_1}\Big((\frac{\sin^2\psi}{M_1}+\frac{\cos^2\psi}{4M_0})^{3/2}-\frac{|\sin^3\psi|}{M_1^{-3/2}}\Big)}{r^2M_1^{-3/2}\sin^2\psi(M_1^{-1}\sin^2\psi+\frac{\cos^2\psi}{4M_0})^{3/2}}-\frac{\frac{3}{2M_1^2k_1^2M_0}\cos^2\psi\sin^2\psi(v|\sin\psi|+sign(\psi)w\cos\psi)^2}{r^2(M_1^{-1}\sin^2\psi+\frac{\cos^2\psi}{4M_0})^{5/2}}.
\end{aligned}\]

\end{Lm}
\subsection{The second order derivatives  of $U_2$ with respect to $G_i$ and $g_i$, $i=1,2$}\label{AP-B2}
When evaluating at the isosceles limit, we have 
\begin{equation}
\label{U2G1s}\begin{aligned}\frac{\partial^2U_2}{\partial G_1\partial g_1}&=-\frac{2m_1m_2}{m_1+2}\frac{1}{L_1}(\frac{r_1r_2}{(r_2+\frac{2}{m_1+2}r_1)^3}-\frac{r_1r_2}{((r_2-\frac{m_1r_1}{m_1+2})^2+\frac{r_0^2}{4})^{3/2}}),\\
\frac{\partial^2U_2}{\partial g_1^2}&=-\frac{2m_1m_2}{m_1+2}(\frac{r_1r_2}{(r_2+\frac{2}{m_1+2}r_1)^3}-\frac{r_1r_2}{((r_2-\frac{m_1r_1}{m_1+2})^2+\frac{r_0^2}{4})^{3/2}}),\\
\frac{\partial^2U_2}{\partial G_1^2}&=-\frac{2m_1m_2}{m_2+2}\frac{1}{L_1^2}(\frac{r_1r_2}{(r_2+\frac{2}{m_1+2}r_1)^3}-\frac{r_1r_2}{((r_2-\frac{m_1r_1}{m_1+2})^2+\frac{r_0^2}{4})^{3/2}}).
\end{aligned}
\end{equation}
and (we only indicate the leading term here)
\begin{equation}\label{u2-estimates}
\begin{aligned}
&\frac{\partial^2 U_2}{\partial G_2^2}\sim\frac{\langle\partial_{G_2}^2\bx_2,\bx_2\rangle+|\partial_{G_2}\bx_2|^2}{|\bx_2|^3}\simeq\frac{1}{L_2^2r_2}\simeq\frac{\beta^2}{\chi}, \quad
\frac{\partial^2 U_2}{\partial g_2\partial G_2}\sim\frac{\langle\partial^2_{G_2g_2}\bx_2,\frac{2}{m_2}\bx_1\rangle}{|\bx_2+\frac{2\bx_1}{m_2+2}|^3}\simeq\frac1{L_2r_2}\simeq\frac{\beta}{\chi}.\\
&\frac{\partial^2 U_2}{\partial g_2^2}\sim\frac{\langle \partial_{g_2}^2\bx_2,\frac{2}{m_1+2}\bx_1\rangle}{|\bx_2+\frac{2\bx_1}{m_1+2}|^{3}}\simeq\frac{1}{r_2}\simeq\frac{1}{\chi},\quad
\frac{\partial^2 U_2}{\partial G_1^2}\sim \frac{\langle\partial_{G_1}^2\bx_1,\bx_2+\bx_1\rangle+|\partial_{G_1}\bx_1|^2}{|\bx_1+\bx_2|^3}\sim\frac{r_1r_2+r_1^2}{(r_1+r_2)^3},\\
&\frac{\partial^2U_2}{\partial G_1\partial g_1}\sim\frac{\langle\partial_{g_1G_1}\bx_1,\bx_2+\bx_1\rangle+\langle \partial_{g_1}\bx_1,\partial_{G_1}\bx_1\rangle}{|\bx_1+\bx_2|^3}\sim\frac{r_1r_2+r_1^2}{(r_1+r_2)^3},\\
 &\frac{\partial^2 U_2}{\partial g_1^2}\sim \frac{\langle\partial_{g_1}^2\bx_1,\bx_2+\bx_1\rangle+|\partial_{g_1}\bx_1|^2}{|\bx_1+\bx_2|^3}\sim\frac{r_1r_2+r_1^2}{(r_1+r_2)^3},\\
&\frac{\partial^2U_2}{\partial G_1\partial G_2}\sim\frac{\partial^2U_2}{\partial G_2\partial g_1}\simeq\frac{r_2r_1}{(r_1+r_2)^3L_2},\quad\frac{\partial^2U_2}{\partial G_1\partial g_2}\sim\frac{\partial^2U_2}{\partial g_2\partial g_1}\simeq\frac{r_2r_1}{(r_1+r_2)^3}\sim\frac{1}{\chi}.
\end{aligned}
\end{equation}
\subsection{The mixed second order derivatives of $U=U_{01}+U_2$ with respect to $L_i,\ell_i$ and $G_i,g_i$, $i=0,1,2$}\label{mix-derivatives} Let us assume $r_2=|\bx_2|=O(\chi)\gg1$, $r_1=|\bx_1|\gg1$, $|\bx_0|\sim1$ and 
\[|G_0|,|G_1|,|G_2|,|g_0-g_1|,|g_0-g_2|,|g_1-g_2|\ll1.\]
Let $\Upsilon=(|g_0-g_1-\frac{\pi}{2}|+|g_0-g_2+\frac{\pi}{2}|+|g_1-g_2-\pi|+\sum_{i=0}^2|G_{i}|)$. % Since in Delauney coordinates we have
%$|\bx_2|=\frac{1}{M_2k_2}L_2^2(e_2\cosh u_2-1),\quad |\bx_1|=\frac{1}{M_1k_1}L_1^2(e_1\cosh u_1-1),$ and
%\[\begin{split}\langle \bx_2,\bx_1\rangle=&\frac{1}{M_1k_1M_2k_2}\Big(\big(L_1^2(\cosh u_1-e_1)L_2^2(\cosh u_2-e_2)+L_1G_1\sinh u_1 L_2G_2\sinh u_2\big)\cos(g_1-g_2)\\
%&+\big(-L_2^2(\cosh u_2-e_2)L_1G_1\sinh u_1+L_1^2(\cosh u_1-e_1)L_2G_2\sinh u_2\big)\sin(g_1-g_2)\Big),
%\end{split}\]
%then 
 Direct calculation gives us the following estimates.
 \[\begin{aligned}\frac{\partial U}{\partial g_0}&=(\frac{1}{r_1^3}+\frac{1}{\chi^3})O( \Upsilon),\quad 
\frac{\partial U}{\partial G_0}=(\frac{1}{r_1^3}+\frac{1}{\chi^3} )O(\Upsilon),\quad
\frac{\partial U}{\partial g_1}=(\frac{1}{r_1^3}+\frac{r_1r_1}{|r_2-\frac{2}{m_1+2}r_1|^3})O(\Upsilon),\\
\frac{\partial U}{\partial G_1}&=(\frac{1}{r_1^3}+\frac{r_2r_1}{|r_2-\frac{m_1}{m_1+2}r_1|^3})O(\Upsilon),\quad
\frac{\partial U}{\partial g_2}=\frac{r_1r_2}{|r_1+r_2|^3}O(\Upsilon),\quad
\frac{\partial U}{\partial G_2}=\frac{\beta}{\chi}O(\Upsilon),
\end{aligned}
\]
and
\[\begin{aligned}
\partial_{g_1L_1}U_2&= (\frac{r_1r_2}{|r_1+r_2|^3}+\frac{r^2_2r_1^2}{|r_1+r_2|^5})O(\Upsilon),
\quad \partial_{g_1\ell_1}U_2=(\frac{r_2}{|r_1+r_2|^3}+\frac{r_2^2r_1}{|r_1+r_2|^5})O(\Upsilon),\\
\partial_{G_1 L_1}U_2&=( \frac{r_1r_2}{|r_1+r_2|^3}+\frac{r_1^2r_2^2}{|r_1+r_2|^5})O(\Upsilon),\quad
\partial_{G_1\ell_1}U_2=(\frac{r_2}{|r_1+r_2|^3}+\frac{r_1r_2^2}{|r_1+r_2|^5})O(\Upsilon),\\
\partial_{g_2L_2}U_2&=(\frac{\beta r_1r_2}{|r_1+r_2|^3}+\frac{\beta r_1^2r_2^2}{|r_1+r_2|^5})O(\Upsilon),\quad 
\partial_{g_2\ell_2}U_2=(\frac{\beta^{-2} r_1}{|r_1+r_2|^3}+\frac{\beta^{-2} r_1^2r_2}{|r_1+r_2|^5})O(\Upsilon),\\
\partial_{G_2 L_2}U_2&=O(\frac{\beta^3}{r_2}G_2)+(\frac{\beta^3 r_1r_2}{|r_1+r_2|^3}+\frac{\beta^3r_1^2r_2^2}{|r_1+r_2|^5})O(\Upsilon),\quad
\partial_{G_2\ell_2}U_2=\frac{1}{r_2^2}G_2+(\frac{r_1}{|r_1+r_2|^3}+\frac{r_1^2r_2}{|r_1+r_2|^5})O(\Upsilon),\\
\partial_{G_2 L_1}U_2&=(\frac{\beta^2 r_1r_2}{|r_1+r_2|^3}+\frac{\beta^2r_1^2r_2^2}{|r_1+r_2|^5})O(\Upsilon),\quad
\partial_{G_2\ell_1}U_2=(\frac{\beta^2r_2}{|r_1+r_2|^3}+\frac{\beta^2r_1r_2^2}{|r_1+r_2|^5})O(\Upsilon),\\
\end{aligned}\]\[
\begin{aligned}
\partial_{g_2\ell_1}U_2 &=(\frac{r_2}{|r_1+r_2|^3}+\frac{r_1r_2^2}{|r_1+r_2|^5})O(\Upsilon),\quad
\partial_{g_2 L_1}U_2=(\frac{ r_1r_2}{|r_1+r_2|^3}+\frac{r_1^2r_2^2}{|r_1+r_2|^5})O(\Upsilon),\\
\partial_{G_1L_2}U_2&=(\frac{ \beta r_1r_2}{|r_1+r_2|^3}+\frac{\beta r_1^2r_2^2}{|r_1+r_2|^5}O(\Upsilon),\quad
\partial_{G_1\ell_2}U_2=( \frac{ r_1}{|r_1+r_2|^3}+\frac{r_1^2r_2}{|r_1+r_2|^5})O(\Upsilon)\\
\partial_{g_1 L_2}U_2&=\frac{ (\beta r_1r_2}{|r_1+r_2|^3}+\frac{\beta r_1^2r_2^2}{|r_1+r_2|^5})O(\Upsilon),\quad
\partial_{g_1\ell_2}U_2= (\frac{ r_1}{|r_1+r_2|^3}+\frac{r_1^2r_2}{|r_1+r_2|^5})O(\Upsilon),
\end{aligned}
\]
and
\[\frac{\partial^2 U_{01}}{\partial \eta\partial\xi}=\frac{1}{r_1^3}O(\Upsilon),\quad \eta\in\{G_i,g_i,i=0,1\},\;\xi\in\{L_0,\ell_0\}.\]
\section{The derivatives of  the  Delaunay variables in the left-right transition  }\label{appendix-D}
When making transition from the left Jacobi coordinates to the right Jacobi coordinates on the middle section $\cS^M=\{-\beta\bx_1^L=\bx_2^R\}$, we have the following statement for the differentials.
\begin{Lm}\label{LmdL}
 On the middle section $\cS^M$, we have
\[\left[\begin{array}{c}
dL_1^R\\
dL_2^R
\end{array}\right]=
\left[\begin{array}{cc}
\frac{M_1^Lk^L_1}{M_1^Rk_1^R}\frac{(L_1^R)^2}{(L_1^L)^2}&\frac{M^L_1M_2^Lk^L_2}{M_2^Rk_2^R} \frac{(L_1^R)^2}{2(L_2^L)^2}\\
-\frac{M^L_1M_2^Lk^L_2}{M_1^Rk_1^R}  \frac{(L_2^R)^2}{2(L_1^L)^2}&(1-\frac{M_1^LM_2^R}{4})\frac{M^L_2k_2^L}{M^R_1k^R_1}\frac{(L_2^R)^2}{(L_2^L)^2}
\end{array}\right]\left[\begin{array}{c}
dL_1^L\\
dL_2^L
\end{array}\right]\]
\[2L_1^R\ell_1^RdL_1^R+(L_1^R)^{2}d\ell_1^R=(1-\frac{M_1^LM_2^R}{4})\frac{M_1^Rk_1^R}{M^L_1k^L_1}\left(2L_1^L\ell_1^LdL_1^L+(L_1^L)^{2}d\ell_1^L\right)- \frac{M_1^Rk_1^R}{M^L_2k^L_2}\left(2L_2^L\ell_2^LdL_2^L-(L_2^L)^{2}d\ell_2^L\right),\]
\[
2L_2^R\ell_2^RdL_2^R+(L_2^R)^{2}d\ell_2^R= \frac{M_2^Rk_2^R}{2k^L_1}\left(2L_1^L\ell_1^LdL_1^L+(L_1^L)^{2}d\ell_1^L\right)+\frac{M_2^Rk_2^R}{M_2^Lk_2^L}\left( 2L_2^L\ell_2^LdL_2^L+(L_2^L)^{2}d\ell_2^L\right).\\
\]

\end{Lm}

\begin{proof}[Proof of Lemma \ref{LmdL}] In the isosceles case, we have $\bx_1^L=(x_1^L,0)$, $\bp_1^L=(p_1^L,0)$, $\bx^L_2=(x_2^L,0)$ and $\bx^L_2=(p_2^L,0)$. The same formulas hold true for the right Jacobi coordinates.
In terms of Delaunay coordinates (see Appendix \ref{App-delaunay}), we have
\begin{equation}
\begin{aligned}
x_1^L&=-\left[\frac{(L_1)^2}{M_1k_1}(\cosh u_1-1)\right]^L<0,\quad p_1^L=\left[\frac{M_1k_1}{L_1}\frac{\sinh u_1}{1-\cosh u_1}\right]^L<0,\\
x_2^L&=\left[\frac{(L_2)^2}{M_2k_2}(\cosh u_2 -1)\right]^L>0,\quad p_2^L=\left[\frac{M_2k_2}{L_2}\frac{\sinh u_2}{1-\cosh u_2}\right]^L>0,\\
x_1^R&=-\left[\frac{(L_1)^2}{M_2k_2}(\cosh u_1-1)\right]^R<0,\quad p_1^R=\left[\frac{M_2k_2}{L_1}\frac{\sinh u_1}{1-\cosh u_1}\right]^R<0,\\
x_2^R&=\left[\frac{(L_2)^2}{M_1k_1}(\cosh u_2 -1)\right]^R>0,\quad p_2^R=\left[\frac{M_1k_1}{L_2}\frac{\sinh u_2}{1-\cosh u_2}\right]^R>0.\\
\end{aligned}
\end{equation}
We remark that there are two subtleties here. First, since $\dot \ell<0$ decreases, $u$ and $\ell$ have different signs if $|\ell|$ is large, so $u$ changes from negative to positive along the orbit. We choose $u_1>0$ and $u_2<0$. When we change coordinates from $L$ to $R$, we change $M_1^L,k_1^L$ to $M_2^R,k_2^R$ and vice versa.
%\begin{equation}
%\begin{aligned}
%-\frac{(L_1^R)^{2}}{m_2k_2}(\cosh u_1^R-1)&=-\frac{4\mu(1+\mu)}{(1+2\mu)^2}\frac{(L_1^L)^{2}}{m_1k_1}(\cosh u_1^L-1)-\frac{1}{1+2\mu}\frac{(L_2^L)^{2}}{m_2k_2}(\cosh u_2^L-1),\\
%\frac{(L_2^R)^{2}}{m_1k_1}(\cosh u_2^R-1)&=-\frac{1}{1+2\mu}\frac{(L_1^L)^{2}}{m_1k_1}(\cosh u_1^L-1)+ \frac{(L_2^L)^{2}}{m_2k_2}(\cosh u_2^L-1),\\
%\frac{m_2k_2}{L^R_1}\frac{\sinh u_1^R}{1-\cosh u_1^R}&=\frac{m_1k_1}{L^L_1}\frac{\sinh u_1^L}{1-\cosh u_1^L}-\frac{1}{1+2\mu}\frac{m_2k_2}{L^L_2}\frac{\sinh u_2^L}{1-\cosh u_2^L},\\
%\frac{m_1k_1}{L^R_2}\frac{\sinh u_2^R}{1-\cosh u_2^R}&=\frac{1}{1+2\mu}\frac{m_1k_1}{L^L_1}\frac{\sinh u_1^L}{1-\cosh u_1^L}+\frac{4\mu(1+\mu)}{(1+2\mu)^2}\frac{m_2k_2}{L^L_2}\frac{\sinh u_2^L}{1-\cosh u_2^L}.
%\end{aligned}
%\end{equation}

Using Lemma \ref{Lm: simplify} in the Appendix, we can simplify formula \eqref{EqLR}  as
\begin{equation}
\begin{aligned}
\frac{(L_1^R)^{2}\ell_1^R}{M_1^Rk_1^R}&=(1-\frac{M_1^LM_2^R}{4})\frac{(L_1^L)^{2}}{M_1^Lk_1^L}\ell_1^L- \frac{(L_2^L)^{2}}{M_2^Lk_2^L}\ell_2^L,\quad
\frac{(L_2^R)^{2}}{M_2^Rk_2^R}\ell_2^R= \frac{(L_1^L)^{2}}{2k^L_1}\ell_1^L+ \frac{(L_2^L)^{2}}{M^L_2k^L_2}\ell_2^L,\\
-\frac{M^R_1k^R_1}{L^R_1} &=-\frac{M^L_1k^L_1}{L^L_1} - \frac{M^L_1M^L_2k^L_2}{2L^L_2} ,\quad
\frac{M^R_2k^R_2}{L^R_2} =- \frac{M^L_1M^L_2k^L_1}{2L^L_1} +(1-\frac{M_1^LM_2^R}{4})\frac{M_2^Lk_2^L}{L^L_2} .
\end{aligned}
\end{equation}
The  lemma is then proved by differentiating the above expressions.
\end{proof}

We now prove Lemma \ref{LmTransitionGg}.
\begin{proof}[Proof of Lemma \ref{LmTransitionGg}]
We use formula \eqref{EqLR}.  We differentiate the expressions $\bx_1^R,\bx_2^R$, then we evaluate them in the isosceles limit. By Lemma \ref{LmD2G}, we have that  that in the isosceles limit,
\begin{equation}\label{llggbot}
\{\frac{\partial \bx_i}{\partial G_i},\;\frac{\partial \bx_i}{\partial g_i},1,2\}\bot\{\frac{\partial \bx_i}{\partial L_i},\;\frac{\partial \bx_i}{\partial \ell_i},i=1,2\}.
\end{equation}Therefore,  we obtain the following equality,
\[
\begin{aligned}
\frac{\partial \bx_1^R}{\partial g_1^R}dg_1^R+\frac{\partial \bx_1^R}{\partial G_1^R}dG_1^R &=(1-\frac{1}{4}M_1^LM_2^R )\left(\frac{\partial \bx_1^L}{\partial g_1^L}dg_1^L+\frac{\partial \bx_1^L}{\partial G_1^L}dG_1^L\right)-\frac{M_2^R}{2}\left(\frac{\partial \bx_2^L}{\partial g_2^L}dg_2^L+\frac{\partial \bx_2^L}{\partial G_2^L}dG_2^L\right),\\
\frac{\partial \bx_2^R}{\partial g_2^R}dg_2^R+\frac{\partial \bx_2^R}{\partial G_2^R}dG_2^R&=\frac{M_1^L}{2}\left(\frac{\partial \bx_1^L}{\partial g_1^L}dg_1^L+\frac{\partial \bx_1^L}{\partial G_1^L}dG_1^L\right)+ \frac{\partial \bx_2^L}{\partial g_2^L}dg_2^L+\frac{\partial \bx_2^L}{\partial G_2^L}dG_2^L.\\
\end{aligned}
\]
From the formula in Lemma \ref{LmD2G},  we have
$$\frac{\partial \bx^R_i}{\partial G^R_i}=-\frac{L_i}{M_ik_i}(0,\sinh u_i),
\ \frac{\partial \bx^R_i}{\partial g_i^R}=-\frac{L_i^2}{M_ik_i}(0,\cosh u_i-1)),$$
and on  $\cS^M$ we have the following estimates for the coefficients,$$\frac{\partial \bx_1^R}{\partial g_1^R}\sim \frac{\partial \bx_1^R}{\partial G_1^R}\sim \frac{\partial \bx_1^L}{\partial g_1^L}\sim\frac{\partial \bx_1^L}{\partial G_1^L} \sim  \frac{\partial \bx_2^L}{\partial g_2^L}\sim\chi, \quad \frac{\partial \bx_2^R}{\partial g_2^R}\sim  \frac{\partial \bx_2^R}{\partial G_2^R}\sim\frac{\partial \bx_2^L}{\partial G_2^L}\sim\beta\chi.$$
This gives
$$dg_1^R=\bar c_1 dG_1^R+O(1)d g_1^L+O(1)dG_1^L+O(1)dg_2^L+O(\beta)dG_2^L,\quad $$
$$dg_2^R=\bar c_2 dG_2^R+O(\frac{1}\beta)d g_1^L+O(\frac{1}\beta)dG_1^L+O(\frac{1}\beta)dg_2^L+O(1)dG_2^L,\quad $$
where we have $$\bar c_1=-\frac{1}{L^R_1}\frac{\sinh u_1}{\cosh u_1-1}\text{ and }\bar c_2=-\frac{1}{L^R_2}\frac{\sinh u_2}{\cosh u_2-1}.$$
We next show that both $dG_1^R$ and $dG_2^R$ are of order $O(\chi)$, so they give the leading contribution.

For the quantity  $G_1^R$, using again \eqref{EqLR}, we have
\begin{equation}\label{G1R}
\begin{aligned}
G_1^R&=\bp_1^R\times \bx_1^R=\left(\bp_1^L-\frac{M^L_1}{2}\bp_2^L\right)\times \left((1-\frac{1}{4}M_1^LM_2^R )\bx_1^L-\frac{M^L_2}{2}\bx_2^L\right)\\
&=(1-\frac{M_1^LM_2^R}{4} ) G_1^L+\frac{M_1^LM_1^R}{4}G_2^L-(1-\frac{1}{4}M_1^LM_2^R ) \frac{M_1^L}{2}\bp_2^L\times \bx_1^L-\frac{M_2^L}{2} \bp_1^L\times \bx_2^L,
%\\
%G_2^R&= p_2^R\times x_2^R=\left(\frac{M_2^R}{2}p_1^L+(1-\frac{1}{4}M_1^LM_2^R )p_2^L\right)\times\left(\frac{M_1^L}{2}x_1^L+ x_2^L\right)\\
%&= \frac{M_1^LM_2^R}{4}G_1^L+(1-\frac{1}{4}M_1^LM_2^R )G_2^L+(1-\frac{1}{4}M_1^LM_2^R ) \frac{M_1^L}{2}p_2^L\times x_1^L+\frac{M_2^R}{2} p_1^L\times x_2^L
\end{aligned}
\end{equation}
where the operator $\times$ is the standard outer product on $\mathbb{R}^2$.
We then differentiate the $G_1^R$ expression and evaluate it in the isosceles limit. Note that in the isosceles limit, 
$$\dfrac{\partial \bp}{\partial G}\perp \bp,\ \dfrac{\partial \bp}{\partial g}\perp \bp,\  \dfrac{\partial \bp}{\partial L}\parallel \bp,\ \dfrac{\partial \bp}{\partial \ell}(=O(1/\chi))\parallel p,\ \dfrac{\partial \bx}{\partial G}\perp \bx,\ \dfrac{\partial \bx}{\partial g}\perp \bx,\  \dfrac{\partial \bx}{\partial L}\parallel \bx,\ \dfrac{\partial \bx}{\partial \ell}(=O(1))\parallel \bx.$$
Combining this with \eqref{llggbot}, we have
\begin{equation}\label{DG1R}
\begin{aligned}
dG_1^R&=(1-\frac{M_1^LM_2^R}{4} )dG_1^L+\frac{M_1^LM_1^R}{4}dG_2^L\\
&\quad-(1-\frac{M_1^LM_2^R}{4} ) \frac{M_1^L}{2}\left[\left(\frac{\partial \bp_2^L}{\partial G_2^L} dG_2^L+\frac{\partial \bp_2^L}{\partial g_2^L} dg_2^L\right)\times \bx_1^L+ \bp_2^L\times\left(\frac{\partial \bx_1^L}{\partial G_1^L} dG_1^L+\frac{\partial \bx_1^L}{\partial g_1^L} dg_1^L\right)\right]\\
&\quad-\frac{M_2^L}{2} \left[\left(\frac{\partial \bp_1^L}{\partial G_1^L} dG_1^L+\frac{\partial \bp_1^L}{\partial g_1^L} dg_1^L\right)\times \bx_2^L+ \bp_1^L\times\left(\frac{\partial \bx_2^L}{\partial G_2^L} dG_2^L+\frac{\partial \bx_2^L}{\partial g_2^L} dg_2^L\right)\right]\\
&=O(\chi,\chi,\beta\chi,\chi),
\end{aligned}
\end{equation}
Here the leading $O(\chi)$ terms are give by $\frac{\partial \bp_1^L}{\partial G_1^L}\times \bx_2^L dG_1^L+\frac{\partial \bp_1^L}{\partial g_1^L}\times \bx_2^L dg_1^L+\bp_1^L\times\frac{\partial \bx_2^L}{\partial g_2^L} dg_2^L$.
By Lemma \ref{LmD2G}, we get
$$\bx_2^L=\left(-\frac{(L_2^L)^2}{M_2^L k_2^L} (\cosh u_1^L-1),0\right),\ \bp_1^L=\left(\frac{M_1^Lk_1^L}{L_1^L}\frac{\sinh u_1^L}{\cosh u_1^L-1},0\right),\ \frac{\partial \bp_1^L}{\partial G_1^L}=\left(0,\frac{M_1^L k_1^L}{(L_1^L)^2}\frac{\cosh u_1^L}{\cosh u_1^L-1}\right)$$
Note that $\frac{\partial \bp_1^L}{\partial g_1^L}\times \bx_2^L $ and $\bp_1^L\times\frac{\partial \bx_2^L}{\partial g_2^L} $ have opposite signs, and $\frac{\partial \bp_1^L}{\partial G_1^L}$ and $\frac{\partial \bp_1^L}{\partial g_1^L}$ point in the same direction.
So we get $$\nabla_{Y_L}G_1^R=  \left(\frac{1}{L_1^L},1,\beta,-1\right)|\bx_2^L|+O(1),$$
where $Y_L=(G_1^L,g_1^L,G_2^L,g_2^L)$.

Now we consider the quantity $G_2^R$. Recall that $G_0$ remains untouched in the transition, so 
by the total angular momentum conservationt $G_1^R+G_2^R+G_0=G_1^L+G_2^L+G_0$ we get  $dG_2^R=-dG_1^R+dG_1^L+dG_2^L. $
We then have  $\nabla_{Y_L}G_2^R=-\nabla_{Y_L}G_1^R+O(1)$. This finishes the proof of Lemma \ref{LmTransitionGg}.
\end{proof}
\begin{proof}[Proof of \eqref{non-iso-transit}]
For the non-isosceles cases with $|G_1^L|, |G_2^L|, |g_1^L-g_2^L-\pi|<\nu$ and $|G_1^R|, |G_2^R|, |g_1^R-g_R^L-\pi|<\nu_0$, from \eqref{G1R}, we know that for  $dG_1^R$, there are additional terms, besides those in \eqref{DG1R}, namely
\[-(1-\frac{M_1^LM_2^R}{4})\frac{M_1^L}{2}\Big[(\frac{\partial \bp_2^L}{\partial L_2^L}dL_2^L+\frac{\partial \bp_2^L}{\partial \ell_2^L}d\ell_2)\times \bx_1^L+\bp_2^L\times(\frac{\partial \bx_1^L}{\partial L_1^L}dL_1^L+\frac{\partial \bx_1^L}{\partial\ell_1^L}d\ell_1^L)\Big],\]
and 
\[-\frac{M_2^L}{2}\Big[(\frac{\partial \bp_1^L}{\partial L_1^L}dL_1^L+\frac{\partial \bp_1^L}{\partial \ell_1^L}d\ell_1)\times \bx_2^L+\bp_1^L\times(\frac{\partial \bx_2^L}{\partial L_2^L}dL_2^L+\frac{\partial \bx_2^L}{\partial\ell_2^L}d\ell_2^L)\Big].\]
This leads to the term $O(\nu\chi)dL_1^L+O(\nu)d\ell_1+O(\nu\beta\chi)dL_2^L+O(\beta^{-2}\nu)d\ell_2$. 
For the quantities $L_1^R$ and $\ell_1^R$, from \eqref{eq: delaunay4} and Lemma \ref{Lm: simplify}, we know that there will be additional terms $\nu_0dG_1^R+O(\nu\chi)(dG_1^L+dG_2^L+dg_1^L+dg_2^L)=O(\nu_0\chi)(dG_1^L+dG_2^L+dg_1^L+dg_2^L)$. Then we have the asserted form of the transition matrix \eqref{non-iso-transit}. 
\end{proof}

\section{Numeric verification of the transversality conditions}\label{STrans}

Since the dynamics near triple collision is  nonperturbative, in particular, it is not easy to explicitly find $\gamma_I$ and $\gamma_O$ and integrate the variational equations along them, we choose to verify the transversality condition using a computer, as people have done in the literature (see \cite{S,SM} etc). We only consider the case with 4 equal masses, that is, $m_1=m_2=m_3=m_4=1$, so that both the equations of motion and the variational equations are concrete ODEs. The numeric task is to solve for $\gamma_I$ and $\gamma_O$, then integrate along them the variational equation in Lemma~\ref{LmFundLoc} with the stable and unstable directions at the fixed point as initial conditions. 

Since the masses are equal, the triple collision systems on the  left and right sides are the same. We just  need to check the condition for one side, say the left side. We pick the Lagrangian fixed point $O_1$ with $v<0$ and $\psi<0$. 

{\bf Step 1,} {\it the piece of orbit $\gamma$ shadowing $\gamma_{O}$ and renormalization. }

We need to numerically compute\footnote{We used the ODE solver 'ode113' in Matlab with 'RelTol'=$10^{-13}$ and 'AbsTol'=$10^{-15}$, and with the initial condition $O_1+10^{-14}\mathbf{e}^u$, where $\mathbf{e}^u$ is the unstable direction of the fixed point $O_1$ on the collision manifold.  }  $\gamma^1_{O}$ (going to the side $\psi=-\frac{\pi}{2}$) by solving the regularized system  \eqref{EqSundman} on the collision manifold  from  a very small neighborhood of a fixed point along the unstable direction to the section $\cS^L_+$, and integrate\footnote{We used the ODE solver 'ode113' in Matlab with 'RelTol'=$10^{-7}$ and 'AbsTol'=$10^{-8}$.} the variational equation in Lemma~\ref{LmFundLoc}  along this orbit for the variables $(dw_0,d\mathsf{g}_{01},dw_2,d\mathsf{g}_{12})$ with  initial conditions   $\mathbf{e}_1(O_1)$ and $\mathbf{e}_2(O_1)$, which are the unit eigenvectors corresponding to the two largest eigenvalues $\mu_u>0$ and $-\frac{v_*}{2}$, starting from a neighborhood of $O_1$. We perform the renormalization on the section $\cS_+^L$ and then change from the coordinate $(dw_0,d\mathsf{g}_{01},dw_2,\mathsf{g}_{12})$ to $(dG_0,d\mathsf{g}_{01},dG_2,d\mathsf{g}_{12})$. Numerically, we get $E^u_{4}(\gamma(0))=\text{span}\{\hat{\mathbf{e}}_1,\hat{\mathbf{e}}_2\}$, where $\hat{\mathbf{e}}_1=(0.5847,0.5608,0,0.5862)$ and  $\hat{\mathbf{e}}_2=( 0.4774,0.5282, -0.0701, 0.6987)$.

% {\bf Step 2:} {\it the piece of orbit shadowing $\gamma_I$. }

%  The transversality amounts to show that the vector $\hat{\be}_{0}$, when pushed forward along $\gamma^2_I$ approaching the Lagrange fixed point, does not coincides with the stable $E_0^s$ direction at the fixed point $O_2$ in the $(w_0,g_0)$-plane. To prove this, we first integrate the variational equation in Lemma~\ref{LmVar0} with reversed time from the fixed point to the section $\cS_-^R=\{r^R\simeq10^4\}$. The fundamental solution gives rise to a linear transform, sending the vector $E_0^s$ at the fixed point to a vector $\hat \be_1$ at $\cS^R_-$. Then on $\cS_-^R$, we change coordinates from $(w_0,g_0)$ to $(G_0,g_0)$. The derivative map of this change of coordinates on the section $\cS^R_-$ results in  multiplying the $\delta w_0$-component by $(r^R)^{-1/2}$. Then we get a vector $\hat{\mathbf{e}}_2$ in the $(G_0,g_0)$-coordinates.

 % {\bf The numeric result:}  
 
Next, we find the orbit $\gamma^1_I$ numerically. Note that for the I3BP, at the Lagrange fixed point $O_1$ with $v_*<0$, there are two linearly independent stable directions. One is on the collision manifold, which we denote as~$\mathbf{e}_1^s$, the other is $\mathbf{e}_2^s=(v_*,E,0,0)$, where $E$ is the total energy\footnote{We take the total energy to be $-8/9$.}. So we  solve\footnote{We use the ODE solver `ode113' in Matlab with 'RelTol'=$10^{-13}$ and 'AbsTol'=$10^{-14}$.}  the regularized system \eqref{EqSundman} backward  along the stable direction with initial conditions of the form $$O_1+10^{-6}\big(-\cos(k10^{-3}\frac{\pi}{2})\mathbf{e}_2^s+\sin(k10^{-3}\frac{\pi}{2})\mathbf{e}_1^s\big), \quad k\in[-1000,1000]\cap\mathbb{Z}.$$ From the numeric result, there exists one initial condition ($k=-58$) for which the corresponding backward solution satisfies the energy partition condition (with error $10^{-3}$) for $\gamma^1_I$ when the variable $r$ is of order $10^4$ and goes to the side $\psi=-\frac{\pi}{2}$. This is the numerical realization of the argument in the proof of Proposition \ref{ThmGammaI}. 
 We then  evaluate the variational equation in in Lemma \ref{LmFundLoc}   backward along the corresponding orbit  with $k=-58$,  with initial condition  $E_0^s(O_1)$ and then change from the coordinate $(dw_0,d\mathsf{g}_{01},dw_2,\mathsf{g}_{12})$ to $(dG_0,d\mathsf{g}_{01},dG_2,d\mathsf{g}_{12})$ on the section $\cS^R_-$. Numerically, we get that $E_{4}^s(\gamma_I(0))$ is spanned by $\{\hat{\mathbf{e}}_3,\hat{\mathbf{e}}_4\}$, where  $\hat{\mathbf{e}}_3=(0,0,0,1)$ and  $\hat{\mathbf{e}}_4=( -0.8549, 0.3736, 0, 0.3599)$. We also have $E^s_3(\gamma_I(0))=\text{span}\{\bar{\mathbf{e}}_4:=( -0.9163,0.4004)\}$

  {\bf Step 3:} {\it Checking the cone conditions  in Proposition \ref{PropDG2}. } Since $\bar{\mathbf{l}}_2=(0,0,\frac{1}{L_1^L},1)+O(\beta)$, from the definition of $\bar{\mathbf{u}}_3$ and the numerics in Step 1, we have $\mathbf{u}=\frac{0.0701-0.6607L_1^L}{0.5862L_1^L}\hat{\mathbf{e}}_1+\hat{\mathbf{e}}_2+O(\beta)$ (see Section \ref{SSTrans} for the definition of ${\mathbf u}$) and setting the last two entries of $\mathbf{u}$ to zero we obtain $\bar{\mathbf{u}}_3$. Since $L_1^L\in(e^{-1},e)$\footnote{Choosing  a different renormalization factor could make this interval into an arbitrary small neighborhood containing 1}, we know that $\text{span}\{\bar{\mathbf{u}}_3\}\cap\mathcal{C}_{\frac{1}{100}}(E^s_{3}(\gamma_I(0)))=\{0\}$.
  By the relation \eqref{LRL1L2}, we have
\[c_1=\frac{-1}{L_1^R}=\frac{-4(m_1+m_2+2)}{(m_1+2)(m_2+2)(m_2+1)}\frac{1}{L_1^L}=-\frac{8}{9}\frac{1}{L_1^L},\quad c_2=\frac{1}{L_2^R}=\frac{m_1^2}{4m^2_2(m_1+2)}\frac{1}{L_1^L}=\frac{1}{12}\frac{1}{L_1^L}.\]
Since $\bar{\mathbf{u}}_2=(0,-c_2,1,c_1-c_2)$,  with the numeric results from Step  2, we clearly have $\text{span}\{\bar{\mathbf{u}}_2,\bar{\mathbf{u}}_3\}\cap\mathcal{C}_{\frac{1}{100}}(E^s_{4}(\gamma_I(0)))=\{0\}$. 

\begin{figure}[ht]
\begin{subfigure}{0.48\textwidth}
\includegraphics[height=4.5cm,width=7cm]{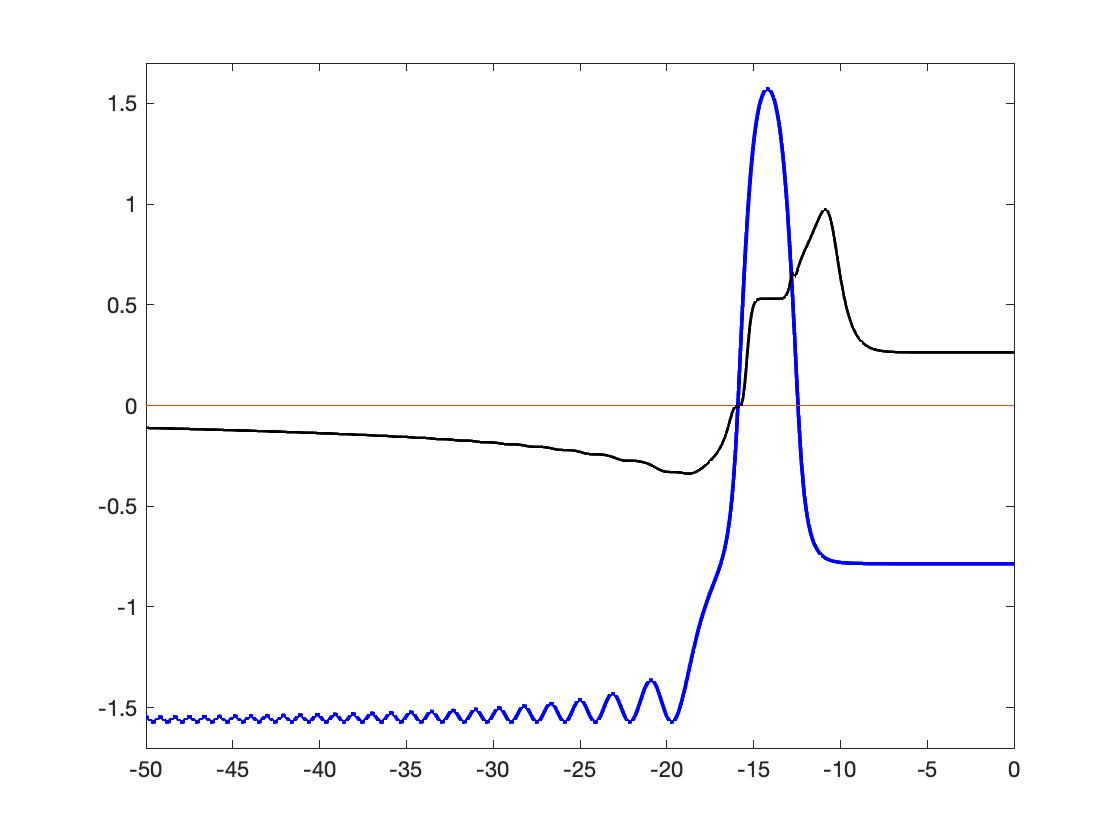}
\end{subfigure}
\begin{subfigure}{0.48\textwidth}
\includegraphics[height=4.5cm,width=7cm]{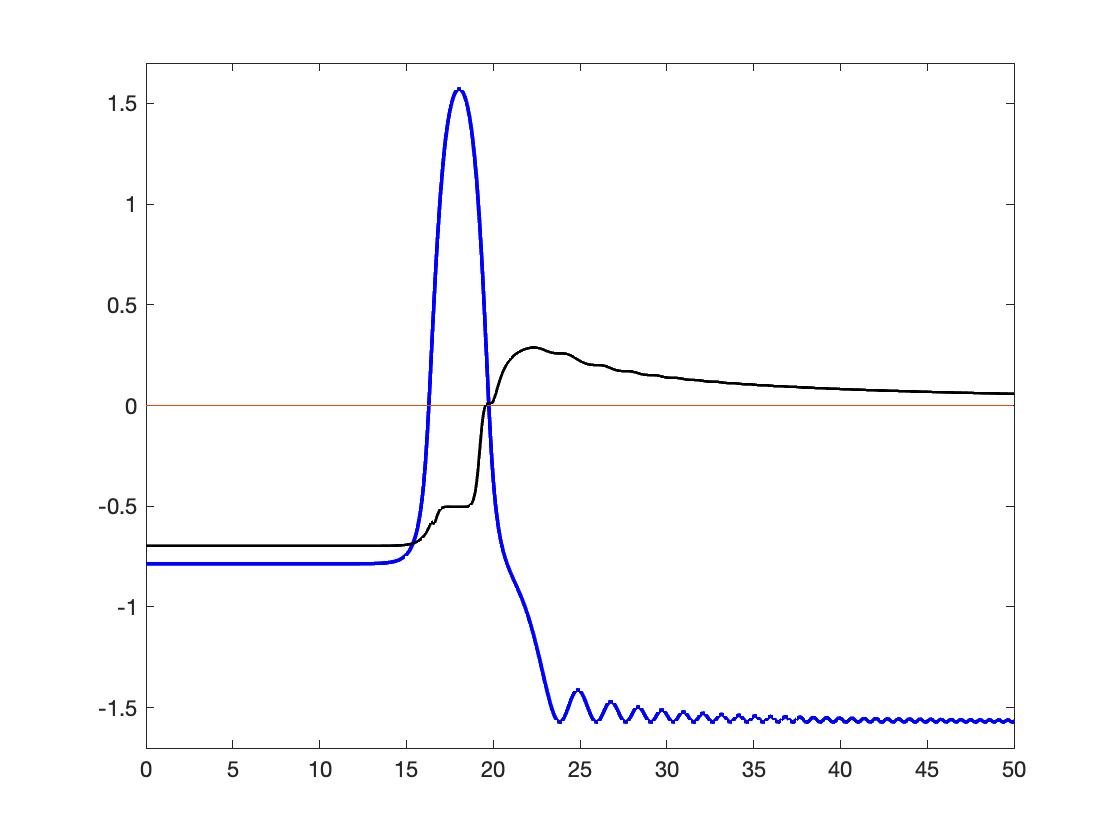}
\end{subfigure}
\caption{The left: backward solution along $\gamma_I^1$ and  the right: forward solution along $\gamma_O^1$; the bright curves are the evolution of the $\psi$-variable and the dark curves are that of the $dw_0$-variable. }
\label{fig-dw0}
\end{figure}

  {\bf Step 4:} {\it Verifying Proposition \ref{PropTrans1}. } We need to keep track of $dw_0$ components of $E^u_0(\gamma_O(t))$, $t<0$ and $E^s_0(\gamma_I(t))$, $t>0$. As we already find $\gamma_I$ and $\gamma_O$ in Steps 1 and 2, we just need to evaluate the variational equation \eqref{variation-w0g0} in Lemma \ref{LmVar0} along $\gamma_O$ forward in time,  with initial condition $E^u(O)$ and along $\gamma_I$ backward in time with initial condition $E^s(O)$, from a small neighborhood of the Lagrangian fixed point. The evolution of $dw_0$ and $\psi$ is illustrated in Figure \ref{fig-dw0}. The values are close to nonzero constants when the orbits are in a neighborhood of the Lagrange fixed point as well as close to the sections $\cS_\pm$, so they can only become zero at finitely many points. On the other hand, along $\gamma_I$ and $\gamma_O$, the variable $\psi$ also becomes zero at finitely many points, in particular it does not vanish in a neighborhood of the Lagrange fixed point. This reduces the problem to checking the non-coincidence of the zeros of $\psi\pm\pi/2$ and $w_0$ on two compact time intervals. Hence we have numerically justified the non-degeneracy  in Proposition  \ref{PropTrans1}.
\section{A Lemma on analytic functions of two variables}\label{SUlam}

\begin{Lm}\label{LmUlam} Let $f:\ E\subset \R^2\to \R$ be an analytic function of two variables on an open set $E$.  Suppose $f$ is nonvanishing at $(x,y)=(0,0)$, then there exists a residual set $\mathcal R\subset E$  such that for each $(x,y)\in \mathcal R$, we have $f(x,y)\neq 0. $
\end{Lm}
\begin{proof}
We apply the Kuratowski-Ulam theorem (c.f. Theorem 15.4 of \cite{Ox}: If $E\subset X\times Y$ has the Baire property, and if $E_x=E\cap \{x\}\times Y$ is of first category for all $x$ but a set of first category, then $E$ is of first category). By assumption, the function $f(x,0)$ is vanishing only on a set of first category denoted by $\hat X$, since $f(0,0)\neq 0$ and the analytic function $f(x,0)$ is either constant or nonconstant, and in the latter case its zero is isolated. For each $x\notin \hat X$, we consider the function $f(x,y)$ as a function of $y$, the same argument shows that there is a set $E_x$ of first category such that for each $y\notin E_x$, we have $f(x,y)\neq 0$. Thus, the statement follows from the Kuratowski-Ulam theorem. 
\end{proof}

\section*{Acknowledgement}
This paper grew out of an offhand remark by Saari to Gerver in 2012.  We also thank Professors Dolgopyat, Moeckel, Montgomery and G. Yu for their illuminating discussions and  kind suggestions. G. H and J.X. are supported by the grant NSFC (Significant project No.11790273) in China. J. X. is in addition supported by the Beijing Natural Science Foundation (Z180003).

\end{document}